\documentclass[11pt]{amsart}
\usepackage[left=0.9in,right=0.9in,top=0.75in,bottom=0.7in,includeheadfoot]{geometry}
\usepackage{enumerate}
\usepackage{comment}
\usepackage{mathpazo}
\usepackage{microtype}

\usepackage{bm}
\usepackage{tikz-cd}
\IfFileExists{stmaryrd.sty}{\usepackage{stmaryrd}}{}
\usepackage{amsfonts}
\usepackage{amssymb}
\usepackage{amsmath,amsthm}
\usepackage{latexsym}
\usepackage{amscd}
\usepackage{esint}
\usepackage{mathrsfs}
\usepackage{indentfirst}
\usepackage{enumitem}
\newlist{Aenumerate}{enumerate}{1}
\setlist[Aenumerate]{label=(S\arabic*)}

\usepackage{enumerate}
\usepackage{slashed}
\usepackage{fancyhdr}
\usepackage{aliascnt}
\usepackage{pdfsync}
\usepackage{setspace}
\usepackage{color}
\usepackage{xcolor}
\definecolor{citeblue}{RGB}{0,0,139}
 \usepackage[colorlinks = true,
            linkcolor = citeblue,
            urlcolor  =  citeblue,
            citecolor =  citeblue,
            anchorcolor = citeblue]{hyperref}
\newtheorem{theorem}{Theorem}[section]

\newtheorem{lemma}[theorem]{Lemma}

\theoremstyle{definition}

\theoremstyle{remark}
\newtheorem{rmk}[theorem]{Remark}
\numberwithin{equation}{subsection}
\numberwithin{theorem}{subsection}
\newcommand{\norm}[1]{\left\Vert#1\right\Vert}

\newcommand{\Vol}{\mathrm{Vol}}

\newcounter{stepnum}[section]

\newcounter{claimnum}[subsection]

\newcounter{casenum}[section]

\usepackage{etoolbox}

\makeatletter
\newcommand{\reduceoperator}[1]{%
  \@for\next:=#1\do{\expandafter\reduceoperator@\expandafter{\next}}%
}
\newcommand{\reduceoperator@}[1]{%
  \csletcs{normal@#1@}{#1@}%
  \csedef{#1@}{\noexpand\reduceoperator@@\csname normal@#1@\endcsname}%
}
\newcommand{\reduceoperator@@}[1]{%
  \mathop{\mathpalette\reduceoperator@@@{#1}}%
}
\newcommand{\reduceoperator@@@}[2]{%
  \ifx#1\displaystyle\textstyle\fi#2%
}
\makeatother

\reduceoperator{bigcup,bigcap,bigotimes}

\def\R{\mathbb R}
\def\S{\mathbb S}

\selectfont
\newcommand{\sgn}[1]{\text{sgn}\left(#1\right)}

\makeatletter
\def\@tocline#1#2#3#4#5#6#7{\relax
	\ifnum #1>\c@tocdepth %
	\else
	\par \addpenalty\@secpenalty\addvspace{#2}%
	\begingroup \hyphenpenalty\@M
	\@ifempty{#4}{%
		\@tempdima\csname r@tocindent\number#1\endcsname\relax
	}{%
		\@tempdima#4\relax
	}%
	\parindent\z@ \leftskip#3\relax \advance\leftskip\@tempdima\relax
	\rightskip\@pnumwidth plus4em \parfillskip-\@pnumwidth
	#5\leavevmode\hskip-\@tempdima
	\ifcase #1
	\or\or \hskip 2em \vspace{0.5mm}\or \hskip 3em \else \hskip 3em \fi%
	#6\nobreak\relax
	\hfill\hbox to\@pnumwidth{\@tocpagenum{#7}}\par%
	\nobreak
	\endgroup
	\fi}
\makeatother
\usepackage{mathtools}
\usepackage{bookmark}

\makeatletter
\@removefromreset{equation}{subsection}
\@removefromreset{theorem}{subsection}
\makeatother
\numberwithin{equation}{section}
\numberwithin{theorem}{section}

\theoremstyle{plain}
\newtheorem{proposition}[theorem]{Proposition}
\newtheorem{corollary}[theorem]{Corollary}
\theoremstyle{definition}

\theoremstyle{remark}
\newtheorem{remark}[theorem]{Remark}

\DeclareMathOperator{\Div}{div}
\DeclareMathOperator{\Hess}{Hess}
\DeclareMathOperator{\tr}{tr}
\DeclareMathOperator{\cof}{cof}
\DeclareMathOperator{\ind}{ind}
\DeclareMathOperator{\nul}{nul}

\DeclareMathOperator{\Id}{Id}
\let\sgn\relax
\DeclareMathOperator{\sgn}{sgn}
\DeclareMathOperator{\Span}{span}

\newcommand{\Sph}{\mathbb S}
\newcommand{\dd}{\,d}
\newcommand{\ip}[2]{\langle #1,#2\rangle}
\newcommand{\Aop}{\mathcal A}
\newcommand{\Killing}{\mathcal K}
\newcommand{\Hone}{\mathcal H_1}
\setlist[enumerate,1]{label=(\arabic*),leftmargin=2.2em,itemsep=3pt}
\hypersetup{
  pdftitle={Spectral and convex geometry of alpha-harmonic maps and sharp Morse index rigidity of spheres},
  pdfsubject={Research manuscript; H-disk LaTeX template adaptation},
  pdfkeywords={Sacks-Uhlenbeck functional, Morse index, rigidity, bubble tree, Codazzi tensor}
}
\begin{document}
\title[Spectral and convex geometry of $\alpha$-harmonic maps]{Spectral and convex geometry of $\alpha$-harmonic maps and sharp Morse index rigidity}
\author[Gao]{Rui Gao}
\address{School of Mathematics and Statistics, Shaanxi Normal University, 620 West Chang'an Avenue, Xi'an, Shaanxi, 710119, P. R. China}
\email{rui.gao@snnu.edu.cn}

\author[Zhu]{Miaomiao Zhu}
\address{School of Mathematical Sciences, Shanghai Jiao Tong University\\ 800 Dongchuan Road \\ Shanghai, 200240 \\ P. R. China}
\email{mizhu@sjtu.edu.cn}
\date{\today}
\allowdisplaybreaks[4]

\begin{abstract}
We develop a weighted scalar spectral geometric description of Sacks-Uhlenbeck $\alpha$-harmonic maps into round spheres and identify its convex geometry realization for $\alpha$-harmonic diffeomorphisms onto the round $2$-sphere. More precisely, we establish a one-to-one correspondence, modulo orientation preserving isometries of $\S^2$, between orientation preserving $\alpha$-harmonic diffeomorphisms onto the round $2$-sphere and positive Codazzi tensors satisfying a determinant bound and a compatibility condition with the prescribed metric on $\S^2$. These tensors are the radii tensors of auxiliary strictly convex surfaces, and the weighted scalar forms in spectral geometry equal $-\alpha$ times the second variations of the enclosed volumes with respect to the support functions. Using these spectral geometry and convex geometry results, we prove a sharp Morse index characterization of isometric equatorial parametrizations among $\alpha$-harmonic $2$-spheres in the round $n$-sphere. We also prove that every $\alpha$-harmonic diffeomorphism between unit round $2$-spheres is orthogonal for all $\alpha>1$. Equivalently, the only admissible positive Codazzi tensor whose reconstructed metric on $\S^2$ has Gaussian curvature 1 is $(\kappa+2)^{\alpha-1}\bar g$, where $\bar g$ is the unit round metric.
\end{abstract}
\maketitle

\tableofcontents
\section{Introduction and main results}\label{sec:introduction}

\subsection{Background and motivation}\label{subsec:background}

Let $(\Sigma,g)$ be a smooth closed oriented Riemannian surface, and let $(N,h)$ be a closed Riemannian manifold isometrically embedded in a fixed Euclidean space $\R^K$. For $\alpha>1$ and a fixed constant $\kappa>0$, define the $\alpha$-energy and the Dirichlet energy by
\begin{equation}\label{eq:energy-intro}
 E_{\alpha,\kappa}(u)=\frac12\int_\Sigma(\kappa+|du|_g^2)^\alpha\dd A_g,
 \qquad E(u)=\frac12\int_\Sigma|du|_g^2\dd A_g.
\end{equation}
% If $N = \S^n$, we use $\bar g$ to denote the unit round metric. When the
% domain is the unit round $\Sph^2$, $x:\Sph^2\hookrightarrow\R^3$
% denotes its standard embedding.

The Dirichlet energy $E$ does not in general satisfy the Palais-Smale condition on $W^{1,2}(\S^2,N)$, which obstructs a direct application of classical variational methods in \cite{palais1966foundations} to the construction of harmonic 2-spheres. To overcome this difficulty, Sacks-Uhlenbeck \cite{SacksUhlenbeck1981} introduced the standard $\alpha$-energy $E_\alpha:=E_{\alpha,1}$. The $\alpha$-energy approximation restores the compactness needed for variational constructions at each fixed $\alpha>1$. Passing to the limit as $\alpha\searrow1$ then yields possibly constant harmonic maps being the weak limit, with nonconstant harmonic 2-spheres arising as bubbles at points of energy concentration. By combining this approximation with topological variational constructions, they established the existence of branched minimal 2-spheres of arbitrary codimension in closed Riemannian manifolds under suitable nontriviality assumptions on the homotopy groups of the target. 
% For closed orientable surfaces of positive genus, they \cite{SacksUhlenbeck1982} further allowed the conformal structure of $\Sigma$ to vary and proved that, provided the induced homomorphism on fundamental groups has no nontrivial element represented by a simple closed curve in its kernel, there exists a conformal branched minimal immersion of least area among all branched immersions inducing the same action on the fundamental group. 
Building on the $\alpha$-energy perturbation, Micallef-Moore \cite{MicallefMoore1988} combined the existence of minimal 2-spheres admitting controlled Morse index with Morse index lower bounds under positive isotropic curvature to prove the vanishing of the homotopy groups in degrees from $2$ to $\lfloor n/2\rfloor$, and consequently proved that every closed simply connected manifold of dimension $n\geq4$ with positive isotropic curvature is homeomorphic to $\S^n$.

Beyond the $\alpha$-energy approximation, several other methods have been developed for the construction of harmonic maps, in particular, for the construction of branched immersed minimal 2-spheres. Struwe \cite{Struwe1985} developed the heat flow approach on surfaces, while Chen-Struwe \cite{ChenStruwe1989} introduced a Ginzburg-Landau type approximation in their study of the harmonic map heat flow. 
% The corresponding compactness and bubbling theory for critical points of the Ginzburg-Landau approximation was subsequently developed by Lin-Wang \cite{LinWang1999,LinWang2002}. 
Colding-Minicozzi \cite{ColdingMinicozzi2008} introduced a min-max construction based on harmonic replacement. Since these constructions use different approximations and different variational families, it is natural to compare not only the resulting min-max values, but also the limiting critical maps. A closely related question arises in the construction of minimal hypersurfaces by the Almgren-Pitts min-max theory \cite{Pitts1981} and the Allen-Cahn min-max approach developed by Guaraco \cite{Guaraco2018} and subsequent works. Dey \cite{Dey2022} proved equality of the normalized Almgren-Pitts and Allen-Cahn min-max values together with an inclusion of the associated critical sets, whereas Mantoulidis \cite{Mantoulidis2024} showed that phase transition limits can satisfy additional rigidity and form a proper subset of the Almgren-Pitts critical set.  

For the $\alpha$-energy approximation, the corresponding selection problem is a main motivation of the present paper. It is natural to consider the following 
\\

{\bf Question}: Are the harmonic $2$-spheres arising from $\alpha$-energy approximation single out a distinguished subclass within the broader family of harmonic maps? 
\\

Such a comparison needs essentially distinguish parametrized maps from their geometric images and specify the $\alpha$-energy and Morse index bounds. Our main observation is that the $\alpha$-harmonic 2-sphere into round $n$-sphere exhibits a genuine selection phenomenon. More precisely, under suitable Morse index bounds, the limiting harmonic $2$-spheres arising from $\alpha$-harmonic $2$-spheres satisfy additional rigidity in their parametrizations, a feature that cannot be detected by comparing min-max values or geometric images alone.

A basic reason to expect greater rigidity for $\alpha$-harmonic maps is that the perturbation breaks the conformal reparametrization symmetry of the Dirichlet energy. On a surface, the Dirichlet energy is invariant under conformal changes of parametrization, so a harmonic map from $\S^2$ generates an entire M\"obius orbit of harmonic maps with the same energy. In contrast, for $\alpha>1$, the functional $E_{\alpha,\kappa}$ is only invariant under isometries of the domain but not under general conformal reparametrizations. Thus the noncompact M\"obius symmetry of the Dirichlet energy is reduced to the compact isometry group of the round sphere, suggesting that critical points of the $\alpha$-energy functional should satisfy stronger restrictions on their parametrizations. For maps between round 2-spheres, Lamm-Malchiodi-Micallef \cite{LammMalchiodiMicallef2020,LammMalchiodiMicallef2021} discovered a selection phenomenon for the $\alpha$-energy approximation. In \cite[Theorem 1.4]{LammMalchiodiMicallef2021}, they show that degree $\pm1$ $\alpha$-harmonic maps satisfying the low energy threshold must be orthogonal when $\alpha$ is sufficiently close to 1. Moreover, they explicitly asked whether this low energy restriction could be replaced by a controlled Morse index condition. Theorem~\ref{thm:main-full-gap} addresses this question for higher dimensional round target $\Sph^n$, $n\geq3$, under a uniform bound on the $\alpha$-energy. Within that $\alpha$-energy bound, the Theorem \ref{thm:main-full-gap} replaces the low energy threshold by a sharp bound on the Morse index, without assuming a degree or equatorial $\S^2$ containment of the image. For closed non-spherical $\Sigma$ equipped with the normalized constant curvature metrics, Sharp \cite{Sharp2024} classifies low energy $\alpha$-harmonic maps into $\S^2$, showing in particular that bubbling may occur and that the bubbling location and scale are determined by the holomorphic 1-forms of the $\Sigma$.

A further motivation of present paper comes from spectral geometry. The classical work of Takahashi \cite{Takahashi1966}, Hersch \cite{Hersch1970}, and Li-Yau \cite{LiYau1982} relates sphere valued harmonic maps and minimal immersions to coordinate eigenfunctions, extremal eigenvalues, and conformal geometry. More recent developments include the index estimates of harmonic spheres by Karpukhin \cite{Karpukhin2021}, the sharp eigenvalue inequalities by Karpukhin-Nadirashvili-Penskoi-Polterovich \cite{KarpukhinNadirashviliPenskoiPolterovich2021}, the min-max characterization of conformal eigenvalues by Karpukhin-Stern \cite{KarpukhinStern2024}, and the stability results of Karpukhin-Nahon-Polterovich-Stern \cite{KarpukhinNahonPolterovichStern2025} etc. The connections between harmonic maps, Morse index, and eigenvalue optimization suggest studying rigidity through the scalar equation satisfied by the coordinate functions of sphere valued harmonic maps. For $\alpha$-harmonic maps, the associated scalar form \eqref{eq:scalar-form-intro} carries the conductivity $(\kappa+|du|_g^2)^{\alpha-1}$, which cannot generally be absorbed into a conformal change of the domain metric in dimension 2. Thus, the weighted spectral theory for $\alpha$-harmonic maps must be newly developed.

In this paper, we establish a sharp Morse index criterion of rigidity for the Sacks-Uhlenbeck $\alpha$-energy approximation. Under a fixed $\alpha$-energy bound, Theorem~\ref{thm:main-full-gap} shows that every smooth nonconstant $\alpha$-harmonic map from the unit round $2$-sphere into $\Sph^n$, $n\geq3$, with Morse index strictly below $3(n-2)$ is an isometric equatorial parametrization when $\alpha>1$ is sufficiently close to 1. In particular, no assumption on degree or equatorial containment is required. Equatorial inclusions of degree 2 minimizers attain this Morse index threshold, proving sharpness for the rigidity conclusion. Theorems~\ref{thm:main-spectral} and~\ref{thm:main-convex-energy} establish the weighted scalar spectral theory and its realization through positive Codazzi tensors and mixed volumes in convex geometry. On the unit round domain, Theorem~\ref{cor:main-rigidity} further proves rigidity at scalar index 1 for every $\alpha>1$.

The convex geometry aspect of our method is a correspondence between $\alpha$-harmonic diffeomorphisms and positive Codazzi tensors on the target round $2$-sphere. This correspondence gives an exact geometric description of the weighted scalar form \eqref{eq:scalar-form-intro}, although its conductivity cannot be absorbed into a conformal change of the domain metric. A separate feature of the sharp Morse index argument is the use of global matching variations associated with the rotational Jacobi fields of the bubbles. These variations detect negative directions that are not accounted for by adding the Morse indices of the limiting bubble tree components as $\alpha \searrow 1$. We describe the main ideas of proofs for these results together in Subsection~\ref{subsec:proof-outline}.

\subsection{Sharp Morse index rigidity on the round sphere}\label{subsec:main-index}

For an $\alpha$-harmonic map $u\in W^{1,2\alpha}(\Sigma,\S^n)$ from a closed surface $(\Sigma,g)$, by the standard inclusion $\S^n\hookrightarrow\R^{n+1}$, $u = (u^1,\dots,u^{n+1})$ solves the scalar equation
\begin{equation}\label{eq:equation-intro}
 -\Div_g\bigl((\kappa+|du|_g^2)^{\alpha-1}\nabla u\bigr)
 =(\kappa+|du|_g^2)^{\alpha-1}|du|_g^2u.
\end{equation}
Consequently, every Euclidean coordinate function of $u$ lies in the kernel of the self adjoint scalar operator
\begin{equation}\label{eq:operator-intro}
 J_{\alpha,u}f
 =-(\kappa+|du|_g^2)^{1-\alpha}
 \Div_g\bigl((\kappa+|du|_g^2)^{\alpha-1}\nabla f\bigr)
 -|du|_g^2f
 \quad\text{on }
 L^2\bigl(\Sigma,(\kappa+|du|_g^2)^{\alpha-1}\dd A_g\bigr).
\end{equation}
Since $u$ is smooth, the corresponding scalar bilinear form is defined for real valued functions $f,h\in W^{1,2}(\Sigma)$ by
\begin{equation}\label{eq:scalar-form-intro}
 Q_{\alpha,u}(f,h)
 =\alpha\int_\Sigma
 (\kappa+|du|_g^2)^{\alpha-1}
 \bigl(\ip{\nabla f}{\nabla h}_g-|du|_g^2fh\bigr)\dd A_g.
\end{equation}
Integration by parts identifies its associated scalar operator on the weighted $L^2$ space in \eqref{eq:operator-intro} as $\alpha J_{\alpha,u}$.

We denote by $I_\alpha(u)$ the number of negative eigenvalues of $J_{\alpha,u}$, counted with multiplicity, and by $N_\alpha(u)$ the dimension of its kernel. Equivalently, these are the index and nullity of $Q_{\alpha,u}$, respectively. These scalar spectral quantities are distinct from $\ind_{\Sph^n}E_{\alpha,\kappa}(u)$ and $\nul_{\Sph^n}E_{\alpha,\kappa}(u)$, which denote the Morse index and nullity of the second variation on sections of $u^{-1}T\S^n$. At $\alpha=1$, the same notation refers to the scalar operator $J_{1,u}=\Delta_g-|du|_g^2$ and to the Jacobi operator associated to the second variation of Dirichlet energy. The connection between these quantities is given by an exact splitting of the second variation, that is, when the image of $u$ lies in an equatorial subsphere, each independent normal direction contributes $I_\alpha(u)$ to the Morse index and $N_\alpha(u)$ to the nullity, see Theorem \ref{thm:main-spectral} below for more details.

As the main application of our spectral and convex geometric theory of $\alpha$-harmonic maps into round spheres, we establish a sharp Morse index characterization of isometric equatorial parametrizations among $\alpha$-harmonic $2$-spheres in the round $n$-sphere. No assumption on the degree, scalar index, immersion, or equatorial containment of the image is imposed for the $\alpha$-harmonic $2$-sphere.

\begin{theorem}\label{thm:main-full-gap}
Fix $n\geq3$, $\kappa>0$, and $\Lambda<\infty$. There is $\alpha_0=\alpha_0(n,\kappa,\Lambda)>1$ such that every smooth nonconstant $\alpha$-harmonic map $u:\Sph^2\to\Sph^n$ from the unit round 2-sphere into the unit round $n$-sphere satisfying
\begin{equation}\label{eq:main-sharp-hyp}
 1<\alpha<\alpha_0,\qquad E_{\alpha,\kappa}(u)\leq\Lambda,\qquad
 \ind_{\Sph^n}E_{\alpha,\kappa}(u)<3(n-2)
\end{equation}
is an isometric parametrization of an equatorial 2-sphere. More precisely, $u(x)=Lx$ for a linear map $L:\R^3\rightarrow\R^{n+1}$ satisfying $L^TL=\Id$. The Morse index and nullity of the isometric parametrization of an equatorial 2-sphere are $n-2$ and $3n-3$, respectively.

The strict Morse index threshold $3(n-2)$ for the rigidity conclusion is sharp. For every fixed $\kappa>0$ and every $\alpha>1$ sufficiently close to 1, there is a smooth degree 2 minimizer $v_\alpha:\Sph^2\to\Sph^2$ of $E_{\alpha,\kappa}$ in its degree class. If $\iota:\Sph^2\hookrightarrow\Sph^n$ is equatorial inclusion, then
\begin{equation}\label{eq:main-sharpness}
 \ind_{\Sph^n}E_{\alpha,\kappa}(\iota\circ v_\alpha)=3(n-2), \qquad E_{\alpha,\kappa}(v_\alpha)\longrightarrow2\pi\kappa+8\pi.
\end{equation}
Since the equatorially included degree 2 minimizers are not isometric parametrizations, for every $\Lambda>2\pi\kappa+8\pi$, the strict Morse index upper bound in \eqref{eq:main-sharp-hyp} cannot be replaced by a nonstrict bound.
\end{theorem} 

\begin{rmk}
    Combining Karpukhin's \cite{Karpukhin2021} comparison between the Morse index and the scalar spectral index  with Ejiri's \cite{Ejiri1998} harmonic degree estimate  yields the following classification for harmonic $2$-spheres: a nonconstant harmonic map \(v:S^{2}\to S^{n}\), \(n\geq3\), between unit round spheres has Morse index strictly less than \(3(n-2)\) if and only if it is a M\"{o}bius or anti-M\"{o}bius parametrization of an equatorial 2-sphere. Every such map has Morse index \(n-2\) by the normal splitting formula \eqref{eq:index-splitting} at $\alpha=1$.  In contrast, Theorem \ref{thm:main-full-gap} shows that, under a uniform bound on the $\alpha$-energy and for \(\alpha>1\) sufficiently close to 1, the same strict Morse index bound forces an $\alpha$-harmonic map to be an isometric equatorial parametrization. Thus the $\alpha$-energy perturbation eliminates the nonisometric conformal reparametrizations permitted by the harmonic map equation, giving stronger rigidity of the parametrization. Theorem \ref{thm:main-full-gap} therefore provides a complete and sharp Morse index characterization for $\alpha$-harmonic 2-spheres with uniformly bounded $\alpha$-energy, answering the open question posed by Lamm–Malchiodi–Micallef \cite{LammMalchiodiMicallef2021} concerning the replacement of their low energy restriction by a Morse index bound. Moreover, our result applies to the high dimensional round targets \(S^{n}\), \(n\geq3\), without assuming a low $\alpha$-energy threshold, a degree requirement, or $2$-sphere equatorial containment of the image. 
\end{rmk}

% The additional index contribution in the configuration of 2 bubbles comes from variations coupling their rotational kernels. The construction retains the positive Sacks-Uhlenbeck term in the Hessian and gives exact orthogonality to the negative directions supported in the bubble cores. These are global variations, consistent with the positivity of the Hessian for variations supported on the neck with zero boundary values in \cite{DaLioRiviereSchlagenhauf2025}.

\subsection{Spectral geometry on closed surfaces and its convex realization}\label{subsec:main-spectral}

Let $\mathfrak g$ denote the genus of the domain $\Sigma$ and let $r(u)$ be the dimension of the smallest equatorial sphere containing $u(\Sigma)$. Theorem~\ref{thm:main-spectral} below describes the scalar spectral structure of $\alpha$-harmonic maps from $(\Sigma,g)$ into round $n$-sphere and its relation to the Morse index and nullity, and classifies $\alpha$-harmonic maps of scalar index 1 on domains of genus 0.

\begin{theorem}[Scalar spectral structure]\label{thm:main-spectral}
Let $u:(\Sigma,g)\to\Sph^n$ be a smooth nonconstant $\alpha$-harmonic map, with $\alpha>1$ and $\kappa>0$.
\begin{enumerate}
\item\label{thm:main-spectral 1} The weighted variational eigenvalue problem
\begin{equation}\label{eq:main-weighted-eigenproblem}
\begin{aligned}
 \int_\Sigma(\kappa+|du|_g^2)^{\alpha-1}
 \ip{\nabla f}{\nabla h}\dd A_g
 &=\lambda\int_\Sigma fh\dd\nu_{\alpha,u}
 \quad\text{for every }h\in W^{1,2}(\Sigma),\\
 \dd\nu_{\alpha,u}
 &:=\frac12(\kappa+|du|_g^2)^{\alpha-1}|du|_g^2\dd A_g,
\end{aligned}
\end{equation}
has discrete eigenvalues, with simple first eigenvalue $\lambda_0=0$. Every nonzero coordinate function of $u$ has eigenvalue 2. The number of eigenvalues below 2 is $I_\alpha(u)$ and the multiplicity of eigenvalue 2 is $N_\alpha(u)$. In particular, $I_\alpha(u)=1$ precisely when the first positive eigenvalue is 2.
\item\label{thm:main-spectral 2} Every nonzero function in $\ker J_{\alpha,u}$ has at most $I_\alpha(u)+1$ many nodal domains. Moreover,
\begin{equation}\label{eq:main-general-multiplicity}
 N_\alpha(u)\leq2I_\alpha(u)+4\mathfrak g+1,
 \qquad r(u)\leq2I_\alpha(u)+4\mathfrak g.
\end{equation}
\item\label{thm:main-spectral 3} If $u$ takes values in an equatorial $\Sph^r\subset\Sph^m$, then
\begin{align*}
 \ind_{\Sph^m}E_{\alpha,\kappa}(u)
 &=\ind_{\Sph^r}E_{\alpha,\kappa}(u)+(m-r)I_\alpha(u),\\
 \nul_{\Sph^m}E_{\alpha,\kappa}(u)
 &=\nul_{\Sph^r}E_{\alpha,\kappa}(u)+(m-r)N_\alpha(u).
\end{align*}
\end{enumerate}
If $\Sigma$ is diffeomorphic to $\Sph^2$  with any smooth metric $g$, then $I_\alpha(u)=1$ if and only if $u$ is a diffeomorphism onto an equatorial 2-sphere. In this case $N_\alpha(u)=3$ and the scalar kernel $\ker J_{\alpha,u}$ is spanned by the 3 coordinate functions of $u$ in the equatorial target $2$-sphere.
\end{theorem}
The density of $\nu_{\alpha,u}$ in \eqref{eq:main-weighted-eigenproblem} may vanish where $du=0$. In Proposition~\ref{prop:weighted-spectrum}, we construct the spectrum using the weighted Dirichlet form $\mathfrak d_{\alpha,u}$ and the identities
\[
\begin{aligned}
 \mathfrak h(f,h)=\mathfrak d_{\alpha,u}(f,h)+\int_\Sigma fh\dd\nu_{\alpha,u},\qquad
 \mathfrak h(Tf,h)=\int_\Sigma fh\dd\nu_{\alpha,u}.
\end{aligned}
\]
The strict positivity of $(\kappa+|du|_g^2)^{\alpha-1}$ and the positive total mass of $\nu_{\alpha,u}$ make the norm induced by $\mathfrak h$ equivalent to the $W^{1,2}(\Sigma)$ norm. The operator $T$ is nonnegative and self adjoint with respect to $\mathfrak h$, and compact by Rellich's theorem, moreover its positive eigenvalues $t$ give the discrete eigenvalues $\lambda=t^{-1}-1$ in \eqref{eq:weighted-variational-spectrum}. This construction requires no immersion assumption or condition on the zero set of $du$. At $\alpha=1$, $\mathfrak d_{1,u}$ is the Dirichlet form of $g$ and $\nu_{1,u}$ is the area measure of the possibly degenerate metric $\frac12|du|_g^2g$, so \eqref{eq:weighted-variational-spectrum} is its variational Laplace eigenvalue problem. For $\alpha>1$, the conductivity $(\kappa+|du|_g^2)^{\alpha-1}$ in $\mathfrak d_{\alpha,u}$ cannot be removed by a conformal change of $g$, which prevents many methods established in the harmonic map setting from being applied directly.

The characterization of $\alpha$-harmonic maps at scalar index $I_\alpha(u)=1$ in Theorem~\ref{thm:main-spectral} holds for every smooth metric $g$ on a closed domain $\Sigma$ of genus 0. We now relate the scalar form \eqref{eq:scalar-form-intro} to convex geometry by expressing it on the unit round target $(\Sph^2,\bar g)$.

Fix $\alpha>1$ and $\kappa>0$, and let $u:(\Sigma,g)\to(\Sph^2,\bar g)$ be a smooth orientation preserving diffeomorphism with positive Jacobian $J$ defined by $u^*\dd A_{\bar g}=J\dd A_g$. Using $\bar g$ to identify symmetric endomorphisms with covariant 2-tensors, define at $y=u(x)$
\begin{equation}\label{eq:tensors-intro}
 A_y=\frac{\bigl(\kappa+|du_x|_g^2\bigr)^{\alpha-1}}{J(x)}du_xdu_x^*,\qquad
 C=\cof_{\bar g}A=(\tr_{\bar g}A)\bar g-A.
\end{equation}
Here $du_x^*:T_y\Sph^2\to T_x\Sigma$ is the adjoint of $du_x$ with respect to $g$ and $\bar g$. For a positive definite endomorphism $A$, the cofactor is $(\det_{\bar g}A)A^{-1}$. Here, \eqref{eq:tensors-intro} writes the same operation as a covariant 2-tensor. The tensor $A$ is the coefficient of the gradient term when $\alpha^{-1}Q_{\alpha,u}$ is transported to the target, as in \eqref{eq:pushforward}. If $u$ is $\alpha$-harmonic, Proposition~\ref{prop:target-conductivity} gives $\Div_{\bar g}A=0$, which by \eqref{eq:codazzi-equivalence} is equivalent to the Codazzi condition $\nabla_iC_{jk}=\nabla_jC_{ik}$ of $C$, where $\nabla$ is the Levi-Civita connection of $\bar g$. For a smooth positive definite tensor $C$ on $(\Sph^2,\bar g)$ satisfying $\det_{\bar g}C>\kappa^{2(\alpha-1)}$, define the Riemannian metric
\begin{equation}\label{eq:GC-intro}
 G_C=
 \frac{\tr_{\bar g}C}
 {(\det_{\bar g}C)\bigl((\det_{\bar g}C)^{1/(2(\alpha-1))}-\kappa\bigr)}\,C.
\end{equation}
When $C$ is obtained from $u$ by \eqref{eq:tensors-intro}, the identities \eqref{eq:C-invariants} give $G_C=(u^{-1})^*g$. Theorem~\ref{thm:main-convex-energy} below establishes a bijection between orientation preserving $\alpha$-harmonic diffeomorphisms modulo precomposition by orientation preserving isometries of $(\Sigma,g)$ and smooth positive definite Codazzi tensors $C$ satisfying the determinant bound $\det_{\bar g}C>\kappa^{2(\alpha-1)}$ for which $(\Sph^2,G_C)$ is orientation preservingly isometric to $(\Sigma,g)$.
\begin{theorem}[Convex realization and reconstruction]\label{thm:main-convex-energy}
For every orientation preserving $\alpha$-harmonic diffeomorphism $u:(\Sigma,g)\to(\Sph^2,\bar g)$, the tensor $A$ is divergence free and $C=\cof A$ is a positive Codazzi tensor. There is a smooth function $h$, unique modulo first spherical harmonics such that
\[
 C=\Hess_{\bar g}h+h\bar g.
\]
After adding a first spherical harmonic, $h>0$, and the map
\begin{equation}\label{eq:convex-X}
 X(y)=\nabla_{\bar g}h(y)+h(y)y
\end{equation}
is a strictly convex embedding in $\R^3$ with Gauss map $y$ and principal
radii equal to the eigenvalues of $C$.

The scalar form \eqref{eq:scalar-form-intro} transported by $u^{-1}$ equals $-\alpha B_C$, where
\[
 B_C(f,k)=\int_{\Sph^2}\left((\tr C)fk-\cof C(\nabla f,\nabla k)\right)\dd A_{\bar g}.
\]
The bilinear form $B_C$ has exactly 1 positive direction, its kernel consists of the first spherical harmonics, and
\[
 B_C(f,f)B_C(1,1)\leq B_C(f,1)^2,
\]
with equality precisely for constants plus first spherical harmonics.

Conversely, positive definite Codazzi tensors $C$ with $\det_{\bar g}C>\kappa^{2(\alpha-1)}$ for which $(\Sph^2,G_C)$ is orientation preservingly isometric to $(\Sigma,g)$ correspond to orientation preserving $\alpha$-harmonic diffeomorphisms onto the round target, modulo orientation preserving  isometries of the domain 2-sphere. For the unit round domain 2-sphere, the metric isometric condition is exactly a prescribed Gaussian curvature condition $K_{G_C}=1$.
\end{theorem}

The signature of $B_C$ in Theorem~\ref{thm:main-convex-energy} is the classical smooth 2-dimensional Alexandrov-Fenchel signature, that is, the form has exactly 1 positive direction and its kernel consists of the first spherical harmonics, see \cite{Alexandrov1938,GuanMaTrudingerZhu2010,ShenfeldVanHandel2019}. Our contribution is the identity \eqref{eq:Q-minus-B}, which identifies the transported scalar form \eqref{eq:scalar-form-intro} with $-\alpha B_C$, together with the reconstruction of the map from the Codazzi tensor in Proposition~\ref{prop:general-reconstruction}. This identification makes the classical signature theorem applicable to the weighted scalar equation \eqref{eq:equation-intro} and gives the characterization $I_\alpha(u)=1$ and $N_\alpha(u)=3$ for $\alpha$-harmonic diffeomorphisms.

For the unit round domain, Theorem \ref{thm:main-convex-energy} and Theorem~\ref{thm:reconstruction} expresses the metric isometric condition of $(\S^2, G_C)$ and $(\Sigma,g)$ as
\[
 \mathrm{det}_{\bar g}C>\kappa^{2(\alpha-1)},\qquad K_{G_C}=1,
\]
where $C$ is a smooth positive definite Codazzi tensor and $G_C$ is defined by \eqref{eq:GC-intro}. This is a prescribed Gaussian curvature problem with $C$ as its unknown. It is analogous to the 2-dimensional constant curvature problem underlying the Yamabe problem, but the metric here varies through the explicit family $G_C$, rather than through conformal changes of a fixed metric.

Both this problem and the classical Minkowski problem \cite{Nirenberg1953} can be written in terms of a support function $h$ with positive definite radii tensor $C=\Hess_{\bar g}h+h\bar g$. In the Minkowski problem, one prescribes the Gaussian curvature $(\det_{\bar g}C)^{-1}$ of the convex surface $X$ in \eqref{eq:convex-X} as a function of its outward unit normal. Here, the support function $h$ determines $G_C$ through \eqref{eq:GC-intro}, and the prescribed curvature $K_{G_C}=1$ is that of this reconstructed domain metric $G_C$.

The following Theorem~\ref{cor:main-rigidity} establishes both rigidity of $\alpha$-harmonic diffeomorphisms between unit round $2$-spheres and uniqueness of the positive Codazzi tensor solving this prescribed Gaussian curvature problem, for every $\alpha>1$ and $\kappa>0$.

% The signature in Theorem \ref{thm:main-convex-energy} is the classical smooth 2-dimensional Alexandrov-Fenchel signature \cite{Alexandrov1938,GuanMaTrudingerZhu2010,ShenfeldVanHandel2019}. 
% The convex surface associated with the transported scalar equation is auxiliary and is distinct from the image of $u$, which lies in a round sphere. 
% A global diffeomorphism to $\Sph^2$ necessarily has domain of genus 0. Thus Theorem~\ref{thm:main-convex-energy} is not asserted for arbitrary maps from surfaces of positive genus or branched covers; the general spectral statements in Theorem~\ref{thm:main-spectral} have no such restriction. 
% If $\Sigma$ the unit round $2$-sphere, we obtain the following rigidity statement for the entire range $\alpha>1$.
\begin{theorem}[Rigidity at scalar index 1]\label{cor:main-rigidity}
For every $\alpha>1$ and $\kappa>0$, every smooth $\alpha$-harmonic diffeomorphism between unit round 2-spheres is orthogonal. Consequently, a smooth nonconstant $\alpha$-harmonic map $u:\Sph^2\to\Sph^n$ has $I_\alpha(u)=1$ precisely when the map is an isometric equatorial parametrization. The only positive Codazzi tensor satisfying $\det_{\bar g}C>\kappa^{2(\alpha-1)}$ and $K_{G_C}=1$ is $C=(\kappa+2)^{\alpha-1}\bar g$.
\end{theorem}
Here, we note that Theorem~\ref{cor:main-rigidity} holds for every $\alpha>1$ without an energy bound, only under the assumption that the $\alpha$-harmonic map is a diffeomorphism. While Theorem~\ref{thm:main-full-gap} concerns arbitrary smooth nonconstant maps into $\Sph^n$, $n\geq3$, satisfying the energy and Morse index bounds in \eqref{eq:main-sharp-hyp}, for exponents $\alpha$ sufficiently close to 1.

\begin{remark}\label{rem:determinant-nondegeneracy}
Fix $\alpha>1$. Since $C$ is smooth and positive definite on the compact sphere, we have
\[
 \min_{\Sph^2}(\mathrm{det}_{\bar g}C)^{1/(2(\alpha-1))}>0.
\]
Thus every such tensor satisfies the determinant bound for all sufficiently small $\kappa>0$. For fixed $\kappa$, the condition $\det_{\bar g}C>\kappa^{2(\alpha-1)}$ gives a strictly positive minimum for $(\det_{\bar g}C)^{1/(2(\alpha-1))}-\kappa$, so \eqref{eq:GC-intro} defines a smooth nondegenerate Riemannian metric $G_C$. Together with $C>0$, this condition ensures that the curvature equation $K_{G_C}=1$, expressed in terms of the support function $h$ from Lemma~\ref{lem:support}, has a uniformly elliptic linearization at each fixed admissible $h$. The determinant bound $\det_{\bar g}C>\kappa^{2(\alpha-1)}$ therefore specifies the range in which the reconstruction is nondegenerate and the prescribed Gaussian curvature equation is elliptic. 
% The parameter $\kappa$ remains fixed in each problem, since changing $\kappa$ changes $G_C$ and hence the curvature equation.
\end{remark}

\subsection{Outline of the proofs}\label{subsec:proof-outline}
We now outline the main ideas behind the proofs of our main Theorems \ref{thm:main-full-gap}, \ref{thm:main-spectral}, \ref{thm:main-convex-energy} and \ref{cor:main-rigidity}.

The proof of Theorem~\ref{thm:main-full-gap} requires both rigidity of $\alpha$-harmonic 2-sphere at scalar index 1 and an additional Morse index estimate that excludes the remaining configuration of 2 bubbles.

The starting point is the weighted scalar equation \eqref{eq:equation-intro} satisfied by the coordinate functions of the $\alpha$-harmonic map $u$. Proposition~\ref{prop:weighted-spectrum} identifies the scalar index with the number of variational eigenvalues below 2 and the scalar nullity with the multiplicity of 2. Combining the nodal domain bound in Lemma~\ref{lem:nodal-count} with the multiplicity estimate \eqref{eq:general-multiplicity}, obtained from the conjugation \eqref{eq:conjugation}, shows that on a domain of genus 0, scalar index 1 places the image of $u$ in an equatorial $2$-sphere. Observe that every nonzero height function has exactly 2 nodal domains and, by Lemma~\ref{lem:nodal-critical}, has no critical zero point. A loss of rank would produce such a zero, so Proposition~\ref{prop:index-one-forward} identifies the $\alpha$-harmonic map with scalar index 1 as a diffeomorphism onto its image. Under equatorial inclusion, Proposition~\ref{prop:normal-splitting} then turns this scalar information into Morse index and nullity formulas.

The converse for the last assertion of Theorem \ref{thm:main-spectral} requires scalar index 1 for every $\alpha$-harmonic diffeomorphism. We prove it by transporting the scalar form \eqref{eq:scalar-form-intro} to the round target $2$-sphere. Proposition~\ref{prop:target-conductivity} gives a divergence free tensor $A$, whose cofactor $C$ is a positive definite Codazzi tensor by \eqref{eq:codazzi-equivalence}. Lemma~\ref{lem:support} realizes $C$ as the radii tensor of an auxiliary strictly convex surface. The new link with convex geometry is the identity \eqref{eq:Q-minus-B}, which identifies the transported scalar form with $-\alpha B_C$. The classical signature and translation kernel of $B_C$ in Proposition~\ref{prop:signature} therefore give scalar index 1 and nullity 3, as stated in Proposition~\ref{prop:automatic-index}. These arguments complete the classification in Theorem~\ref{thm:main-spectral}, while $C$ also determines the transported domain metric $G_C$. For positive Codazzi tensors with $\det C>\kappa^{2(\alpha-1)}$, Proposition~\ref{prop:general-reconstruction} recovers the map when $(\Sph^2,G_C)$ is isometric to the domain $(\Sigma,g)$, proving the correspondence in Theorem~\ref{thm:main-convex-energy}. For the unit round $2$-sphere domain, Theorem~\ref{thm:reconstruction} expresses this compatibility as $K_{G_C}=1$.

The convex description also provides the uniform $\alpha$-energy bound needed to prove that every $\alpha$-harmonic diffeomorphism of the unit round $2$-sphere is an orthogonal map when $\alpha>1$ is sufficiently close to 1. More precisely,  combining the mixed-volume estimate in Proposition~\ref{prop:barycenter} with the domain balance identity \eqref{eq:MN} yields the uniform energy bound in Theorem~\ref{thm:energy-bound}. For a sequence of orientation preserving diffeomorphisms with $\alpha_j\searrow1$, the energy identity proved in \cite{LiZhu2019}, see also Theorem~\ref{thm:compactness-input}, and nonnegativity of the limiting Jacobians give degree splitting with exactly 1 nonconstant harmonic component, of degree 1. Lemma~\ref{lem:noncollapse} rules out the possibility that this component forms a bubble, giving the smooth convergence of $\alpha$-harmonic diffeomorphisms in Proposition~\ref{prop:no-bubbling-diffeo}. The normalized balance identity then identifies the M\"obius limit as an orthogonal map in Proposition~\ref{prop:rotational-limit}. Smooth convergence alone does not give rigidity of the approximating $\alpha$-harmonic maps, because the Hessian degenerates in the nonrotational conformal directions as $\alpha\searrow1$. Lemma~\ref{lem:sequential-local} excludes these directions by applying the linearized balance identity to normalized displacement fields, proving Theorem~\ref{thm:near-one-rigidity}.

To extend rigidity to every $\alpha>1$, we must control the Jacobi kernel at a hypothetical nonorthogonal diffeomorphism. The index and nullity bound in Proposition~\ref{prop:six-dimensional-bound}, together with the exclusion of continuous isotropy using Proposition~\ref{prop:equivariant} and Lemma~\ref{lem:continuous-isotropy}, shows that the 6 vector fields generated by domain and target rotations exhaust this kernel, as proved in Proposition~\ref{prop:nonround-nondegenerate}. The implicit function argument in Lemma~\ref{lem:parameter-continuation} continues such maps in the exponent and gives local rigidity near orthogonal maps. Together with the smooth compactness in Corollary~\ref{cor:compact-diffeo-branch} for exponents $\alpha$ in compact subsets of $(1,\infty)$, these statements make the set of exponents $\alpha$ admitting nonorthogonal $\alpha$-harmonic diffeomorphisms open and relatively closed in $(1,\infty)$. Theorem~\ref{thm:near-one-rigidity} forces this set to be empty, proving Theorem~\ref{thm:global-diffeo-rigidity}. The spectral index 1 classification and reconstruction of $\alpha$-harmonic 2-spheres then give Theorem~\ref{cor:main-rigidity}, as recorded in Corollary~\ref{cor:global-scalar-rigidity}.

We now consider a sequence of counterexamples $u_{\alpha_j}:\Sph^2\to\Sph^n$ to Theorem~\ref{thm:main-full-gap}, with $\alpha_j\searrow1$ and the $\alpha$-energy and Morse index bounds in \eqref{eq:main-sharp-hyp}. Bubble tree convergence in Theorem~\ref{thm:compactness-input}, applied to $u_{\alpha_j}$, gives a harmonic base map and finitely many harmonic bubbles after passing to a subsequence. Morse index lower semicontinuity in \cite{Moore-Ream2019}, see also Theorem~\ref{thm:full-index-input}, together with \eqref{eq:harmonic-index-comparison} and Lemmas~\ref{lem:harmonic-scalar-gap} and~\ref{lem:harmonic-sphere-index-one}, leaves at most 2 nonconstant components, each an equatorial $2$-sphere of harmonic degree 1. The integrated coordinate equation \eqref{eq:hemisphere-balance} for $u_{\alpha_j}$ excludes a constant weak limit in the absence of bubbles. To determine how these components can occur, we use the balance identity and mass bounds in Lemma~\ref{lem:balance-atoms}, which exclude both a single bubble over a constant base and a bubble accompanying a nonconstant base. In the absence of bubbles, spectral convergence gives $I_{\alpha_j}(u_{\alpha_j})=1$ for all sufficiently large $j$, and Corollary~\ref{cor:global-scalar-rigidity} identifies $u_{\alpha_j}$ as an isometric equatorial parametrization. Thus Proposition~\ref{prop:two-bubble-reduction} reduces the remaining case for $u_{\alpha_j}$ to a constant base and 2 bubbles at distinct antipodal points.

The Morse indices of the 2 bubbles alone account for only $2(n-2)$ many directions, so they cannot reach the threshold $3(n-2)$. The additional $(n-2)$ directions in Theorem~\ref{thm:two-bubble-index} come from rotational Jacobi fields along $u_{\alpha_j}$ joined across the neck region. To evaluate the second variation of $E_{\alpha_j,\kappa}$ at $u_{\alpha_j}$ on these variations, we need the neck expansion to first order in Proposition~\ref{prop:first-neck-expansion} inspired by \cite{Yin2022}. The conductivity and flux identities in Lemma~\ref{lem:gnomonic-flux}, together with the estimates in Lemmas~\ref{lem:neck-decay} and~\ref{lem:gnomonic-regularity} and the comparison between $\alpha_j-1$ and the bubble scales in Lemma~\ref{lem:neck-Hamiltonian}, yield the desired expansion. Lemma~\ref{lem:odd-rotational-mode} then identifies rotational vector fields whose leading profiles change sign across the middle neck cylinder. Reversing an exact rotational Jacobi field on one bubbling cap and joining Jacobi fields on two bubbling caps by harmonic interpolation followed by tangent projection across that cylinder gives the variations in Lemma~\ref{lem:negative-bridge}. The limit \eqref{eq:bridge-negative-limit} shows that the second variation of $E_{\alpha_j,\kappa}$ at $u_{\alpha_j}$ is negative definite on their $(n-2)$-dimensional span. To add these directions to those supported in the bubbling regions, the Jacobi identity in Lemma~\ref{lem:Killing-kernel} gives the exact Hessian orthogonality \eqref{eq:exact-core-orthogonality}. Theorem~\ref{thm:two-bubble-index} therefore gives a Morse index lower bound of $3(n-2)$ for $u_{\alpha_j}$ when $j$ is sufficiently large, contradicting \eqref{eq:main-sharp-hyp} and proving the rigidity assertion of Theorem~\ref{thm:main-full-gap}.

To prove the sharpness part of Theorem~\ref{thm:main-full-gap} within the same energy restriction, we use the degree 2 minimizers $v_\alpha:\Sph^2\to\Sph^2$ of $E_{\alpha,\kappa}$ in Proposition~\ref{prop:degree-two-minimizers}. Comparison with a fixed holomorphic map $\Sph^2\to\Sph^2$ of degree 2 gives the entropy estimate \eqref{eq:degree-two-entropy}, and the uniform integrability in \eqref{eq:degree-two-set-energy} excludes bubbling. For any sequence $\alpha_j\searrow1$, the minimizers $v_{\alpha_j}$ therefore have a subsequence converging smoothly to a harmonic map of degree 2. By Lemma~\ref{lem:harmonic-scalar-gap}, spectral convergence, and the 3 independent coordinate kernel functions of $v_\alpha$, we have $I_\alpha(v_\alpha)=N_\alpha(v_\alpha)=3$ for all $\alpha>1$ sufficiently close to 1. Intrinsic stability of $v_\alpha$ and Proposition~\ref{prop:normal-splitting} then give Morse index $3(n-2)$ for $u_\alpha=\iota\circ v_\alpha$, where $\iota:\Sph^2\hookrightarrow\Sph^n$ is an equatorial inclusion. Since equatorial inclusion preserves $\alpha$-energy, the limit \eqref{eq:degree-two-energy-limit} gives $E_{\alpha,\kappa}(u_\alpha)<\Lambda$ for every fixed $\Lambda>2\pi\kappa+8\pi$ and all $\alpha>1$ sufficiently close to 1. These nonisometric maps $u_\alpha$ therefore prove sharpness under the stated  $\alpha$-energy upper bound.

\subsection{Organization of the paper}\label{subsec:organization}
The remaining part of this paper is organized as follows.

In Section~\ref{sec:spectral}, we develop the weighted scalar spectrum, nodal domain and multiplicity estimates, and equatorial normal splitting. Subsection~\ref{subsec:scalar-operator} introduces $J_{\alpha,u}$ from \eqref{eq:operator-intro} and defines the weighted variational spectrum through \eqref{eq:main-weighted-eigenproblem}, identifying $I_\alpha(u)$ with the number of eigenvalues below 2 and $N_\alpha(u)$ with the multiplicity of 2. Subsection~\ref{subsec:normal-splitting} derives the second variation \eqref{eq:hessian} and the Morse index and nullity formulas \eqref{eq:index-splitting} and~\eqref{eq:nullity-splitting} under equatorial inclusion. Subsection~\ref{subsec:nodal-classification} establishes nodal domain and multiplicity bounds for $\ker J_{\alpha,u}$ and proves that $I_\alpha(u)=1$ on a domain of genus 0 forces $u$ to be a diffeomorphism onto an equatorial $2$-sphere.

In Section~\ref{sec:codazzi}, we express the transported scalar form as $-\alpha B_C$ and reconstruct $\alpha$-harmonic diffeomorphisms from Codazzi tensors, completing the proofs of Theorems~\ref{thm:main-spectral} and~\ref{thm:main-convex-energy}. Subsection~\ref{subsec:scalar-transport} transports the scalar form \eqref{eq:scalar-form-intro} to the round target and proves that the conductivity tensor $A$ in \eqref{eq:tensors-intro} is divergence free. Subsection~\ref{subsec:support-signature} represents the Codazzi cofactor $C=\cof_{\bar g}A$ by a support function of a strictly convex surface and determines the signature and kernel of $B_C$ in \eqref{eq:mixed-form}. Subsection~\ref{subsec:codazzi-reconstruction} reconstructs $u$ from a positive definite Codazzi tensor $C$ satisfying $\det_{\bar g}C>\kappa^{2(\alpha-1)}$ when $(\Sph^2,G_C)$ is orientation preservingly isometric to $(\Sigma,g)$, with $G_C$ defined in \eqref{eq:GC-intro}. For the unit round domain, this isometry condition is equivalent to $K_{G_C}=1$.

In Section~\ref{sec:energy}, we establish the domain balance identity, estimates for $E_{\alpha,\kappa}$, and local and equivariant rigidity. Subsection~\ref{subsec:domain-balance-identity} derives \eqref{eq:MN} from the stress tensor and rewrites it as \eqref{eq:Omega-balance} with the nonnegative weight $\Omega_{\alpha,\kappa}(|du|_g^2)$ defined in \eqref{eq:Omega}. Subsection~\ref{subsec:mixed-form-energy} proves the upper bound \eqref{eq:energy-upper} for $\alpha$-harmonic diffeomorphisms between unit round $2$-spheres and the sharp lower bound \eqref{eq:energy-lower} for maps of degree $\pm1$. Subsection~\ref{sec:fixed} computes the Hessian of $E_{\alpha,\kappa}$ at orthogonal maps, proves that sufficiently close $\alpha$-harmonic maps are orthogonal, and records the Morse index and nullity of isometric equatorial embeddings. Subsection~\ref{subsec:twisted-equivariant} shows that the twist $\gamma$ in \eqref{eq:twisted-map} is constant and that every such $\alpha$-harmonic diffeomorphism is orthogonal.

In Section~\ref{sec:compactness-rigidity}, we prove compactness and rigidity for $\alpha$-harmonic diffeomorphisms between unit round $2$-spheres, yielding Theorem~\ref{cor:main-rigidity}. Subsection~\ref{subsec:diffeo-compactness} recalls bubble tree compactness and excludes bubbling for sequences $u_{\alpha_j}$ of $\alpha_j$-harmonic diffeomorphisms as $\alpha_j\searrow1$. Subsection~\ref{subsec:diffeo-near-one} identifies the harmonic limits of $u_{\alpha_j}$ as orthogonal maps and proves that all $\alpha$-harmonic diffeomorphisms are orthogonal for $\alpha>1$ sufficiently close to 1. Subsection~\ref{sec:full-index-rigidity} bounds the sum of the Morse index and nullity by 6 and identifies the Jacobi kernel of a hypothetical nonorthogonal $\alpha$-harmonic diffeomorphism $u$ as the space \eqref{eq:symmetry-kernel} generated by domain and target rotations. Subsection~\ref{subsec:uniform-exponent-compactness} establishes uniform gradient bounds for $\alpha$-harmonic maps with bounded $E_{\alpha,\kappa}$ and smooth compactness of $\alpha$-harmonic diffeomorphisms when $\alpha$ ranges over a compact subset of $(1,\infty)$. Subsection~\ref{subsec:exponent-global-rigidity} continues a hypothetical nonorthogonal $\alpha$-harmonic diffeomorphism to nearby exponents and proves that every $\alpha$-harmonic diffeomorphism is orthogonal for all $\alpha>1$.

In Section~\ref{sec:bubble-index}, we record the harmonic index estimates and prove the alternative in Proposition~\ref{prop:two-bubble-reduction}. Subsection~\ref{subsec:harmonic-index-gap} recalls Morse index lower semicontinuity and establishes the alternative $I_1(v)=1$ or $I_1(v)\geq3$ in \eqref{eq:harmonic-scalar-gap} for nonconstant harmonic spheres $v$. Subsection~\ref{sec:two-bubble-reduction} identifies the weak limit of the balance measures $\mu_{\alpha_j}$ defined in \eqref{eq:neck-balance-measure} and proves that, for sequences $u_{\alpha_j}$ satisfying \eqref{eq:two-bubble-budget}, after passing to a subsequence either $u_{\alpha_j}$ is eventually an isometric equatorial parametrization, or the base map is constant and exactly 2 bubbles of harmonic degree 1 form at distinct antipodal points.

In Section~\ref{sec:two-bubble-necks}, we establish the Morse index lower bound $3(n-2)$ for the configuration of 2 bubbles in Theorem~\ref{thm:two-bubble-index} and construct the degree 2 minimizers in Proposition~\ref{prop:degree-two-minimizers}, completing the proof of Theorem~\ref{thm:main-full-gap}. Subsection~\ref{subsec:neck-conservation} introduces cylinder and gnomonic coordinates for $u_{\alpha_j}$ and derives the equation \eqref{eq:gnomonic-conductivity} and flux identities \eqref{eq:gnomonic-fluxes} for its gnomonic coordinate functions $z_{\alpha_j}$. Subsection~\ref{subsec:neck-decay-hamiltonian} proves decay and regularity estimates for $z_{\alpha_j}$ on the neck and compares $\alpha_j-1$ with the scale $\delta_{\alpha_j}$ from \eqref{eq:neck-scales} through the bound \eqref{eq:epsilon-delta}. Subsection~\ref{subsec:relative-neck-expansion} derives the expansion of $z_{\alpha_j}$ to first order in \eqref{eq:first-neck-expansion} and identifies the rotational Jacobi fields $K_\xi u_{\alpha_j}$ whose leading profiles change sign as in \eqref{eq:odd-rotational-mode}. Subsection~\ref{subsec:negative-variations} constructs an $(n-2)$-dimensional space of global tangent variations along $u_{\alpha_j}$ on which $D^2E_{\alpha_j,\kappa}(u_{\alpha_j})$ is negative definite, and proves its exact Hessian orthogonality \eqref{eq:exact-core-orthogonality} to the variations supported in the bubble cores. Subsection~\ref{subsec:degree-two-sharpness} proves existence of degree 2 minimizers $v_\alpha$, establishes \eqref{eq:degree-two-spectral-data} for $v_\alpha$ and \eqref{eq:degree-two-full-index} for their equatorial inclusions $\iota\circ v_\alpha$, and completes the rigidity and sharpness arguments.

In Appendix~\ref{app:poles}, we derive the pole expansions \eqref{eq:north-pole} and~\eqref{eq:south-pole} used in the proof of Proposition~\ref{prop:equivariant}.
\\

\noindent \textbf{Disclosure of AI.} The authors designed all research questions, theorem statements, and the main ideas for proving the main theorems. AI-assisted tool was used during the preparation of some technical proof details and it was not used to generate the main research ideas or the manuscript text. The manuscript was written by the authors, who take full responsibility for the content of the manuscript, including all mathematical arguments and conclusions.

\vskip1cm

\section{The weighted scalar equation and the first spectral level}\label{sec:spectral}

In this section, we develop the scalar spectral formulation of the $\alpha$-harmonic map equation and relate it to geometric rigidity. We first identify the scalar index and nullity through a weighted eigenvalue problem in Subsection~\ref{subsec:scalar-operator}, then determine their contributions to the second variation under equatorial inclusion in Subsection~\ref{subsec:normal-splitting}. We use this formulation in Subsection~\ref{subsec:nodal-classification} to classify maps of scalar index 1 on surfaces of genus 0.

Throughout this section, $(\Sigma,g)$ is a smooth closed connected oriented surface, the target sphere has sectional curvature 1, and all maps are smooth. We retain the Laplace convention fixed in Subsection~\ref{subsec:background}.
For a vector field $V$ along $u$, $dV$ denotes its ambient Euclidean
derivative. In local domain coordinates $(x^1,x^2)$, the indices $i,j$ range from 1 to 2 and repeated coordinate indices are summed. Thus
\begin{equation}\label{eq:notation}
 |du|_g^2=g^{ij}\langle\partial_i u,\partial_j u\rangle,
 \qquad
 \langle du,dV\rangle_g=g^{ij}\langle\partial_i u,\partial_j V\rangle.
\end{equation}
The contractions use the Euclidean inner product in the ambient target space.
We denote by $\nabla^u$ the pullback of the target Levi-Civita connection
and by $\nabla^g$ the domain connection. For vector fields $X,Y$ on
$\Sigma$, the second covariant differential is
\[
 (\nabla du)(X,Y)
 =\nabla^u_X(du(Y))-du(\nabla^g_XY).
\]
Torsion freeness makes this tensor symmetric in $X,Y$. For a covariant
2-tensor $S$, its divergence is the 1-form with coordinate components
\[
 (\Div_gS)_j
 =\sum_{i,k=1}^2g^{ik}(\nabla^g_{\partial_k}S)(\partial_i,\partial_j).
\]
The covariant derivative here acts on the full tensor, including its
covariant slots. Divergence of a vector field is scalar; divergence of
an ambient-vector-valued 1-form is computed componentwise. Throughout this section, $\alpha\geq1$ and $\kappa>0$. At $\alpha=1$, criticality means harmonicity, since $E_{1,\kappa}(u)=E(u)+\frac{\kappa}{2}\operatorname{Area}_g(\Sigma)$. At the harmonic exponent, the conductivity is 1 and the additional term in the full Hessian vanishes. We use the same endpoint convention for the harmonic components of a bubble tree. Unless indicated otherwise, gradients and area measures use the domain metric $g$.

\subsection{The scalar operator and its first positive eigenvalue}\label{subsec:scalar-operator}
In this subsection, we formulate the scalar equation as a nonnegative weighted eigenvalue problem in which the coordinate functions have eigenvalue 2. This formulation identifies the scalar index with the number of eigenvalues below 2 and the scalar nullity with the multiplicity of 2, as stated in Proposition~\ref{prop:weighted-spectrum}.

For a smooth map $u:\Sigma\to\Sph^n$ and a smooth vector field $V$ along $u$
satisfying $\ip{V}{u}=0$, the first variation of $E_{\alpha,\kappa}$ is
\[
 DE_{\alpha,\kappa}(u)[V]
 =\alpha\int_\Sigma
 (\kappa+|du|_g^2)^{\alpha-1}\ip{du}{dV}_g\dd A_g.
\]
Integration by parts and the constraint $|u|^2=1$ show that $u$ is
$\alpha$-harmonic if and only if
\begin{equation}\label{eq:sphere}
 -\Div_g\bigl((\kappa+|du|_g^2)^{\alpha-1}\nabla u\bigr)
 =(\kappa+|du|_g^2)^{\alpha-1}|du|_g^2u.
\end{equation}
Equivalently,
\begin{equation}\label{eq:drift}
 \Delta_g u
 -\frac{\alpha-1}{\kappa+|du|_g^2}
 du\bigl(\nabla|du|_g^2\bigr)
 =|du|_g^2u.
\end{equation}
The coefficient of order 0 in \eqref{eq:sphere} satisfies
\begin{equation}\label{eq:linked}
 (\kappa+|du|_g^2)^{\alpha-1}|du|_g^2
 =(\kappa+|du|_g^2)^\alpha
 -\kappa(\kappa+|du|_g^2)^{\alpha-1}.
\end{equation}

Fix a smooth nonconstant $\alpha$-harmonic map $u:\Sigma\to\Sph^n$. For
real valued $f,h\in W^{1,2}(\Sigma)$, define the scalar bilinear form
\begin{equation}\label{eq:scalar-form}
 Q_{\alpha,u}(f,h)
 =\alpha\int_\Sigma
 (\kappa+|du|_g^2)^{\alpha-1}
 \bigl(\ip{\nabla f}{\nabla h}_g-|du|_g^2fh\bigr)\dd A_g.
\end{equation}
The associated self adjoint operator on
$L^2\bigl(\Sigma,(\kappa+|du|_g^2)^{\alpha-1}\dd A_g\bigr)$
is $\alpha J_{\alpha,u}$, with $J_{\alpha,u}$ as in
\eqref{eq:operator-intro}. Thus the scalar quantities $I_\alpha(u)$ and
$N_\alpha(u)$ introduced in Subsection~\ref{subsec:main-index} are,
equivalently, the index and nullity of $Q_{\alpha,u}$.

By \eqref{eq:sphere}, every Euclidean coordinate function of $u$
belongs to $\ker J_{\alpha,u}$. Moreover, nonconstancy of $u$ gives
\[
 Q_{\alpha,u}(1,1)
 =-\alpha\int_\Sigma
 (\kappa+|du|_g^2)^{\alpha-1}|du|_g^2\dd A_g<0,
\]
which implies
\begin{equation}\label{eq:basic-index}
 I_\alpha(u)\geq1,\qquad N_\alpha(u)\geq r(u)+1.
\end{equation}
The nullity lower bound follows from the linear independence of the
coordinate functions in the smallest equatorial sphere, whose dimension
$r(u)$ was fixed in Subsection~\ref{subsec:main-spectral}.

Multiplication by $(\kappa+|du|_g^2)^{(\alpha-1)/2}$ is unitary from
$L^2\bigl(\Sigma,(\kappa+|du|_g^2)^{\alpha-1}\dd A_g\bigr)$
to $L^2(\Sigma,\dd A_g)$ and conjugates $J_{\alpha,u}$ to a
Schr\"odinger operator with smooth potential. Indeed, for
$\psi\in C^\infty(\Sigma)$, a direct computation gives
\begin{equation}\label{eq:conjugation}
\begin{aligned}
 (\kappa+|du|_g^2)^{(\alpha-1)/2}
 J_{\alpha,u}\bigl((\kappa+|du|_g^2)^{(1-\alpha)/2}\psi\bigr)
 &=\Delta_g\psi-|du|_g^2\psi
 -\frac{\alpha-1}{2}\Delta_g\log(\kappa+|du|_g^2)\,\psi\\
 &\quad+\frac{(\alpha-1)^2}{4}
 |\nabla\log(\kappa+|du|_g^2)|_g^2\psi.
\end{aligned}
\end{equation}

Thus $J_{\alpha,u}$ has discrete spectrum, and its lowest eigenvalue is simple with a strictly positive eigenfunction. The conjugating factor is smooth and positive, so it preserves nodal sets and vanishing orders. Unique continuation and the local nodal set
structure for eigenfunctions of Schr\"odinger operators therefore apply to eigenfunctions of $J_{\alpha,u}$, see
\cite[Section~2 and Appendix]{Kokarev2014}.

The scalar equation \eqref{eq:sphere} also admits a nonnegative
variational spectral formulation. With the weighted energy measure
as reference, the coordinate functions of $u$ have eigenvalue $2$, while
the scalar index and nullity count the eigenvalues with multiplicity below and at $2$, respectively.

\begin{proposition}\label{prop:weighted-spectrum}
Let $u:(\Sigma,g)\to\Sph^n$ be a smooth nonconstant
$\alpha$-harmonic map, with $\alpha\geq1$ and $\kappa>0$. Use the weighted
energy measure $\nu_{\alpha,u}$ from \eqref{eq:main-weighted-eigenproblem}
and define the nonnegative bilinear form
\[
 \mathfrak d_{\alpha,u}(f,h)
 =\int_\Sigma(\kappa+|du|_g^2)^{\alpha-1}
 \ip{\nabla f}{\nabla h}_g\dd A_g,
 \qquad f,h\in W^{1,2}(\Sigma).
\]
The finite eigenvalues of the variational problem
\begin{equation}\label{eq:weighted-variational-spectrum}
 \mathfrak d_{\alpha,u}(f,h)
 =\lambda\int_\Sigma fh\dd\nu_{\alpha,u}
 \qquad\text{for all }h\in W^{1,2}(\Sigma),
\end{equation}
with $f\in W^{1,2}(\Sigma)\setminus\{0\}$, form a discrete sequence
\[
 0=\lambda_0<\lambda_1\leq\lambda_2\leq\cdots,
 \qquad \lambda_k\longrightarrow+\infty,
\]
where eigenvalues are repeated according to multiplicity.
The 0-eigenspace consists of the constant functions, and every
nonzero Euclidean coordinate function of $u$ is an eigenfunction
with eigenvalue $2$. There are exactly $I_\alpha(u)$ eigenvalues
strictly below $2$, and the multiplicity of $2$ is $N_\alpha(u)$.
\end{proposition}

\begin{proof}
Define
\[
 \mathfrak h(f,h)
 =\mathfrak d_{\alpha,u}(f,h)
 +\int_\Sigma fh\dd\nu_{\alpha,u}.
\]
Since the conductivity $(\kappa+|du|_g^2)^{\alpha-1}$ is smooth and strictly positive and $\nu_{\alpha,u}$ has bounded density and positive total mass, $\mathfrak h$ is a Hilbert inner product inducing a norm equivalent to the standard
$W^{1,2}(\Sigma)$ norm.

Let $T:W^{1,2}(\Sigma)\to W^{1,2}(\Sigma)$ be a linear operator determined by
\[
 \mathfrak h(Tf,h)=\int_\Sigma fh\dd\nu_{\alpha,u}.
\]
The operator $T$ is nonnegative and self adjoint with respect to
$\mathfrak h$, and compact by Rellich's theorem. Since
$\mathfrak h(Tf,f)\leq\mathfrak h(f,f)$, every positive eigenvalue
$t$ of $T$ lies in $(0,1]$. Moreover, the density of
$\nu_{\alpha,u}$ is positive on a nonempty open set, so $T$ has
infinite rank and its positive eigenvalues tend to 0.
For $t>0$, the equation $Tf=tf$ is equivalent to
\eqref{eq:weighted-variational-spectrum} with
\[
 \lambda=t^{-1}-1.
\]
The compact spectral theorem therefore gives the asserted
sequence of finite eigenvalues $\lambda_k \nearrow \infty$. The eigenvalue $\lambda=0$ is simple, because $\mathfrak d_{\alpha,u}(f,f)=0$ holds precisely for constant functions on the connected surface $\Sigma$.

By equation~\eqref{eq:sphere}, we see that coordinate function of $u$ has eigenvalue $2$ for \eqref{eq:weighted-variational-spectrum}.
To identify the corresponding index and nullity interpretation, observe that
\begin{equation*}
    \frac1\alpha Q_{\alpha,u}(f,h) =\mathfrak d_{\alpha,u}(f,h)
   -2\int_\Sigma fh\dd\nu_{\alpha,u} =\mathfrak h(f,h)-3\mathfrak h(Tf,h).
\end{equation*}
On the eigenspace of $T$ corresponding to $t>0$, the form
$Q_{\alpha,u}$ is $\alpha(1-3t)\mathfrak h$, and
\[
 1-3t=\frac{\lambda-2}{\lambda+1}.
\]
On $\ker T$, the form $Q_{\alpha,u}=\alpha\mathfrak h$ is
positive definite. The spectral decomposition of $T$ thus
identifies the index and nullity of $Q_{\alpha,u}$
with the number of variational eigenvalues below $2$ and the
multiplicity of $2$, respectively. This completes the proof of Proposition \ref{prop:weighted-spectrum}.
\end{proof}

\begin{remark}\label{rem:weighted-not-metric}
We call $u$ a first eigenmap for \eqref{eq:weighted-variational-spectrum} if its coordinate functions belong to the first positive eigenspace of $\mathfrak d_{\alpha,u}$.
By Proposition~\ref{prop:weighted-spectrum}, this condition is
equivalent to $\lambda_1=2$, or equivalently to $I_\alpha(u)=1$.

At $\alpha=1$, the conductivity $(\kappa+|du|_g^2)^{\alpha-1}$ is identically 1, and
\eqref{eq:weighted-variational-spectrum} is the variational
Laplace eigenvalue problem for the energy metric $\frac12|du|_g^2g$, which may degenerate at finitely many points. Indeed, in a local conformal coordinate the harmonic map equation gives $\nabla^u_{\partial_{\bar z}}(du(\partial_z))=0$, with the pullback connection complexified, so $\partial u$ is a holomorphic $(1,0)$-form with values in the complexified pullback tangent bundle, equipped with the holomorphic structure induced by its pullback connection. Since $u$ is nonconstant, the zeros of this holomorphic section are isolated and hence finite on $\Sigma$. For a nonconstant conductivity, however, the weighted problem in \eqref{eq:weighted-variational-spectrum} is not in general the Laplace
eigenvalue problem for the energy metric. On a surface, the ordinary Dirichlet form is conformally invariant, so the
nonconstant factor $(\kappa+|du|_g^2)^{\alpha-1}$ in
$\mathfrak d_{\alpha,u}$ cannot be absorbed into a conformal change of $g$.
\end{remark}

\subsection{Second variation and equatorial normal splitting}\label{subsec:normal-splitting}
In this subsection, we compute the second variation of $E_{\alpha,\kappa}$ and isolate the contribution from directions normal to an equatorial target sphere. Under equatorial inclusion, each independent normal direction contributes a copy of the scalar form, giving the Morse index and nullity identities in Proposition~\ref{prop:normal-splitting}.

Let $u:\Sigma\to\Sph^n$ be a smooth $\alpha$-harmonic map. For a smooth
section $V$ of $u^{-1}T\Sph^n$, the ambient derivative $dV$ used in
\eqref{eq:notation} satisfies, since $\langle V,u\rangle=0$,
\[
 dV(X)=\nabla^u_XV-\langle V,du(X)\rangle u.
\]
Thus $dV$ is distinguished from the intrinsic derivative $\nabla^uV$.
The standard second variation formula of $E_{\alpha,\kappa}$ 
(cf.\ Micallef-Moore \cite{MicallefMoore1988}) reads
\begin{equation}\label{eq:hessian}
\begin{split}
 D^2E_{\alpha,\kappa}(u)[V,V]
 ={}&\alpha\int_\Sigma
 (\kappa+|du|_g^2)^{\alpha-1}
 \bigl(|dV|_g^2-|du|_g^2|V|^2\bigr)\dd A_g\\
 &+2\alpha(\alpha-1)\int_\Sigma
 (\kappa+|du|_g^2)^{\alpha-2}
 \ip{du}{dV}_g^2\dd A_g.
\end{split}
\end{equation}
The Hessian of $E_{\alpha,\kappa}$ extends continuously to $W^{1,2}$ sections of $u^{-1}T \S^n$. For the target $\Sph^m$, we use the Morse index and nullity notation
introduced in Subsection~\ref{subsec:main-index}. Recall that the Morse index is the maximal dimension of a subspace of
$W^{1,2}(\Sigma,u^{-1}T\Sph^m)$ on which the Hessian $D^2E_{\alpha,\kappa}(u)$ is negative
definite, and the nullity is the dimension of the kernel of the
associated Jacobi operator. These definitions agree with those using smooth sections by density of Sobolev space $W^{1,2}(\Sigma,u^{-1}T\Sph^m)$ and elliptic regularity for the associated Jacobi operator.

For $\alpha$-harmonic maps taking values in an equatorial subsphere of $\S^n$, the second variation of $E_{\alpha,\kappa}$ splits into the Hessian for the smaller target and copies of the scalar form $Q_{\alpha,u}$, that is, the following holds:

\begin{proposition}\label{prop:normal-splitting}
Let $u:\Sigma\to\Sph^r\subset\Sph^m$ be $\alpha$-harmonic,
where the inclusion is equatorial. Then
\begin{align}
 \ind_{\Sph^m}E_{\alpha,\kappa}(u)
 &=\ind_{\Sph^r}E_{\alpha,\kappa}(u)
   +(m-r)I_\alpha(u),
 \label{eq:index-splitting}\\
 \nul_{\Sph^m}E_{\alpha,\kappa}(u)
 &=\nul_{\Sph^r}E_{\alpha,\kappa}(u)
   +(m-r)N_\alpha(u).
 \label{eq:nullity-splitting}
\end{align}
Both identities hold at $\alpha=1$, with $E_{\alpha,\kappa}$
replaced by the Dirichlet energy.
\end{proposition}

\begin{proof}
After an orthogonal change of ambient coordinates, assume that
$\Sph^r=\Sph^m\cap\R^{r+1}$. Choose a constant orthonormal basis
$\nu_1,\ldots,\nu_{m-r}$ of the orthogonal complement of
$\R^{r+1}$ in $\R^{m+1}$. The pullback tangent bundle has the
orthogonal decomposition
\[
 u^{-1}T\Sph^m
 =u^{-1}T\Sph^r
 \oplus\bigoplus_{j=1}^{m-r}(\Sigma\times\R\nu_j).
\]
Since $\nu_j$ is constant and perpendicular to the image of $du$, for any $f \in W^{1,2}(\Sigma)$ there holds
\[
 d(f\nu_j)=df\otimes\nu_j,
 \qquad
 \ip{du}{d(f\nu_j)}_g=0.
\]
Thus polarization of \eqref{eq:hessian} gives
\[
 D^2E_{\alpha,\kappa}(u)[f\nu_j,h\nu_j]
 =Q_{\alpha,u}(f,h),
 \qquad f,h\in W^{1,2}(\Sigma).
\]

If $V$ is a section of $u^{-1}T\Sph^r$, then both $V$ and its
ambient derivative $dV$ take values in $\R^{r+1}$. Hence
\[
 \ip{V}{f\nu_j}=0,\qquad
 \ip{dV}{d(f\nu_j)}_g=0,
\]
and the mixed Hessian terms between sections in
$u^{-1}T\Sph^r$ and each normal summand $f \nu_j$ vanish. Orthogonality of
the vectors $\nu_j$ likewise eliminates the mixed terms between
distinct normal summands.

Consequently, the Hessian $D^2 E_{\alpha,\kappa}$  is the orthogonal direct sum of its restriction to $W^{1,2}(\Sigma,u^{-1}T\Sph^r)$ and $m-r$ copies of $Q_{\alpha,u}$. The first restriction is precisely the Hessian
of $E_{\alpha,\kappa}$ acting on $\alpha$-harmonic maps with target $\Sph^r$. Additivity of the index and nullity under this direct sum proves
\eqref{eq:index-splitting} and \eqref{eq:nullity-splitting}. At $\alpha=1$, the second integral in \eqref{eq:hessian}
vanishes, and $E_{1,\kappa}$ differs from the Dirichlet energy
by a constant independent of the map which does not affect their Hessian. Thus, we complete the proof of Proposition \ref{prop:normal-splitting}.
\end{proof}

\subsection{Nodal domains and the classification in genus 0}\label{subsec:nodal-classification} 
In this subsection, we relate the nodal domains of scalar kernel functions to the scalar index, obtaining the bound in Lemma~\ref{lem:nodal-count}. On a surface of genus 0, this control shows that a nonconstant $\alpha$-harmonic map of scalar index 1 is a diffeomorphism onto an equatorial $2$-sphere, as proved in Proposition~\ref{prop:index-one-forward}. Throughout this subsection, $u:(\Sigma,g)\to\Sph^n$ is a smooth nonconstant $\alpha$-harmonic map, where $\Sigma$ is closed, connected and oriented, $\alpha\geq1$, and $\kappa>0$. 

For a smooth function $f$ on $\Sigma$, its nodal set is $f^{-1}(0)$, and its
nodal domains are the connected components of
$\Sigma\setminus f^{-1}(0)$. The function has a fixed sign on each nodal domain.  Near a regular zero point of $f$, the implicit function theorem identifies the nodal set with a smooth embedded curve.

We recall the corresponding description near a critical zero point for solutions of a scalar Schr\"odinger equation with smooth real potential, see \cite[Section~2 and Appendix]{Kokarev2014}. Strong unique continuation implies that a nontrivial solution has finite vanishing order $q\geq1$ at every zero point.  In local conformal coordinates, the first nonzero homogeneous term in the Taylor expansion is a harmonic polynomial of degree $q$. Cheng's \cite{cheng1976eigenfunctions} local nodal set theorem shows that the nodal set consists of exactly $2q$ incident arcs in a sufficiently small neighborhood of the zero point. These arcs have equally spaced tangent directions, and the sign of the solution alternates between adjacent sectors. In particular, a critical zero point has $q\geq2$ and is isolated. Unique continuation also implies that a solution vanishing on a nonempty open set of $\Sigma$ vanishes identically.

By \eqref{eq:conjugation}, multiplication by
$(\kappa+|du|_g^2)^{(\alpha-1)/2}$ transforms every element of $\ker J_{\alpha,u}$ into a solution of a Schr\"odinger equation with smooth potential. The multiplying factor is smooth and strictly positive, so it preserves signs, nodal domains, critical zero points, and their vanishing orders. The preceding local description and unique continuation therefore apply to functions in $\ker J_{\alpha,u}$.

For a nontrivial solution of a Schr\"odinger equation with smooth potential on a closed surface, the critical zero points are finite. Removing small neighborhoods of these points leaves a compact 1-dimensional nodal set with boundary and finitely many components. Thus the full nodal set consists of a finite embedded graph, whose vertices are the critical zeros and whose edges are the intervening nodal arcs, together with finitely many smooth closed curve components. Its complementary
components are precisely the nodal domains.

With the genus notation $\mathfrak g$ from
Subsection~\ref{subsec:main-spectral}, Besson's multiplicity bound \cite[Theorem~2.1]{Besson1980}, in the smooth Schr\"odinger formulation recalled in \cite[Theorem~1.1]{Kokarev2014}, bounds the multiplicity of the $k$-th eigenvalue for a Schr\"odinger operator by $2(2-\chi(\Sigma))+2k+1$, where eigenvalues are numbered from 0 and repeated according to multiplicity. The conjugation
\eqref{eq:conjugation} is unitary and therefore preserves eigenvalues and their multiplicities. Since exactly $I_\alpha(u)$ many eigenvalues of $J_{\alpha,u}$ are negative, the first occurrence of the eigenvalue~$0$ has index $k=I_\alpha(u)$. Using $\chi(\Sigma)=2-2\mathfrak g$ and the kernel bound
\eqref{eq:basic-index}, we obtain
\begin{equation}\label{eq:general-multiplicity}
 N_\alpha(u)\leq2I_\alpha(u)+4\mathfrak g+1,
 \qquad r(u)\leq2I_\alpha(u)+4\mathfrak g.
\end{equation}
In particular, when $\Sigma$ has genus 0, we have 
\begin{equation}\label{eq:multiplicity}
 N_\alpha(u)\leq2I_\alpha(u)+1,
 \qquad
 r(u)\leq2I_\alpha(u).
\end{equation}

Courant's nodal domain bound takes the following form at the eigenvalue~$0$, see \cite[Section~2]{Kokarev2014}. We give the detailed argument directly for the scalar form $Q_{\alpha,u}$.

\begin{lemma}\label{lem:nodal-count}
Every nonzero element of $\ker J_{\alpha,u}$ has at most
$I_\alpha(u)+1$ nodal domains.
\end{lemma}

\begin{proof}
Fix $f\in\ker J_{\alpha,u}\setminus\{0\}$, and let
$\Omega_1,\ldots,\Omega_s$ be its nodal domains. For
$1\leq j\leq s$, define
\[
 f_j(x)=
 \begin{cases}
  f(x),&x\in\Omega_j,\\
  0,&x\notin\Omega_j.
 \end{cases}
\]
We first justify using $f_j$ as a test function. If $f>0$ on
$\Omega_j$, extend
$\max\{f-\delta,0\}|_{\Omega_j}$ by 0 to $\Sigma$.
For every $\delta>0$, this function is supported in a compact
subset of $\Omega_j$ and belongs to $W^{1,2}(\Sigma)$. As $\delta\searrow0$, dominated convergence theorem for these functions and their weak gradients gives convergence in $W^{1,2}(\Sigma)$ to $f_j$. On a negative nodal domain, apply the same construction to $-f$ and then change the sign. Thus, we have that 
\[
 f_j\in W^{1,2}(\Sigma),
 \qquad
 \nabla f_j=\mathbf{1}_{\Omega_j}\nabla f
 \quad\text{almost everywhere on }\Sigma.
\]
Testing $J_{\alpha,u}f=0$ against $f_j$ gives
\[
 Q_{\alpha,u}(f_j,f_j)=Q_{\alpha,u}(f,f_j)=0.
\]
For $j\neq k$, disjointness of $\Omega_j$ and $\Omega_k$
implies that both $f_jf_k$ and
$\ip{\nabla f_j}{\nabla f_k}_g$ vanish almost everywhere, which implies
\[
 Q_{\alpha,u}(f_j,f_k)=0
 \qquad  1\leq j,k\leq s.
\]

Suppose that $s>I_\alpha(u)+1$. Choose coefficients
$c_1,\ldots,c_s$, not all 0, such that $c_1=0$ and
\[
 F=\sum_{j=1}^s c_jf_j
\]
is orthogonal to every negative eigenfunction of $J_{\alpha,u}$
in
$L^2\bigl(\Sigma,(\kappa+|du|_g^2)^{\alpha-1}\dd A_g\bigr)$.
Such a choice is possible because the coefficient constraint $c_1 = 0$ and the orthogonality conditions impose at most $I_\alpha(u)+1$ many homogeneous linear conditions on coefficients $c_j$. The functions $f_j$ are linearly independent,
so $F\neq0$, and the preceding identities give $Q_{\alpha,u}(F,F)=0$.

The spectral decomposition of $J_{\alpha,u}$ shows that
$Q_{\alpha,u}$ is nonnegative on the orthogonal complement
of its negative eigenspaces. We verify that the equality $Q_{\alpha,u}(F,F)=0$ forces $F$ to be a kernel function of $J_{\alpha, u}$. For every $h\in W^{1,2}(\Sigma)$ in the same orthogonal complement and every $t\in\R$, we have
\[
 0 \leq Q_{\alpha,u}(F+th,F+th) =2tQ_{\alpha,u}(F,h)+t^2Q_{\alpha,u}(h,h).
\]
Taking $t$ of either sign and letting $t\to0$ gives $Q_{\alpha,u}(F,h)=0$. Pairing $F$ with every negative eigenfunction using $Q_{\alpha,u}$ also vanishes, by its eigenvalue equation and the weighted $L^2$ orthogonality of $F$. Decomposing an arbitrary $W^{1,2}$ test function into its negative spectral part
and its orthogonal complement in the weighted $L^2$ space, we obtain
\[
 Q_{\alpha,u}(F,h)=0
 \qquad\text{for every }h\in W^{1,2}(\Sigma).
\]
Thus $F$ is a weak solution of $J_{\alpha,u}F=0$, and elliptic
regularity gives $F\in\ker J_{\alpha,u}$.

Since $c_1=0$, the function $F$ vanishes on the nonempty open
set $\Omega_1$. Unique continuation gives $F\equiv0$,
contrary to its construction. This proves Lemma~\ref{lem:nodal-count}.
\end{proof}

The nodal domain bound holds on every closed surface. The next lemma uses the genus 0 assumption to pass from
the local structure at a critical zero point to a global nodal domain count.

\begin{lemma}\label{lem:nodal-critical}
Let $(\Sigma,g)$ be a closed connected surface of genus 0 with smooth metric $g$. A nontrivial solution of a scalar Schr\"odinger equation with smooth real potential and exactly 2 nodal domains has no critical zero. The same conclusion holds for every nonzero element of
$\ker J_{\alpha,u}$ with exactly 2 nodal domains.
\end{lemma}

\begin{proof}
Suppose that a critical zero point exists, and consider the connected component of the nodal graph containing that zero point. Take the critical zero points as vertices and the intervening nodal arcs as edges. An edge may start and end at the same vertex, such an edge is counted twice in the degree of that vertex. Here the degree is the number of incident nodal arcs.

Let $V$, $E$, and $F$ denote the numbers of vertices, edges, and connected components of the complement of this graph, respectively. Every vertex is a zero of order at least 2
and therefore has degree at least 4. Counting the 2 ends of each edge gives $ V\geq1$ and $E\geq2V.$
Since $\Sigma$ is topologically a sphere and the selected graph is connected, Euler's formula gives $ V-E+F=2$, which means $F=2-V+E\geq V+2\geq3$.  The full nodal set has empty interior, so each of these $F$ complementary regions contains at least 1 nodal domain. Nodal domains in different regions cannot be joined without crossing the selected graph. Equivalently, adding the remaining components of the nodal set can subdivide the complementary regions but cannot merge them. The solution therefore has at least 3 nodal domains, a contradiction to the assumption of Lemma \ref{lem:nodal-critical}. For $J_{\alpha,u}$, the conclusion follows from \eqref{eq:conjugation}, since the positive conjugating factor
preserves nodal domains and critical zeros. This proves Lemma~\ref{lem:nodal-critical}.
\end{proof}

We now apply the preceding lemmas to the height functions of the map so as to obtain the following.

\begin{proposition}\label{prop:index-one-forward}
Let $u:(\Sigma,g)\to\Sph^n$ be a smooth nonconstant
$\alpha$-harmonic map, where $\Sigma$ has genus 0 and $g$
is any smooth metric. If $I_\alpha(u)=1$, then $r(u)=2$,
and $u$ is a smooth diffeomorphism onto the smallest equatorial
2-sphere containing its image. Moreover,
\[
 N_\alpha(u)=3,
 \qquad
 \ker J_{\alpha,u}=\Span\{u^1,u^2,u^3\},
\]
where $u^1,u^2,u^3$ are the coordinate functions in the
smallest equatorial target.
\end{proposition}

\begin{proof}
By \eqref{eq:multiplicity}, we have $r(u)\leq2$. Nonconstancy of $u$ implies $r(u)\geq1$. We first exclude the case $r(u)=1$. In that case, the map $u$ takes values in a great circle of $\S^n$. Since $\Sigma$ is simply connected, the map
has a smooth real valued lift $u=(\cos\theta,\sin\theta)$. The tangential component of \eqref{eq:sphere} gives
\[
 -\Div_g\bigl((\kappa+|du|_g^2)^{\alpha-1}
 \nabla\theta\bigr)=0.
\]
Multiplying by $\theta$ and integrating by parts, we obtain
\[
 \int_\Sigma
 (\kappa+|du|_g^2)^{\alpha-1}|\nabla\theta|_g^2
 \dd A_g=0,
\]
which means $\theta$ and hence $u$
are constant, a contradiction. Thus, we have $r(u)=2$.

Identify the smallest equatorial target with
$\Sph^2\subset\R^3$. The coordinate functions
$u^1,u^2,u^3$ are linearly independent, since a nontrivial
linear relation would place the image in a smaller equatorial
sphere. Combining \eqref{eq:basic-index} and
\eqref{eq:multiplicity}, we have $N_\alpha(u)=3$, and the 3 coordinate functions of $u$ span $\ker J_{\alpha,u}$.

For every $b\in\R^3\setminus\{0\}$, the height function $f_b=\ip{b}{u}$ is a nonzero element of $\ker J_{\alpha,u}$.
Let $h>0$ be an eigenfunction corresponding to the lowest
eigenvalue of $J_{\alpha,u}$, where the corresponding eigenvalue is negative, whereas $J_{\alpha,u}f_b=0$. Self adjointness of $J_{\alpha,u}$ therefore gives
\[
 \int_\Sigma
 (\kappa+|du|_g^2)^{\alpha-1}f_bh\dd A_g=0.
\]
Since both the conductivity $(\kappa+|du|_g^2)^{\alpha-1}$ and $h$ are strictly positive, $f_b$ takes both positive and negative values. It has at least 2
nodal domains and, by Lemma~\ref{lem:nodal-count}, at most
$I_\alpha(u)+1=2$ nodal domains. Thus every nonzero height function has
exactly 2 nodal domains and, by
Lemma~\ref{lem:nodal-critical}, $f_b$ has no critical zero points.

Suppose that $du$ has rank less than 2 at some $p\in\Sigma$. We can choose a nonzero vector $b\in T_{u(p)}\Sph^2\subset\R^3$ perpendicular to
$du_p(T_p\Sigma)$. The corresponding height function satisfies
\begin{equation*}
    f_b(p)=\ip{b}{u(p)}=0, \qquad df_b(p)[X]=\ip{b}{du_p(X)}=0
 \quad\text{for every }X\in T_p\Sigma.
\end{equation*}
Thus $p$ is a critical zero point of the nontrivial kernel function $f_b$, a contradiction. Hence $du$ has rank 2 everywhere, and the inverse function theorem makes $u$ a local diffeomorphism. The image $u(\Sigma)$ is open because $u$ is a local diffeomorphism, and closed because $\Sigma$ is compact. Connectedness of the target gives $u(\Sigma)=\Sph^2$. Moreover, a local diffeomorphism from a compact manifold is proper, so $u$ is a covering map. Since the domain $\Sigma$ is connected and the target $\Sph^2$ is simply connected, this covering has 1 sheet. Therefore $u$ is a smooth diffeomorphism, which completes the proof of Proposition~\ref{prop:index-one-forward}. 
\end{proof}

The converse of Proposition~\ref{prop:index-one-forward}, that is,  any smooth $\alpha$-harmonic diffeomorphism on a sphere must admit $I_{\alpha}(u) = 1$ and $N_{\alpha}(u) = 3$, will be proved in the next section by transporting the scalar quadratic form to the equatorial target and applying the tools of convex geometry.

\vspace{1cm}

\section{Convex realization and reconstruction}\label{sec:codazzi}

In this section, we give a convex geometric description of the weighted scalar form associated with an $\alpha$-harmonic diffeomorphism. We first transport the form to the round target in Subsection~\ref{subsec:scalar-transport}, then realize its Codazzi tensor by a strictly convex embedding and determine its signature in Subsection~\ref{subsec:support-signature}. We finally reconstruct the map from the tensor in Subsection~\ref{subsec:codazzi-reconstruction}.

Let $u:(\Sigma,g)\to(\Sph^2,\bar g)$ be a smooth orientation preserving $\alpha$-harmonic diffeomorphism, where $\alpha>1$ and $\kappa>0$. Although $\Sigma$ necessarily has genus 0, no curvature assumption is imposed on $g$ until the curvature 1 reconstruction at the end of this section. We retain the Jacobian $J>0$ from \eqref{eq:tensors-intro}. Unless otherwise stated, gradients, covariant derivatives, traces, determinants, and area measures are taken with respect to the round target metric $\bar g$. Functions defined on the domain appearing in target formulas are understood to be composed with $u^{-1}$.

\subsection{Transport of the scalar form}\label{subsec:scalar-transport}

In this subsection, we transport the scalar form to the round target and express it in terms of a positive conductivity tensor. The $\alpha$-harmonic map equation makes this tensor divergence free, giving the representation in Proposition~\ref{prop:target-conductivity} and the Codazzi description used below.

Consider the positive definite symmetric endomorphism $A$ introduced in
\eqref{eq:tensors-intro}. At $y=u(x)$, its adjoint factor $du_x^*:T_y\Sph^2\to T_x\Sigma$ is the adjoint of $du_x$ with
respect to $g$ and $\bar g$, characterized by
\[
 g(du_x^*V,X)=\bar g(V,du_xX),
 \qquad V\in T_y\Sph^2,\quad X\in T_x\Sigma.
\]
Thus $du_xdu_x^*$ acts on the target tangent plane. We use $\bar g$ to identify symmetric endomorphisms with symmetric covariant 2-tensors. Under this identification, $A$ has the bilinear expression
\[
 A_y(V,W)
 =\frac{\bigl(\kappa+|du_x|_g^2\bigr)^{\alpha-1}}{J(x)}
   g(du_x^*V,du_x^*W),
 \qquad V,W\in T_y\Sph^2.
\]
The tensor $A$ is the coefficient obtained by transporting the weighted Dirichlet scalar form to the target sphere, that is, the composition $du_xdu_x^*$ accounts for the transformation of the gradients of the scalar test functions, and the factor $J(x)^{-1}$ accounts for the change of area measure.

\begin{proposition}\label{prop:target-conductivity}
The tensor $A$ is divergence free. Moreover, for every $f,h\in W^{1,2}(\Sph^2)$, we have
\begin{equation}\label{eq:pushforward}
 \frac1\alpha Q_{\alpha,u}(f\circ u,h\circ u)
 =\int_{\Sph^2}
   \bigl(A(\nabla f,\nabla h)-(\tr A)fh\bigr)\dd A_{\bar g}.
\end{equation}
\end{proposition}

\begin{proof}
For $f,h\in C^\infty(\Sph^2)$, the chain rule gives
\[
 \nabla_g(f\circ u)(x)=du_x^*\nabla f(y),
 \qquad y=u(x),
\]
and hence
\[
 \ip{\nabla_g(f\circ u)}
    {\nabla_g(h\circ u)}_g(x)
 =\bar g\bigl(du_xdu_x^*\nabla f(y),\nabla h(y)\bigr).
\]
Furthermore, note that 
\[
 \tr A(y)
 =\frac{\bigl(\kappa+|du_x|_g^2\bigr)^{\alpha-1}
              |du_x|_g^2}{J(x)}.
\]
Substituting these identities into the definition of $Q_{\alpha,u}$
and changing variables by $u$ proves \eqref{eq:pushforward} for smooth
functions. Composition with $u$ is an isomorphism from
$W^{1,2}(\Sph^2)$ to $W^{1,2}(\Sigma)$, and both bilinear forms are continuous.
The identity therefore extends to all $f,h\in W^{1,2}(\Sph^2)$ by density of Sobolev spaces.

Next, we show $A$ is divergence free. Let $y^i$, $1\leq i\leq3$, denote the restrictions of the ambient Euclidean coordinate functions to $\Sph^2$, and write $u^i=y^i\circ u$. The $\alpha$-harmonic map equation gives
\[
 Q_{\alpha,u}(u^i,f\circ u)=0
 \qquad\text{for every }f\in C^\infty(\Sph^2).
\]
Applying \eqref{eq:pushforward}, we obtain
\[
 \int_{\Sph^2}A(\nabla y^i,\nabla f)\dd A_{\bar g}
 =\int_{\Sph^2}(\tr A)y^if\dd A_{\bar g},
\]
which implies
\[
 -\Div(A\nabla y^i)=(\tr A)y^i.
\]
On the unit round sphere, $\Hess y^i=-y^i\bar g$, so the product rule yields
\[
 -\Div(A\nabla y^i)
 =-(\Div A)(\nabla y^i)+(\tr A)y^i,
\]
where the divergence of the covariant tensor $A$ is regarded as a 1-form. Comparing the 2 expressions gives
$(\Div A)(\nabla y^i)=0$ for $1\leq i\leq3$. Since the 3 coordinate
gradients span the tangent plane at every point of $\Sph^2$, it follows
that $\Div A=0$, which completes the proof of Proposition \ref{prop:target-conductivity}.
\end{proof}

We use the cofactor convention and the tensor $C=\cof A$ from
\eqref{eq:tensors-intro}. 
Metric compatibility gives $\Div(\cof C)=d(\tr C)-\Div C$. More precisely, in a local $\bar g$-orthonormal frame $e_1,e_2$,
write $\nabla_iC_{jk}=(\nabla_{e_i}C)(e_j,e_k)$ for the full
covariant-tensor derivative. The components are
\begin{equation*}
    \bigl(\Div(\cof C)\bigr)_1
   =\nabla_1C_{22}-\nabla_2C_{12}, \qquad \bigl(\Div(\cof C)\bigr)_2
   =\nabla_2C_{11}-\nabla_1C_{12}.
\end{equation*}
The vanishing of these 2 components is precisely the pair of
independent Codazzi equations in dimension 2, that is, 
\begin{equation}\label{eq:codazzi-equivalence}
 \Div(\cof C)=0
 \quad\Longleftrightarrow\quad
 \nabla_iC_{jk}=\nabla_jC_{ik}.
\end{equation}
Since $A$ is positive definite, $C=\cof A$ is also positive definite.
Moreover, $\cof C=A$ and Proposition~\ref{prop:target-conductivity}
gives $\Div A=0$. Hence $C$ is a positive definite Codazzi tensor.

\subsection{Support functions and the scalar signature}\label{subsec:support-signature}
In this subsection, we represent positive Codazzi tensors by support functions of strictly convex embeddings, as in Lemma~\ref{lem:support}. This realization determines the signature and kernel of the associated bilinear form in Proposition~\ref{prop:signature}, and consequently gives scalar index 1 and nullity 3 for $\alpha$-harmonic diffeomorphisms in Proposition~\ref{prop:automatic-index}.

For a smooth real valued function $h$ on the unit round sphere, define
\begin{equation}\label{eq:support-operator}
 \Aop[h]=\Hess h+h\bar g,
 \qquad
 \Hone=\{y\mapsto b\cdot y:b\in\R^3\}.
\end{equation}
Commuting covariant derivatives on the unit round sphere shows that
$\Aop[h]$ is Codazzi, while the nonnegative Laplace convention gives
$\tr\Aop[h]=(2-\Delta)h$.

\begin{lemma}\label{lem:support}
Every smooth symmetric Codazzi tensor $C$ on $\Sph^2$ has the form
$C=\Aop[h]$, where $h$ is unique modulo $\Hone$. If $C$ is positive
definite, then $h$ can be chosen strictly positive, and the map $X$ in
\eqref{eq:convex-X} is a smooth strictly convex embedding in $\R^3$ whose outward Gauss
map is $y$ and whose principal radii are the eigenvalues of $C$.
\end{lemma}

\begin{proof}
Contracting the Codazzi equation \eqref{eq:codazzi-equivalence} gives $\Div C=d(\tr C)$. For
$\ell\in\Hone$, the identities $\Hess\ell=-\ell \bar g$ and
$\Delta\ell=2\ell$ therefore yield
\begin{align*}
 0=\int_{\Sph^2}\Div(C\nabla\ell)\dd A_{\bar g}=\int_{\Sph^2}
    \bigl(\ip{d(\tr C)}{d\ell}-(\tr C)\ell\bigr)\dd A_{\bar g}
   =\int_{\Sph^2}(\tr C)\ell\dd A_{\bar g}.
\end{align*}
Thus $\tr C$ is orthogonal to $\Hone=\ker(2-\Delta)$, and the Fredholm alternative provides a smooth solution of $(2-\Delta)h=\tr C$.
The tensor $D=C-\Aop[h]$ is consequently trace free and Codazzi.

In a local conformal coordinate $z=x+iy$ of $\S^2$, write
\[
 D=D_{11}(dx^2-dy^2)+2D_{12}\,dx\,dy.
\]
Since $D_{22}=-D_{11}$, the terms involving the conformal factor
cancel in the Codazzi equations, leaving
\[
 \partial_xD_{11}=-\partial_yD_{12},
 \qquad
 \partial_yD_{11}=\partial_xD_{12}.
\]
Equivalently,
\[
 D_{zz}:=\frac12(D_{11}-iD_{12}),
 \qquad
 \partial_{\bar z}D_{zz}=0,
\]
so the $(2,0)$ part of $D$ is a holomorphic quadratic differential, which vanishes identically on $\Sph^2$. This proves $C=\Aop[h]$.

For the map $X$ defined in \eqref{eq:convex-X}, differentiation in $\R^3$
gives
\[
 dX_y(V)=\nabla_V\nabla h+h(y)V,
 \qquad V\in T_y\Sph^2.
\]
In particular, $dX_y(V)\in T_y\Sph^2$ and
\[
 \bar g(dX_y(V),W)=\Aop[h]_y(V,W),
 \qquad V,W\in T_y\Sph^2.
\]
Thus $dX=\Aop[h]$ under the identification of symmetric tensors
with endomorphisms through $\bar g$. If $\Aop[h]=0$, then $X$ is constant,
say $X=b$, and taking the scalar product with $y$ gives
$h(y)=b\cdot y$. Conversely, $\Aop[\ell]=0$ for every
$\ell\in\Hone$, which proves uniqueness modulo $\Hone$.

Assume now that $C>0$. The identity $dX=C$ makes $X$ an immersion with unit normal $y$. To establish global convexity of $X$, fix distinct points $y_0,y\in\Sph^2$ and a unit speed minimizing great circle arc
$\gamma:[0,t_0]\to\Sph^2$ joining them, where $0<t_0\leq\pi$.
The function
\[
 k(t)=h(\gamma(t))-X(y_0)\cdot\gamma(t)
\]
satisfies
\[
 k(0)=k'(0)=0,
 \qquad
 k''+k=C(\gamma',\gamma').
\]
Solving this initial value problem gives
\begin{equation}\label{eq:support-positive}
 h(y)-X(y_0)\cdot y
 =\int_0^{t_0}
   \sin(t_0-s)C(\gamma'(s),\gamma'(s))\dd s
 >0.
\end{equation}

Since $X(y)\cdot y=h(y)$, inequality
\eqref{eq:support-positive} shows that $X$ is injective and hence,
by compactness, a smooth embedding. Let $K$ be the convex hull of
$X(\Sph^2)$. For each $y\in\Sph^2$, the linear functional
$z\mapsto z\cdot y$ attains its maximum over $K$ uniquely at $X(y)$ by \eqref{eq:support-positive}.
In particular, $K$ has nonempty interior, since otherwise some nonzero linear functional would be constant on $K$, contradicting the uniqueness of its maximizer. Every boundary point of $K$ lies on a supporting plane, and each supporting plane intersects $K$ in the single point $X(y)$ corresponding to its outward unit normal.
It follows that $\partial K=X(\Sph^2)$ and that $K$ is strictly
convex. The outward Gauss map $N$ satisfies $N\circ X=\mathrm{Id}$,
so $dN$ is represented by $C^{-1}$ under the Gauss map
identification. Consequently, the eigenvalues of $C$ are the
positive principal radii of curvature of $X(\S^2)$.

Finally, set
\[
 b=\frac1{4\pi}\int_{\Sph^2}X(y_0)\dd A_{\bar g}(y_0).
\]
For each fixed $y$, averaging \eqref{eq:support-positive} over
$y_0$ gives $h(y)-b\cdot y>0$, since the integrand is strictly
positive except at $y_0=y$. Since $\nabla(b\cdot y)=b-(b\cdot y)y$ on $\Sph^2$, replacing $h$ by $h-b\cdot y$ translates the corresponding
embedding $X$ by $-b$. This replacement leaves $\Aop[h]=C$ unchanged and yields a strictly positive representative of the support function $h$, which completes the proof of Lemma \ref{lem:support}.
\end{proof}

For a smooth positive definite Codazzi tensor $C$, the form $B_C$ from
Theorem~\ref{thm:main-convex-energy} has the following representation
for $f,h\in C^\infty(\Sph^2)$:
\begin{equation}\label{eq:mixed-form}
  B_C(f,h) =\int_{\Sph^2}f\,\cof C:\Aop[h]\dd A_{\bar g} =\int_{\Sph^2}
    \bigl((\tr C)fh-\cof C(\nabla f,\nabla h)\bigr)\dd A_{\bar g}.
\end{equation}
Here the colon ``:" denotes full contraction with respect to $\bar g$.
The second identity in \eqref{eq:mixed-form} follows by integration by parts, using
$\Div(\cof C)=0$ and $\tr(\cof C)=\tr C$. Its right-hand side
extends $B_C$ to a continuous symmetric bilinear form on
$W^{1,2}(\Sph^2)$. Throughout the following, $\ker B_C$ consists of those
$f\in W^{1,2}(\Sph^2)$ for which $B_C(f,w)=0$ for every
$w\in W^{1,2}(\Sph^2)$.

Let $K\subset\R^3$ be the convex body whose support function $s$ satisfies $C=\Aop[s]$, as in Lemma~\ref{lem:support}. With the standard normalization $V(K,K,K)=\Vol(K)$, the mixed volume formula \cite[Corollary~2.5]{ShenfeldVanHandel2019} gives
\[
 B_C(h_L,h_M)=2\int_{\Sph^2}h_L\dd S(K,M;\cdot)=6V(K,L,M).
\]
Here $h_L,h_M$ are the smooth support functions of convex bodies $L,M$, and $S(K,M;\cdot)$ is the mixed area measure of $K$ and $M$. Thus $B_C$ is the mixed-volume bilinear form associated with $K$. Equivalently, for smooth variations it is the second variation of the enclosed volume with respect to the support function at $s$.

The following proposition is the smooth 2-dimensional
Alexandrov-Fenchel signature theorem. We give the classical kernel
and homotopy argument, as recalled in
\cite[Sections 6.1-6.2]{ShenfeldVanHandel2019}, in the present
normalization, since both the sign and the precise kernel of
$B_C$ will be used below.

\begin{proposition}\label{prop:signature}
Let $C>0$ be a smooth Codazzi tensor on $\Sph^2$. The bilinear form $B_C$ on
$W^{1,2}(\Sph^2)$ has kernel $\Hone$, and the maximal dimension of a
subspace on which $B_C$ is positive definite is 1. Moreover, we have the following
\begin{equation}\label{eq:AF}
 B_C(f,f)B_C(1,1)\leq B_C(f,1)^2,
 \qquad f\in W^{1,2}(\Sph^2),
\end{equation}
with equality if and only if $f$ is a constant plus a first
spherical harmonic function.
\end{proposition}

\begin{proof}
By Lemma~\ref{lem:support}, we can choose a smooth positive function $h>0$ with $C=\Aop[h]$ .
For smooth functions $h_1,h_2,h_3$, the trilinear expression
\[
 \int_{\Sph^2}
 h_1\,\cof\Aop[h_2]:\Aop[h_3]\dd A_{\bar g}
\]
is symmetric in all 3 arguments. Indeed, symmetry in $h_2,h_3$
follows from the 2-dimensional identity
$\cof D:E=\cof E:D$, while integration by parts gives symmetry
in $h_1,h_3$ because
$\Div(\cof\Aop[h_2])=0$.

Suppose that $f\in\ker B_C$. Since $\cof C$ is positive definite,
elliptic regularity applied to the weak equation for $f$ shows that
$f$ is smooth and satisfies
\begin{equation}\label{eq:kernel-C}
 \cof C:\Aop[f]=0.
\end{equation}
Set $D=\Aop[f]$. Under the identification with symmetric
endomorphisms, \eqref{eq:kernel-C} becomes
\[
 0
 =(\det C)\tr(C^{-1}D)
 =(\det C)\tr(C^{-1/2}DC^{-1/2}).
\]
The symmetric endomorphism $C^{-1/2}DC^{-1/2}$ is therefore trace
free. In dimension 2, a trace free symmetric endomorphism has
nonpositive determinant, and its determinant vanishes only when
the endomorphism itself vanishes. Consequently, $\det D\leq0$ pointwise, with equality only where $D=0$.

On the other hand, the trilinear symmetry gives
\[
 \begin{aligned}
  \int_{\Sph^2}h\det D\dd A_{\bar g}
   =\frac12\int_{\Sph^2}h\,\cof D:D\dd A_{\bar g} =\frac12\int_{\Sph^2}f\,\cof C:D\dd A_{\bar g}=0.
 \end{aligned}
\]
As $h>0$ and $\det D\leq0$, continuity implies
$\det D=0$ everywhere, and hence $D=0$. Lemma~\ref{lem:support}
now gives $f\in\Hone$. Conversely, $\Aop[\ell]=0$ for
$\ell\in\Hone$, so \eqref{eq:mixed-form} and the symmetry of $B_C$ imply
$\Hone\subset\ker B_C$. Thus $\ker B_C=\Hone$.

Consider the family of positive definite Codazzi tensors $C_t=(1-t)\bar g+tC$, for $0\leq t\leq1$, and the associated operators
\[
 L_tf
 =\cof C_t:\Aop[f]
 =\Div(\cof C_t\nabla f)+(\tr C_t)f.
\]
Realized on $L^2(\Sph^2,\dd A_{\bar g})$ with common domain $W^{2,2}(\Sph^2)$, the operators $L_t$ are self adjoint and have compact resolvent. Their coefficients depend smoothly on $t$,
and uniform positivity on the compact parameter interval $[0,1]$ gives constants $c_0,c_1>0$, independent of $t$, such that
\[
 B_{C_t}(f,f)
 \leq c_1\norm{f}_{L^2}^2
      -c_0\norm{\nabla f}_{L^2}^2,
 \qquad f\in W^{1,2}(\Sph^2).
\]
The continuous dependence of the forms $B_{C_t}$ on their common domain $W^{1,2}(\Sph^2)$ implies, by the min-max principle, that the eigenvalues of $L_t$, arranged in decreasing order and counted with multiplicity, depend continuously on $t$. The displayed
estimate also shows that these eigenvalues tend to $-\infty$
uniformly in $t$.

By previous arguments, every $L_t$ has kernel $\Hone$. Self adjointness of $L_t$ therefore makes the fixed subspace $\Hone^\perp$ invariant,
where orthogonality is taken in $L^2(\Sph^2,\dd A_{\bar g})$, and the
restriction of $L_t$ to this subspace has no eigenvalue~$0$. Its number of positive eigenvalues is consequently finite and constant on $[0,1]$.

At $t=0$, the operator $L_0=2-\Delta$ has exactly 1 positive eigenvalue, namely $2$ on the constant functions, and kernel $\Hone$.
It follows that $L_1$ has exactly 1 positive eigenvalue,
counted with multiplicity. The spectral characterization of
the index of its bilinear form on $W^{1,2}(\Sph^2)$ proves the asserted
signature of $B_C$.

Next we prove the inequality \eqref{eq:AF} and characterize the equality case. Since
\[
 B_C(1,1)=\int_{\Sph^2}\tr C\dd A_{\bar g}>0,
\]
the signature of $B_C$ just proved implies that $B_C(w,w)\leq0$ whenever
$B_C(w,1)=0$. For an arbitrary $f\in W^{1,2}(\Sph^2)$, write
\[
 f=c+f_0,
 \qquad
 c=\frac{B_C(f,1)}{B_C(1,1)},
 \qquad
 B_C(f_0,1)=0.
\]
The resulting orthogonal decomposition with respect to $B_C$
gives
\[
 B_C(f,f)
 =\frac{B_C(f,1)^2}{B_C(1,1)}+B_C(f_0,f_0)
 \leq\frac{B_C(f,1)^2}{B_C(1,1)},
\]
which proves \eqref{eq:AF} directly on $W^{1,2}(\Sph^2)$.

If equality holds, then $B_C(f_0,f_0)=0$. For every
$w\in W^{1,2}(\Sph^2)$ satisfying $B_C(w,1)=0$ and every
$t\in\R$, the nonpositivity just established yields
\[
 0\geq B_C(f_0+tw,f_0+tw)
 =2tB_C(f_0,w)+t^2B_C(w,w).
\]
Taking $t$ of either sign and letting $t\to0$ gives
$B_C(f_0,w)=0$. Since every function in $W^{1,2}(\Sph^2)$ is a
sum of a constant and such a function $w$, and
$B_C(f_0,1)=0$, we conclude that $f_0\in\ker B_C=\Hone$.
Conversely, adding a first spherical harmonic does not change
any pairing under $B_C$, so every constant plus a first
spherical harmonic function gives equality in \eqref{eq:AF}.
\end{proof}

For the tensor $C$ defined as $\cof A$ and $A$ defined in \eqref{eq:tensors-intro}, combining
\eqref{eq:pushforward} with \eqref{eq:mixed-form} therefore gives
\begin{equation}\label{eq:Q-minus-B}
 Q_{\alpha,u}(f\circ u,h\circ u)
 =-\alpha B_C(f,h),
 \qquad f,h\in W^{1,2}(\Sph^2).
\end{equation}
Since composition with $u$ is an isomorphism from $W^{1,2}(\Sph^2)$
to $W^{1,2}(\Sigma)$, the signature and kernel of $B_C$ give the
converse to Proposition~\ref{prop:index-one-forward}, that is, we have the following.

\begin{proposition}\label{prop:automatic-index}
Every smooth $\alpha$-harmonic diffeomorphism
$u:(\Sigma,g)\to(\Sph^2,\bar g)$, where $\bar g$ is the unit round
metric, satisfies
\[
 I_\alpha(u)=1,
 \qquad
 N_\alpha(u)=3,
 \qquad
 \ker J_{\alpha,u}=\Span\{u^1,u^2,u^3\}.
\]
No curvature assumption is required on $g$, and the
conclusion holds for either orientation.
\end{proposition}

\begin{proof}
If $u$ preserves orientation, Proposition~\ref{prop:signature}
and \eqref{eq:Q-minus-B} show that $Q_{\alpha,u}$ has index 1
and that
\[
 \ker J_{\alpha,u}
 =\{\ell\circ u:\ell\in\Hone\}
 =\Span\{u^1,u^2,u^3\}.
\]
Surjectivity makes the 3 coordinate functions linearly independent, so $N_\alpha(u)=3$.

If $u$ reverses orientation, let $R$ be a reflection of the
target sphere. Then $R\circ u$ is an orientation preserving
$\alpha$-harmonic diffeomorphism, and
$|d(R\circ u)|_g^2=|du|_g^2$ implies
$J_{\alpha,R\circ u}=J_{\alpha,u}$. Moreover, the coordinate
functions of $R\circ u$ span the same space as those of $u$.
The orientation preserving case therefore gives the conclusion.
\end{proof}

\subsection{Reconstruction from the Codazzi tensor}\label{subsec:codazzi-reconstruction}

In this subsection, we recover the transported domain metric from a positive Codazzi tensor and reconstruct an $\alpha$-harmonic diffeomorphism when this metric is isometric to the prescribed domain, under the determinant condition in Proposition~\ref{prop:general-reconstruction}. For the unit round domain, the metric condition becomes a curvature condition, giving the correspondence in Theorem~\ref{thm:reconstruction}.

Fix $\alpha>1$ and $\kappa>0$. For the orientation preserving
$\alpha$-harmonic diffeomorphism $u:(\Sigma,g) \to (\S^2, \bar g)$ considered above, let $v=u^{-1}$ and write
\begin{equation}\label{eq:transported-metric}
 G:=v^*g
\end{equation}
for the transported domain metric on $\Sph^2$. Traces, determinants, and tensor eigenvalues are taken with respect to the target round metric $\bar g$. As before, functions defined on $\Sigma$ in target formulas are evaluated at $v(y)$.

At $y=u(x)$, choose singular frames for $du_x$, with positive singular values $\lambda_1,\lambda_2$. Then
\[
 |du|_g^2=\lambda_1^2+\lambda_2^2,
 \qquad
 J=\lambda_1\lambda_2,
\]
and the corresponding target orthonormal frame gives
\begin{align*}
 A&=\bigl(\kappa+|du|_g^2\bigr)^{\alpha-1}
   \begin{pmatrix}
    \lambda_1/\lambda_2 & 0\\
    0 & \lambda_2/\lambda_1
   \end{pmatrix}, \, C=\bigl(\kappa+|du|_g^2\bigr)^{\alpha-1}
   \begin{pmatrix}
    \lambda_2/\lambda_1 & 0\\
    0 & \lambda_1/\lambda_2
   \end{pmatrix}, \, G=\begin{pmatrix}
    \lambda_1^{-2} & 0\\
    0 & \lambda_2^{-2}
   \end{pmatrix}.
\end{align*}
Taking determinants and traces, and comparing $C$ with $G$, yields
\begin{equation}\label{eq:C-invariants}
 \begin{aligned}
  \det C&=\bigl(\kappa+|du|_g^2\bigr)^{2(\alpha-1)},\quad \tr C=\frac{\bigl(\kappa+|du|_g^2\bigr)^{\alpha-1}|du|_g^2}{J},\quad
  C=\bigl(\kappa+|du|_g^2\bigr)^{\alpha-1}JG.
 \end{aligned}
\end{equation}
The determinant inequality $\det C>\kappa^{2(\alpha-1)}$ ensures that
the tensor $G_C$ introduced in \eqref{eq:GC-intro} is a smooth Riemannian
metric. When $C$ arises from $u$, the quantity $(\det C)^{1/(2(\alpha-1))}-\kappa$ is precisely the transported
energy density $|du|_g^2$, and
\eqref{eq:C-invariants} gives $G_C=G$.

\begin{proposition}\label{prop:general-reconstruction}
Let $(\Sigma,g)$ be a smooth closed connected oriented
Riemannian surface. The assignment of the transported Codazzi tensor induces a bijection between smooth orientation preserving
$\alpha$-harmonic diffeomorphisms
$u:(\Sigma,g)\to(\Sph^2,\bar g)$, modulo precomposition by orientation preserving domain isometries, and smooth positive definite Codazzi tensors $C$ on $(\Sph^2,\bar g)$ satisfying $\det C>\kappa^{2(\alpha-1)}$ and such that $(\Sph^2,G_C)$ is orientation preservingly isometric to $(\Sigma,g)$.
\end{proposition}

\begin{proof}
If $u$ is an orientation preserving $\alpha$-harmonic diffeomorphism, Proposition~\ref{prop:target-conductivity} and
\eqref{eq:codazzi-equivalence} show that its transported tensor
$C$ is positive definite and Codazzi. Since $|du|_g^2>0$,
\eqref{eq:C-invariants} gives the strict determinant inequality $\det C>\kappa^{2(\alpha-1)}$,
while $G_C=(u^{-1})^*g$ proves the metric condition.

Conversely, let the Codazzi tensor $C$ satisfy the stated conditions, choose an
orientation preserving isometry
\[
 v:(\Sph^2,G_C)\longrightarrow(\Sigma,g),
 \qquad
 u=v^{-1},
\]
and fix $y\in\Sph^2$, with $x=v(y)$. Diagonalize $C_y$ in a $\bar g$-orthonormal frame, with eigenvalues $c_1,c_2>0$. The same
frame also diagonalizes $G_C$. Because $v$ is an isometry, the squared
singular values of $du_x$ are the reciprocals of the eigenvalues
of $G_C$ relative to $\bar g$. Formula~\eqref{eq:GC-intro}
therefore gives
\[
 \begin{aligned}
  \lambda_1^2
  =\frac{c_2}{c_1+c_2}
    \bigl((c_1c_2)^{1/(2(\alpha-1))}-\kappa\bigr),\quad
  \lambda_2^2
  =\frac{c_1}{c_1+c_2}
    \bigl((c_1c_2)^{1/(2(\alpha-1))}-\kappa\bigr).
 \end{aligned}
\]
Taking their sum and the positive square root of their product,
we obtain
\[
 \begin{aligned}
  |du_x|_g^2
  =(\det C)^{1/(2(\alpha-1))}-\kappa,\quad
  J(x)
  =\frac{(\det C)^{1/2}
    \bigl((\det C)^{1/(2(\alpha-1))}-\kappa\bigr)}{\tr C},
 \end{aligned}
\]
where $C$ is evaluated at $y$. In particular, we have
\[
 \bigl(\kappa+|du_x|_g^2\bigr)^{\alpha-1}
 =(\det C)^{1/2},
\]
so the same frame gives
\[
 \frac{\bigl(\kappa+|du_x|_g^2\bigr)^{\alpha-1}}{J(x)}
 \,du_xdu_x^*
 =\begin{pmatrix}
   c_2 & 0\\
   0 & c_1
  \end{pmatrix}
 =\cof C_y.
\]
Thus the tensor defined from $u$ by \eqref{eq:tensors-intro} is
$\cof C$, whose cofactor is $C$, while taking traces recovers
the second identity in \eqref{eq:C-invariants}.

To verify the Euler-Lagrange equation for $u$, using the Codazzi condition of $C$
and \eqref{eq:codazzi-equivalence} we obtain
$\Div_{\bar g}(\cof C)=0$. For the round coordinate functions $y^i$,
$1\leq i\leq3$, the identity $\Hess_{\bar g} y^i=-y^i\bar g$ then gives
\[
 -\Div_{\bar g}\bigl(\cof C\,\nabla_{\bar g} y^i\bigr)
 =(\tr C)y^i.
\]
Given $f\in C^\infty(\Sigma)$, test this equation against
$f\circ v$ and change variables by $y=u(x)$. The recovered
conductivity and trace identities yield
\[
 \begin{aligned}
  &\int_\Sigma
   \bigl(\kappa+|du|_g^2\bigr)^{\alpha-1}
   \ip{\nabla_g u^i}{\nabla_g f}_g\dd A_g   =\int_\Sigma
   \bigl(\kappa+|du|_g^2\bigr)^{\alpha-1}
   |du|_g^2u^if\dd A_g,
 \end{aligned}
\]
that is, $u$ satisfies \eqref{eq:sphere} and is $\alpha$-harmonic.

Precomposition by an orientation preserving isometry of the domain preserves the transported tensor, as follows directly from
\eqref{eq:tensors-intro}. To prove injectivity of the assignment on the quotient space,
suppose that $u_1$ and $u_2$ give the same tensor $C$, and let
$v_i=u_i^{-1}$. Both inverse maps are orientation preserving
isometries from $(\Sph^2,G_C)$ to $(\Sigma,g)$, so $\phi=v_2\circ v_1^{-1}$ is an orientation preserving isometry of $\Sigma$ and
$u_2=u_1\circ\phi^{-1}$. Together with the construction above from Codazzi tensor $C$ to the $\alpha$-harmonic diffeomorphism $u$, this proves the bijection of the assignment, hence it completes the proof of Proposition \ref{prop:general-reconstruction}.
\end{proof}

When the domain metric $g$ is the unit round metric, the isometry
condition in Proposition~\ref{prop:general-reconstruction} can
be expressed entirely in terms of the Gaussian curvature of
$G_C$.

\begin{theorem}\label{thm:reconstruction}
Smooth orientation preserving $\alpha$-harmonic diffeomorphisms
of the unit round 2-sphere, modulo precomposition by 
rotations of $\S^2$, correspond bijectively to smooth positive definite
Codazzi tensors $C$ on $(\Sph^2,\bar g)$ satisfying
\begin{equation}\label{eq:curvature-equation}
 \det C>\kappa^{2(\alpha-1)},
 \qquad
 K_{G_C}=1.
\end{equation}
Rotations of $\S^2$ are $\alpha$-harmonic and correspond to $C=(\kappa+2)^{\alpha-1}\bar g$, for which
$G_C=\bar g$. Moreover, the differential of the convex embedding $X$ defined in Lemma~\ref{lem:support} is represented by $C$ after raising an index with $\bar g$; when $X$ is parametrized by the domain $\S^2$ as $X\circ u$, its outward Gauss
map is $u$.
\end{theorem}

\begin{proof}
For a smooth orientation preserving $\alpha$-harmonic diffeomorphism $u$ of the round sphere,
\eqref{eq:C-invariants} and \eqref{eq:GC-intro}
give $G_C=(u^{-1})^*g$, and hence $K_{G_C}=1$.

Conversely, suppose that $C$ is a smooth positive definite
Codazzi tensor on $\S^2$ satisfying
\eqref{eq:curvature-equation}. The metric $G_C$ is smooth and
complete by the compactness of $\Sph^2$. Since $\Sph^2$ is simply
connected and $K_{G_C}=1$, the classification of complete
simply connected surfaces of constant curvature gives an
isometry to the unit round sphere. Postcomposing with a
reflection of the round domain $\S^2$ if necessary, we obtain an
orientation preserving isometry $v:(\Sph^2,G_C)\longrightarrow(\Sph^2,g)$.
Proposition~\ref{prop:general-reconstruction} shows that
$u=v^{-1}$ is $\alpha$-harmonic and has transported Codazzi tensor $C$.
Any 2 such isometries $v$ differ by a rotation of domain, so Proposition~\ref{prop:general-reconstruction} gives the stated bijection between positive Codazzi tensors satisfying \eqref{eq:curvature-equation} and $\alpha$-harmonic diffeomorphisms of round 2-sphere modulo
precomposition by rotations. Substitution of rotations in
\eqref{eq:GC-intro} verifies that
$C=(\kappa+2)^{\alpha-1}\bar g$ and gives $G_C=\bar g$.

Finally, Lemma~\ref{lem:support} provides the convex embedding
$X$ whose differential is represented by $C$ through $\bar g$. Its outward Gauss map $N$ satisfies
$N\circ X=\mathrm{Id}_{\Sph^2}$, whence
$ N\circ(X\circ u)=u$, 
which proves the last assertion in Theorem \ref{thm:reconstruction}.
\end{proof}

\begin{proof}[Proof of Theorems~\ref{thm:main-spectral}
 and~\ref{thm:main-convex-energy}]
The spectral assertions for parts \ref{thm:main-spectral 1}, \ref{thm:main-spectral 2} and \ref{thm:main-spectral 3} on arbitrary closed surfaces follow from
Proposition~\ref{prop:weighted-spectrum}, Lemma~\ref{lem:nodal-count},
\eqref{eq:general-multiplicity}, and
Proposition~\ref{prop:normal-splitting}. For a genus 0
domain with any smooth metric,
Propositions~\ref{prop:index-one-forward} and
\ref{prop:automatic-index} give the 2 implications in the
scalar index 1 classification, together with the stated
nullity and kernel.

The convex realization in Theorem \ref{thm:main-convex-energy} follows from
Proposition~\ref{prop:target-conductivity},
\eqref{eq:codazzi-equivalence}, and Lemma~\ref{lem:support}.
Proposition~\ref{prop:signature} gives the signature and equality
case of the Alexandrov-Fenchel inequality for $B_C$, while \eqref{eq:Q-minus-B}
identifies the transported scalar form as $-\alpha B_C$.
Finally, Proposition~\ref{prop:general-reconstruction} proves
the reconstruction for a prescribed domain metric, and
Theorem~\ref{thm:reconstruction} gives its round domain 2-sphere specialization.
\end{proof}

\vskip1cm

\section{Energy estimates and preliminary rigidity}\label{sec:energy}

In this section, we establish the energy bounds and preliminary rigidity results needed for the compactness argument. We first derive the domain balance identity in Subsection~\ref{subsec:domain-balance-identity} and combine it with the Alexandrov-Fenchel inequality to obtain the energy estimate in Subsection~\ref{subsec:mixed-form-energy}. We then prove local rigidity near orthogonal maps in Subsection~\ref{sec:fixed} and rigidity for twisted $1$-equivariant diffeomorphisms in Subsection~\ref{subsec:twisted-equivariant}.

Throughout this section, $\alpha>1$ and
$\kappa>0$ are fixed and $g$ is the unit round domain metric. We use the
round target metric $\bar g$ and the domain embedding $x$ fixed in
Subsection~\ref{subsec:background}.

\subsection{The domain balance identity}\label{subsec:domain-balance-identity}

In this subsection, we derive the domain balance identity for $\alpha$-harmonic maps from the unit round $2$-sphere and express it in the equivalent forms \eqref{eq:source-balance} and~\eqref{eq:Omega-balance}. The nonnegative weight in the latter identity will allow us to pass this constraint to the limiting energy measures.

For $\kappa=1$, the domain balance identity appears in the conformal variation argument of
Sacks-Uhlenbeck~\cite[proof of Lemma~5.3, p.~20]{SacksUhlenbeck1981}. Differently, we derive the balance identity for arbitrary $\kappa>0$ using the stress tensor $S$ of $E_{\alpha,\kappa}$,
\begin{equation}\label{eq:stress}
 S:=\frac12(\kappa+|du|_g^2)^\alpha g
 -\alpha(\kappa+|du|_g^2)^{\alpha-1}u^*\bar g.
\end{equation}

To verify that $\Div_gS=0$, choose normal orthonormal coordinates
$(x^1,x^2)$ at $p\in\Sph^2$ and write $u_i=\partial_i u$. All component calculations in this paragraph are evaluated
at $p$, where $g_{ij}=\delta_{ij}$ and the domain Christoffel symbols vanish. Since mixed partial derivatives commute and
the first derivatives of $g$ vanish at $p$, we have
\[
 \sum_{i=1}^2\langle u_i,\partial_i\partial_j u\rangle
 =\sum_{i=1}^2\langle u_i,\partial_j\partial_i u\rangle
 =\frac12\partial_j|du|_g^2
 \qquad\text{at }p.
\]
Differentiating \eqref{eq:stress}, we therefore obtain
\begin{align*}
 (\Div_gS)_j
 &=
 \frac{\alpha}{2}(\kappa+|du|_g^2)^{\alpha-1}
 \partial_j|du|_g^2
 -\alpha(\kappa+|du|_g^2)^{\alpha-1}
 \sum_{i=1}^2\langle u_i,\partial_i\partial_j u\rangle - \alpha\left\langle
 \Div_g\bigl((\kappa+|du|_g^2)^{\alpha-1}\nabla u\bigr),
 u_j
 \right\rangle
 \\*
 &=
 -\alpha\left\langle
 \Div_g\bigl((\kappa+|du|_g^2)^{\alpha-1}\nabla u\bigr),
 u_j
 \right\rangle
 =0.
\end{align*}
The last equality follows from \eqref{eq:sphere} and
$\langle u,u_j\rangle=0$. Since the domain has dimension 2, the trace
of $S$ is
\begin{equation}\label{eq:stress-trace}
\begin{aligned}
 \tr_gS =(\kappa+|du|_g^2)^\alpha
   -\alpha(\kappa+|du|_g^2)^{\alpha-1}|du|_g^2 =(\kappa+|du|_g^2)^{\alpha-1}
   \bigl(\kappa-(\alpha-1)|du|_g^2\bigr).
\end{aligned}
\end{equation}

For every $\xi\in\R^3$, the domain coordinate function $\xi\cdot x$
satisfies
\[
 \Hess_g(\xi\cdot x)=-(\xi\cdot x)g.
\]
Testing $\Div_gS=0$ against $\nabla(\xi\cdot x)$ and integrating by parts
on the closed domain $2$-sphere gives
\[
 0
 =-\int_{\Sph^2}
   \langle S,\Hess_g(\xi\cdot x)\rangle_g\,\dd A_g
 =\int_{\Sph^2}(\xi\cdot x)\tr_gS\,\dd A_g.
\]
Because $\xi$ is arbitrary, \eqref{eq:stress-trace} yields the domain
balance identity
\begin{equation}\label{eq:source-balance}
 \int_{\Sph^2}
 (\kappa+|du|_g^2)^{\alpha-1}
 \bigl(\kappa-(\alpha-1)|du|_g^2\bigr)x\,\dd A_g
 =0.
\end{equation}
Equivalently, the domain balance identity takes the form
\begin{equation}\label{eq:MN}
 (\alpha-1)\int_{\Sph^2}
 (\kappa+|du|_g^2)^{\alpha-1}|du|_g^2x\,\dd A_g
 =
 \kappa\int_{\Sph^2}
 (\kappa+|du|_g^2)^{\alpha-1}x\,\dd A_g.
\end{equation}

For later use, we express the domain balance identity in terms of a nonnegative weight. Define
\begin{equation}\label{eq:Omega}
 \Omega_{\alpha,\kappa}(s)
 :=
 \kappa^\alpha
 +\alpha(\kappa+s)^{\alpha-1}s
 -(\kappa+s)^\alpha,
 \qquad s\geq0.
\end{equation}
By \eqref{eq:stress-trace}, we see that 
\[
 \Omega_{\alpha,\kappa}(|du|_g^2)
 =\kappa^\alpha-\tr_gS.
\]
Since $\int_{\Sph^2}x\,\dd A_g=0$, the domain balance identity
\eqref{eq:source-balance} is therefore equivalent to
\begin{equation}\label{eq:Omega-balance}
 \int_{\Sph^2}
 \Omega_{\alpha,\kappa}(|du|_g^2)x\,\dd A_g
 =0.
\end{equation}
Moreover, $\Omega_{\alpha,\kappa}$ satisfies
\begin{equation}\label{eq:Omega-derivative}
 \Omega_{\alpha,\kappa}(0)=0,
 \qquad
 \Omega'_{\alpha,\kappa}(s)
 =\alpha(\alpha-1)s(\kappa+s)^{\alpha-2}.
\end{equation}
Since $\alpha>1$ and $\kappa>0$, the function
$\Omega_{\alpha,\kappa}$ is nonnegative on $[0,\infty)$ and vanishes
only at 0. Consequently, the weight in
\eqref{eq:Omega-balance} is strictly positive precisely at the points
where $du\neq0$.

\subsection{The Alexandrov-Fenchel inequality and the energy estimate}\label{subsec:mixed-form-energy}
In this subsection, we apply the Alexandrov-Fenchel inequality \eqref{eq:AF} to the transported coordinate functions and combine it with the domain balance identity to bound the $\alpha$-energy of $\alpha$-harmonic diffeomorphisms between unit round $2$-spheres. Theorem~\ref{thm:energy-bound} collects this upper bound and the sharp lower bound for maps of degree $\pm1$, whose equality case consists precisely of orthogonal maps.

Assume first that $u:\Sph^2\to\Sph^2$ is an orientation preserving
$\alpha$-harmonic diffeomorphism. Use the transported metric $G$ from
\eqref{eq:transported-metric} and put $z=x\circ u^{-1}$. Thus $z:(\Sph^2,G)\to(\Sph^2,g)$ is an isometry.
With domain functions evaluated at $u^{-1}(y)$ in target formulas,
the area measures satisfy
$\dd A_{\bar g}=J\dd A_G$. Formula \eqref{eq:C-invariants} gives
\begin{equation}\label{eq:mixed-G}
 B_C(f,h)=\int_{\Sph^2}(\kappa+|du|_g^2)^{\alpha-1}
 \bigl(|du|_g^2fh-\ip{\nabla_Gf}{\nabla_Gh}_G\bigr)\dd A_G.
\end{equation}
Since $\sum_i z_i^2=1$ and $\sum_i|\nabla_Gz_i|_G^2=2$, changing
variables back to the domain gives
\begin{equation}\label{eq:PQ}
\begin{aligned}
 B_C(1,1)&=\int_{\Sph^2}(\kappa+|du|_g^2)^{\alpha-1}|du|_g^2\dd A_g, \qquad
 B_C(z_i,1)=\int_{\Sph^2}(\kappa+|du|_g^2)^{\alpha-1}|du|_g^2x_i\dd A_g,\\
 \sum_{i=1}^3B_C(z_i,z_i)
 &=\int_{\Sph^2}(\kappa+|du|_g^2)^{\alpha-1}(|du|_g^2-2)\dd A_g.
\end{aligned}
\end{equation}
Equipped with these notations, we have the following inequality.
\begin{proposition}\label{prop:barycenter}
Every $\alpha$-harmonic diffeomorphism between unit round 2-spheres satisfies
\begin{equation}\label{eq:barycenter}
\begin{aligned}
 0\leq\int_{\Sph^2}(\kappa+|du|_g^2)^{\alpha-1}(|du|_g^2-2)\dd A_g \leq\frac{\left|\displaystyle\int_{\Sph^2}(\kappa+|du|_g^2)^{\alpha-1}|du|_g^2x\dd A_g\right|^2}
 {\displaystyle\int_{\Sph^2}(\kappa+|du|_g^2)^{\alpha-1}|du|_g^2\dd A_g}.
\end{aligned}
\end{equation}
If the numerator on the right-hand side of \eqref{eq:barycenter} vanishes, then $u$ is orthogonal.
\end{proposition}
\begin{proof}
For the orientation preserving $\alpha$-harmonic diffeomorphism $u$, applying \eqref{eq:AF} to $z_1,z_2,z_3$ and utilizing \eqref{eq:PQ}, we obtain the upper bound for \eqref{eq:barycenter}.

For the lower bound of \eqref{eq:barycenter}, $|du|_g^2\geq2J(u)$ gives $\int_{\Sph^2}(|du|_g^2-2)\dd A_g\geq0$. Adding and subtracting the constant $(\kappa+2)^{\alpha-1}$, we have
\begin{equation}\label{eq:positive-defect}
\begin{aligned}
 \int_{\Sph^2}(\kappa+|du|_g^2)^{\alpha-1}(|du|_g^2-2)\dd A_g&=\int_{\Sph^2}
 \bigl((\kappa+|du|_g^2)^{\alpha-1}-(\kappa+2)^{\alpha-1}\bigr)
 (|du|_g^2-2)\dd A_g\\
 &\quad+(\kappa+2)^{\alpha-1}\int_{\Sph^2}(|du|_g^2-2)\dd A_g\geq0.
\end{aligned}
\end{equation}
If the numerator on the right-hand side of \eqref{eq:barycenter} vanishes, equality holds in
\eqref{eq:positive-defect}. Since the first integrand on the right-hand side of \eqref{eq:positive-defect} is
continuous, nonnegative, and 0 only when $|du|_g^2=2$, we must have  $|du|_g^2\equiv2$. If $\lambda_1,\lambda_2$ are the singular values of $du$, then $\lambda_1^2+\lambda_2^2=2$ and $J=\lambda_1\lambda_2\leq1$. The degree
identity forces $J\equiv1$, so both singular values are 1 and $u$ is
an isometry. A target reflection reduces the orientation reversing
case to the orientation preserving case and preserves all the domain
integrals in the estimate, which completes the proof of Proposition \ref{prop:barycenter}.
\end{proof}

\begin{proposition}\label{prop:weighted-mean}
Every $\alpha$-harmonic diffeomorphism between unit round 2-spheres satisfies
\begin{equation}\label{eq:weighted-mean}
 2\leq
 \frac{\displaystyle\int_{\Sph^2}(\kappa+|du|_g^2)^{\alpha-1}|du|_g^2\dd A_g}
 {\displaystyle\int_{\Sph^2}(\kappa+|du|_g^2)^{\alpha-1}\dd A_g}
 \leq\kappa+2.
\end{equation}
\end{proposition}
\begin{proof}
The lower bound in \eqref{eq:weighted-mean} follows from
Proposition~\ref{prop:barycenter}. For the upper bound, the elementary
inequalities $1-e^{-t}\leq t$ and $\log(1+t)\leq\sqrt t$, for
$t\geq0$, give
\begin{align*}
 0 \leq (\kappa+|du|_g^2)^{\alpha-1}-\kappa^{\alpha-1} \leq (\alpha-1)(\kappa+|du|_g^2)^{\alpha-1}
       \log\left(1+\frac{|du|_g^2}{\kappa}\right)\leq \frac{\alpha-1}{\sqrt\kappa}
       (\kappa+|du|_g^2)^{\alpha-1}|du|_g.
\end{align*}
Since $\int_{\Sph^2}x\dd A_g=0$, the balance identity
\eqref{eq:MN} therefore yields
\begin{equation}\label{eq:M-bound}
\begin{aligned}
 &\left|\int_{\Sph^2}
   (\kappa+|du|_g^2)^{\alpha-1}|du|_g^2x\dd A_g\right|\\
 &\quad=\frac{\kappa}{\alpha-1}
   \left|\int_{\Sph^2}
   \left((\kappa+|du|_g^2)^{\alpha-1}-\kappa^{\alpha-1}\right)
   x\dd A_g\right|\\
 &\quad\leq\sqrt\kappa\int_{\Sph^2}
   (\kappa+|du|_g^2)^{\alpha-1}|du|_g\dd A_g\\
 &\quad\leq\sqrt\kappa
   \left(\int_{\Sph^2}
   (\kappa+|du|_g^2)^{\alpha-1}\dd A_g\right)^{1/2}
   \left(\int_{\Sph^2}
   (\kappa+|du|_g^2)^{\alpha-1}|du|_g^2\dd A_g\right)^{1/2},
\end{aligned}
\end{equation}
where in the last inequality we use the Cauchy-Schwarz inequality.

Squaring \eqref{eq:M-bound} and substituting into
\eqref{eq:barycenter}, we obtain
\[
 \int_{\Sph^2}
 (\kappa+|du|_g^2)^{\alpha-1}(|du|_g^2-2)\dd A_g
 \leq
 \kappa\int_{\Sph^2}
 (\kappa+|du|_g^2)^{\alpha-1}\dd A_g,
\]
which is exactly the upper bound of \eqref{eq:weighted-mean}.
\end{proof}

\begin{theorem}\label{thm:energy-bound}
Let $\alpha>1$ and $\kappa>0$. Every smooth $\alpha$-harmonic diffeomorphism $u:\Sph^2\to\Sph^2$ satisfies
\begin{equation}\label{eq:energy-upper}
 E_{\alpha,\kappa}(u)\leq2\pi[2(\kappa+1)]^\alpha.
\end{equation}
Every smooth degree-$\pm1$ map satisfies, without a criticality assumption,
\begin{equation}\label{eq:energy-lower}
 E_{\alpha,\kappa}(u)\geq2\pi(\kappa+2)^\alpha,
\end{equation}
with equality if and only if $u\in O(3)$.
\end{theorem}
\begin{proof}
By Proposition~\ref{prop:weighted-mean} and H\"older's inequality, for the $\alpha$-harmonic diffeomorphism $u$ we have 
\begin{align*}
 E_{\alpha,\kappa}(u)
 =\frac12\int_{\Sph^2}(\kappa+|du|_g^2)^{\alpha-1}
 (\kappa+|du|_g^2)\dd A_g &\leq(\kappa+1)\int_{\Sph^2}(\kappa+|du|_g^2)^{\alpha-1}\dd A_g\\
 &\leq(\kappa+1)(4\pi)^{1/\alpha}
 \bigl(2E_{\alpha,\kappa}(u)\bigr)^{(\alpha-1)/\alpha},
\end{align*}
which gives \eqref{eq:energy-upper}.

For a degree-$\pm1$ map $u$, it satisfies $\int_{\Sph^2}|du|_g^2\dd A_g\geq8\pi$. Jensen's inequality gives
\[
 E_{\alpha,\kappa}(u)
 \geq2\pi\left(\kappa+\frac1{4\pi}\int_{\Sph^2}|du|_g^2\dd A_g\right)^\alpha
 \geq2\pi(\kappa+2)^\alpha.
\]
Equality requires $|du|_g^2 \equiv 2|J(u)| \equiv2$, that is, $u$ is a local isometry. A local isometry between unit round 2-spheres is orthogonal. Conversely, an orthogonal map attains the bound. This completes the proof of Theorem \ref{thm:energy-bound}.
\end{proof}

\subsection{Local rigidity at an orthogonal map}\label{sec:fixed}
In this subsection, we compute the Hessian at an orthogonal map and show that its kernel consists exactly of the rotational Killing fields, which yields the local rigidity in Proposition~\ref{prop:local-rigidity}. We then record the Morse index and nullity of isometric equatorial embeddings in Corollary~\ref{cor:equatorial-index}.

Target isometries reduce the Hessian calculation to $u_0(x)=x$.
For a tangent vector field $V$ on the round     $2$-sphere along $u_0$, we have 
\[
 |dV|^2=|\nabla V|^2+|V|^2,\qquad
 \ip{du_0}{dV}=\Div V.
\]
Thus, \eqref{eq:hessian} becomes
\begin{equation}\label{eq:rotation-hessian}
D^2E_{\alpha,\kappa}(u_0)[V,V]
 =\alpha(\kappa+2)^{\alpha-1}\int\bigl(|\nabla V|^2-|V|^2\bigr)\dd A+2\alpha(\alpha-1)(\kappa+2)^{\alpha-2}\int(\Div V)^2\dd A.
\end{equation}
Let $\Killing$ denote the 3-dimensional space of Killing fields on $\Sph^2$.

\begin{lemma}\label{lem:harmonic-Jacobi}
The bilinear form $\int(|\nabla V|^2-|V|^2) dA$ is nonnegative. Its kernel is
\begin{equation}\label{eq:harmonic-kernel}
 \Killing\oplus\{\nabla(\xi\cdot x):\xi\in\R^3\}.
\end{equation}
\end{lemma}
\begin{proof}
We identify vector fields with their metric dual $1$-forms using $g$.
Let $\nabla^*$ denote the formal adjoint of the round connection with
respect to the $L^2$ inner product induced by $g$ and $\dd A_g$, and let
$\Delta_H$ be the nonnegative Hodge Laplacian. On the unit round sphere, the Weitzenb\"ock formula gives
\[
 \Delta_H=\nabla^*\nabla+\Id,
 \qquad
 \nabla^*\nabla-\Id=\Delta_H-2\Id.
\]
Since there are no harmonic $1$-forms on $\Sph^2$, by Hodge decomposition, every smooth $1$-form is therefore the sum of an exact form and a coexact form. Since $\Delta_H$ commutes with $d$ and the Hodge star, the scalar spherical harmonics give a complete family of eigenforms $dY_\ell$ and $*dY_\ell$, for $\ell\geq1,$ with eigenvalues $\ell(\ell+1)$, where $Y_\ell$ ranges over the
spherical harmonics of degree $\ell$. Thus the smallest eigenvalue of $\Delta_H$ is $2$, and its eigenspace is
\[
 \bigl\{d(\xi\cdot x)+*d(\eta\cdot x):
       \xi,\eta\in\R^3\bigr\}.
\]
It follows that $\nabla^*\nabla-\Id$ is nonnegative. Under
identification induced by the round metric $g$, the 1-forms $d(\xi\cdot x)$ correspond to
$\nabla(\xi\cdot x)$, whereas the metric duals of $*d(\eta\cdot x)$ are precisely the rotational Killing fields. Hence the quadratic form $\int(|\nabla V|^2-|V|^2) dA$ is nonnegative and has kernel
\[
 \Killing\oplus
 \bigl\{\nabla(\xi\cdot x):\xi\in\R^3\bigr\},
\]
as claimed in Lemma \ref{lem:harmonic-Jacobi}.
\end{proof}

\begin{proposition}\label{prop:local-rigidity}
Fix $\alpha>1$, $\kappa>0$ and $0<\gamma<1$. The identity map $u_0$ has Morse index 0 and  nullity 3 for $D^2 E_{\alpha,\kappa}$. There is $\eta=\eta(\alpha,\kappa,\gamma)>0$ such that every $\alpha$-harmonic map $u$ between round $2$-spheres is orthogonal whenever
\[
 \inf_{R\in O(3)}\norm{u-R}_{C^{2,\gamma}(\Sph^2)}<\eta.
\]
\end{proposition}
\begin{proof}
By Lemma~\ref{lem:harmonic-Jacobi}, the first term in
\eqref{eq:rotation-hessian} is nonnegative, and its kernel is
\[
 \Killing\oplus
 \{\nabla(\xi\cdot x):\xi\in\R^3\}.
\]
The second term in \eqref{eq:rotation-hessian} is also nonnegative. Since every Killing field is divergence free and
\[
 \Div_g\nabla(\xi\cdot x)=-2\xi\cdot x,
\]
the 2 terms in \eqref{eq:rotation-hessian} vanish simultaneously precisely on $\Killing$.
Thus $u_0$ has Morse index 0 and nullity 3, and the kernel of $D^2 E_{\alpha,\kappa}$ is the tangent space to the target rotation orbit.

We first prove local rigidity near $u_0$. For map $u$ sufficiently close to $u_0$, it can be represented by the exponential chart
\[
 u(x)=\exp_xV(x),
\]
for some smooth vector field $V$.
We impose the slice condition $V\perp\Killing$ with respect to the $L^2$ inner product on round 2-sphere. To justify this normalization, choose a
basis $K_1,K_2,K_3$ of $\Killing$ and write $K_i(x)=A_ix$,
where $A_i\in\mathfrak{so}(3)$. For $t\in\R^3$ near 0, set
\[
 R(t)=\exp\left(\sum_{j=1}^3t_jA_j\right),\qquad
 \Phi_i(t,u)=\int_{\Sph^2}
 \left\langle\exp_x^{-1}\bigl(R(t)u(x)\bigr),K_i(x)\right\rangle_g
 \dd A_g.
\]
Then $\Phi(0,u_0)=0$ and
\[
 \frac{\partial\Phi_i}{\partial t_j}(0,u_0)
 =\int_{\Sph^2}\langle K_j,K_i\rangle_g\dd A_g.
\]
This Gram matrix is positive definite. The implicit function theorem
therefore gives, for every $u$ sufficiently close to $u_0$, a unique
small parameter $t=t(u)$ such that $\Phi(t(u),u)=0$.
Thus, after replacing $u$ by  $Ru$ for a small target rotation $R$, we have $V \perp \mathcal{K}$.

Consider the Banach spaces
\[
 X=C^{2,\gamma}(T\Sph^2)\cap\Killing^\perp,\qquad
 Y=C^{0,\gamma}(T\Sph^2)\cap\Killing^\perp,
\]
where orthogonality is understood in $L^2$ inner product. Let
$\mathcal F:X\to Y$ be the projected $L^2$ gradient of the $\alpha$-energy in exponential coordinates, characterized by
\[
 \int_{\Sph^2}\langle\mathcal F(V),W\rangle_g\dd A_g
 =
 \left.\frac{d}{ds}\right|_{s=0}
 E_{\alpha,\kappa}\bigl(x\mapsto\exp_x(V(x)+sW(x))\bigr),
 \qquad W\in X,
\]
which is a smooth map near 0, with $\mathcal F(0)=0$.
Its derivative at $0$ is the Hessian operator $D^2 E_{\alpha,\kappa}$ restricted to $X$, which is self adjoint and strongly elliptic, with kernel
$\Killing$. The Fredholm alternative and Schauder estimates therefore
show that $D\mathcal F(0):X\longrightarrow Y$ is an isomorphism. A second application of the implicit function theorem implies that $V=0$ is the unique sufficiently small solution
of $\mathcal F(V)=0$.

Let $R=R(t(u))$ be the target rotation chosen above and set
$\widetilde u=R\circ u$. Then
$\widetilde u(x)=\exp_xV(x)$ for a sufficiently small vector field
$V\in X$. Since $E_{\alpha,\kappa}$ is invariant under target rotations,
$\widetilde u$ is $\alpha$-harmonic whenever $u$ is. Its first variation
therefore vanishes along every variation
$x\mapsto\exp_x(V(x)+sW(x))$, with $W\in X$, which gives $\mathcal F(V)=0$. The local uniqueness established above implies $V=0$, and hence $R\circ u=\widetilde u=u_0$, that is, $u=R^{-1}$ is a rotation.

Finally, if $u$ is sufficiently close to some $R\in O(3)$, then
$R^{-1}\circ u$ is $\alpha$-harmonic and close to $u_0$.
The preceding argument applies, and $u$ is orthogonal.
Since target orthogonal transformations preserve the ambient
$C^{2,\gamma}$ norm, the same neighborhood size works for every
$R\in O(3)$. This proves the assertion of Proposition \ref{prop:local-rigidity}.
\end{proof}

\begin{corollary}\label{cor:equatorial-index}
An isometric equatorial embedding $\Sph^2\to\Sph^n$, $n\geq2$, has scalar index 1, scalar nullity 3, Morse index $n-2$, and nullity $3n-3$ for every $\alpha>1$ and $\kappa>0$.
\end{corollary}
\begin{proof}
At an isometric equatorial embedding, $J_{\alpha,u}=\Delta-2$, whose scalar index and nullity are 1 and 3, respectively. Apply Proposition~\ref{prop:normal-splitting} and the Morse index and nullity from Proposition~\ref{prop:local-rigidity}.
\end{proof}
\begin{remark}
Proposition~\ref{prop:local-rigidity} admits an infinitesimal
interpretation in terms of the target Codazzi tensor. The identity map $u_0$ and its target rotations determine the same round tensor $C_0=(\kappa+2)^{\alpha-1}\bar g$. Let $C_s$ be a smooth family of positive Codazzi tensors with respect to $\bar g$, with expansion $C_s=C_0+sH+O(s^2),$
and let $G_s$ be the corresponding metrics given by
\eqref{eq:GC-intro}. Suppose that $C_s$ satisfies
\eqref{eq:curvature-equation}, so that $K_{G_s}=1$.
All operators and contractions below are taken with respect to
$\bar g$ and dots denote
derivatives at $s=0$.

Since $G_0=\bar g$, differentiating \eqref{eq:GC-intro}
gives
\begin{equation}\label{eq:metric-linearization}
 \dot G
 =\frac{1}{(\kappa+2)^{\alpha-1}}
 \left[
 H-\left(\frac12+\frac{\kappa+2}{4(\alpha-1)}\right)
 (\tr H)\bar g
 \right].
\end{equation}
The Gaussian curvature variation at the unit round metric is
\[
 2\dot K
 =-\tr\dot G
  +\Div\bigl((\Div\dot G)^\sharp\bigr)
  +\Delta(\tr\dot G).
\]
As $H$ is Codazzi, we have $\Div H=d(\tr H)$. Substituting
\eqref{eq:metric-linearization} into the above Gaussian curvature variation formula
therefore reduces the equation $\dot K=0$ to
\begin{equation}\label{eq:trace-linearization}
 \Delta(\tr H)
 =\frac{2(\kappa+2)}{2(\alpha-1)+\kappa+2}\tr H.
\end{equation}
For $\alpha>1$ and $\kappa>0$, the coefficient on the right lies
strictly between $0$ and $2$. Since the first positive eigenvalue of
the round Laplacian is $2$, equation
\eqref{eq:trace-linearization} implies $\tr H=0$. Since trace free
Codazzi tensor on the round sphere vanishes identically, it follows that  $H=0$.

Thus the linearized curvature equation at $C_0$ has trivial kernel
among Codazzi variations. This gives the infinitesimal tensor
interpretation of Proposition~\ref{prop:local-rigidity}, where the
rotational freedom in the map formulation disappears since all rotations determine the same tensor $C_0$.
\end{remark}

\subsection{Twisted equivariant diffeomorphisms}\label{subsec:twisted-equivariant}

In this subsection, we use the angular conservation law to eliminate the twist and then analyze the radial profile to rule out nonorthogonal solutions. This proves the rigidity of smooth twisted $1$-equivariant $\alpha$-harmonic diffeomorphisms stated in Proposition~\ref{prop:equivariant}.

After orthogonal changes of domain and target coordinates, let
$R_\varphi$ denote rotation through angle $\varphi$ about the vertical
axis on each sphere. A map $u$ is called $1$-equivariant if
\[
 u(R_\varphi x)=R_\varphi u(x)
 \qquad\text{for every }x\in\Sph^2
 \text{ and }\varphi\in\R.
\]
Choose geodesic polar coordinates $(r,\theta)$ on the domain, with $g=\dd r^2+\sin^2r\,\dd\theta^2$. Away from the poles, the equivariance relation implies that $u$ has
the form
\begin{equation}\label{eq:twisted-map}
 u(r,\theta)=
 \bigl(
 \sin f(r)\cos(\theta+\gamma(r)),
 \sin f(r)\sin(\theta+\gamma(r)),
 \cos f(r)
 \bigr).
\end{equation}
Here $f(r)$ is the target polar angle of $u(r,0)$, and
$\gamma(r)$ is its longitude, initially defined modulo $2\pi$.
We choose a smooth real lift $\gamma:(0,\pi)\to\R$, unique up to an additive integer multiple of $2\pi$. Thus $u$ maps the parallel at polar angle $r$ to the parallel at polar angle $f(r)$, shifting its longitude by
$\gamma(r)$. We use the terminology \emph{twisted $1$-equivariant} to allow this angular shift to depend on $r$. When $\gamma$ is
constant, a target rotation removes the shift.

For an orientation preserving $1$-equivariant diffeomorphism, the
poles are fixed and
\[
 f(0)=0,\qquad f(\pi)=\pi,\qquad
 f'(r)>0\quad\text{for }0<r<\pi.
\]
Smoothness and invertibility at the poles further imply
\[
 f'(0)>0,\qquad
 \lim_{r\nearrow\pi}
 \frac{\pi-f(r)}{\pi-r}>0.
\]
An orientation reversing $1$-equivariant diffeomorphism interchanges
the poles. Composing with the target reflection across the equatorial
plane preserves equivariance and reduces it to the orientation preserving case.

\begin{proposition}\label{prop:equivariant}
Every smooth twisted $1$-equivariant $\alpha$-harmonic
diffeomorphism of the form \eqref{eq:twisted-map}, up to orthogonal
changes of domain and target coordinates, is orthogonal.
\end{proposition}

\begin{proof}
By composing with a target reflection if necessary, we may assume
that $u$ is orientation preserving. The energy density of $u$ is
\[
 |du|_g^2
 =f'^2+\sin^2f\,\gamma'^2+\frac{\sin^2f}{\sin^2r}.
\]
Taking Variations of $f$ and $\gamma$ for $E_{\alpha,\kappa}$ give
\begin{equation}\label{eq:twist-equations}
\begin{aligned}
 \partial_r\bigl[(\kappa+|du|_g^2)^{\alpha-1}\sin r\,f'\bigr]
 &=(\kappa+|du|_g^2)^{\alpha-1}\sin f\cos f
 \left(\sin r\,\gamma'^2+\frac1{\sin r}\right),\\
 \partial_r\bigl[(\kappa+|du|_g^2)^{\alpha-1}
 \sin r\sin^2f\,\gamma'\bigr]&=0.
\end{aligned}
\end{equation}
The second equation gives a conserved quantity, which vanishes
by smoothness at the poles. Indeed, smoothness of $u$ gives upper bounds for $(\kappa+|du|_g^2)^{\alpha-1}$ and $|u_r|$ near $r=0$, while
$|\sin f\,\gamma'|\leq |u_r|$. Hence, we get
\[
 \left|
 (\kappa+|du|_g^2)^{\alpha-1}\sin r\sin^2f\,\gamma'
 \right|
 \leq
 (\kappa+|du|_g^2)^{\alpha-1}
 \sin r\,|\sin f|\,|u_r|
 \longrightarrow0
 \qquad\text{as }r\searrow0.
\]
Since $\sin r$ and $\sin f$ are strictly positive away from the pole of $\S^2$,
$\gamma'=0$ away from the pole by \eqref{eq:twist-equations}. After taking a rotation on target 2-sphere, we may assume that $\gamma \equiv 0$.

To determine the radial profile of $u$, introduce
\[
 t=\log\tan\frac r2,\qquad
 y(t)=\log\tan\frac{f(r(t))}{2},\qquad  w(t)=y(t)-t.
\]
The first equation in \eqref{eq:twist-equations} becomes
\begin{equation}\label{eq:cylinder-profile}
\begin{aligned}
 y_{tt}
 +\frac{(\alpha-1)\partial_t|du|_g^2}{\kappa+|du|_g^2}\,y_t
 +\tanh y\,(1-y_t^2)=0,\qquad \text{with}\quad
 |du|_g^2=\frac{\cosh^2t}{\cosh^2y}(1+y_t^2).
\end{aligned}
\end{equation}
At a critical point of $w$, we have $y_t=1$. Differentiating
the expression for $|du|_g^2$ and using the first equation in
\eqref{eq:cylinder-profile}, we get
\[
\begin{aligned}
 \frac{\partial_t|du|_g^2}{|du|_g^2}
 =2(\tanh t-\tanh y)+w_{tt},\qquad
 w_{tt} =-\frac{(\alpha-1)\partial_t|du|_g^2}{\kappa+|du|_g^2}.
\end{aligned}
\]
Eliminating $\partial_t|du|_g^2$ from these identities and using
$y=t+w$, we obtain
\begin{equation}\label{eq:profile-maximum}
 w_{tt}
 =\frac{2(\alpha-1)|du|_g^2}{\kappa+\alpha|du|_g^2}
 \bigl(\tanh(t+w)-\tanh t\bigr)
 \qquad\text{when }w_t=0.
\end{equation}
Since the coefficient in \eqref{eq:profile-maximum} is positive, and $\tanh$ is strictly increasing, $w$ has no positive local maximum and no negative local minimum.

Let
\[
 c=f'(0)>0,\qquad
 d=\lim_{r\nearrow\pi}\frac{\pi-f(r)}{\pi-r}>0.
\]
The pole expansions in Appendix~\ref{app:poles} give
\begin{align*}
 &w(-\infty)=\log c,
 &&\sgn w_t=\sgn(c-1)
 \quad\text{near }-\infty\quad(c\ne1),\\
 &w(+\infty)=-\log d,
 &&\sgn w_t=\sgn(d-1)
 \quad\text{near }+\infty\quad(d\ne1).
\end{align*}
If $c>1$, then $w$ is positive and strictly increasing for
sufficiently negative $t$. If $d\geq1$, its limit at $+\infty$
is nonpositive, if $d<1$, it is strictly decreasing for sufficiently positive $t$. In either case, $w$ attains a positive
interior maximum, contradicting to \eqref{eq:profile-maximum}.
Reversing the inequalities similarly excludes $c<1$.
Consequently, we must have  $c=1$, and $w(-\infty)=0$.

Suppose now that $w$ takes a positive value. If $w(+\infty)\leq0$,
the limits at both infinite ends force a positive interior maximum.
If $w(+\infty)>0$, then $d<1$, and $w$ is strictly decreasing
near $+\infty$ together with $w(-\infty)=0$, this again forces a positive interior maximum. Both alternatives contradict to
\eqref{eq:profile-maximum}. The corresponding argument excludes the existence of a negative value, using the absence of negative local minima.
Therefore $w\equiv0$, $y=t$ and $f(r)=r$.

The normalized map is thus the identity map. Restoring the target rotation and the orthogonal changes of coordinates proves that
the original map $u$ is orthogonal, which completes the proof of Proposition \ref{prop:equivariant}.
\end{proof}

\vskip1cm

\section{Compactness and rigidity of diffeomorphisms}\label{sec:near-one}\label{sec:compactness-rigidity}

In this section, we prove rigidity for $\alpha$-harmonic diffeomorphisms of the unit round $2$-sphere for every $\alpha>1$. We first establish smooth compactness and rigidity near 1 in Subsections~\ref{subsec:diffeo-compactness} and~\ref{subsec:diffeo-near-one}, then bound the sum of the Morse index and nullity and identify the rotational Jacobi kernel in Subsection~\ref{sec:full-index-rigidity}. Compactness on closed exponent intervals in Subsection~\ref{subsec:uniform-exponent-compactness} and local dependence on the exponent in Subsection~\ref{subsec:exponent-global-rigidity} allow us to extend the rigidity conclusion to all exponents.

Throughout this section, $\kappa>0$ is fixed, and domain norms and measures
are taken with respect to $g$ unless otherwise specified. 

\subsection{Compactness and exclusion of concentration}\label{subsec:diffeo-compactness}
In this subsection, we recall the bubble tree compactness theorem for spherical targets, including the energy identity and no-neck property (Theorem~\ref{thm:compactness-input}), and combine it with the domain balance identity to exclude concentration for orientation preserving $\alpha_j$-harmonic diffeomorphisms as $\alpha_j\searrow1$. This gives smooth convergence to a harmonic sphere of degree 1, as stated in Proposition~\ref{prop:no-bubbling-diffeo}.

 Let
$(\Sigma,g)$ be a closed Riemannian surface, let $\alpha_j\searrow1$,
and let $u_{\alpha_j}:\Sigma\to\Sph^n$ be smooth $\alpha_j$-harmonic maps with
uniformly bounded $E_{\alpha_j,\kappa}(u_{\alpha_j})$. For a domain $U\subset\Sigma$, the localized Dirichlet energy associated
with \eqref{eq:energy-intro} is denoted by
\[
 E(u;U)=\frac12\int_U|du|_g^2\,\dd A_g.
\]
We index quantities associated with the critical sequence by $\alpha_j$;
the superscripts on the bubble centers and scales label the bubbles.

The weak $W^{1,2}$ limit of $u_{\alpha_j}$ is called the \emph{base map} and is denoted by
$u_\infty$. The finite set where energy concentrates is denoted by
$\mathcal S\subset\Sigma$, and the convergence is smooth on compact
subsets of $\Sigma\setminus\mathcal S$. Near each
$p\in\mathcal S$, choose a fixed local conformal coordinate. A
\emph{bubble} is a nonconstant harmonic sphere
$\omega_\ell:\Sph^2\to\Sph^n$ obtained as the limit of rescaled maps
\begin{equation}\label{eq:bubble-rescaling}
 u_{\alpha_j}(x_{\alpha_j}^\ell+r_{\alpha_j}^\ell z)\longrightarrow\omega_\ell(z),
 \qquad x_{\alpha_j}^\ell\longrightarrow p,\quad r_{\alpha_j}^\ell\longrightarrow0,
\end{equation}
locally smoothly away from finitely many further concentration points.
Here $\Sph^2$ is identified with $\R^2\cup\{\infty\}$, and the
limiting map extends across the punctures and $\infty$ by the removable
singularity theorem in \cite{SacksUhlenbeck1981}. The rescaled equation for $u_{\alpha_j}(x_{\alpha_j}^\ell+r_{\alpha_j}^\ell z)$ is always obtained from the
original $\alpha$-harmonic map equation by a change of variables. Distinct bubbles in the same coordinate chart have separated centers and scales (see \cite{Qing1995}), that is, we have 
\begin{equation}\label{eq:scale separation}
    \frac{r_{\alpha_j}^\ell}{r_{\alpha_j}^k}+\frac{r_{\alpha_j}^k}{r_{\alpha_j}^\ell}
 +\frac{|x_{\alpha_j}^\ell-x_{\alpha_j}^k|}{r_{\alpha_j}^\ell+r_{\alpha_j}^k}
 \longrightarrow\infty\qquad \text{for any }\,\, \ell\ne k .
\end{equation}
The relative positions and scales determine the attachments of the
bubble tree, including bubbles formed at nested scales.

The \emph{energy identity} means that the base map and the bubbles
account for the entire limiting Dirichlet energy:
\begin{equation}\label{eq:energy-identity}
 \lim_{j\to\infty}E(u_{\alpha_j})
 =E(u_\infty)+\sum_{\ell=1}^LE(\omega_\ell).
\end{equation}
Equivalently, no Dirichlet energy is lost on the neck regions between
adjacent components. To state the corresponding local condition,
consider a neck containing no further bubbling region, written in the
coordinates of its parent component as
\[
 A_{\alpha_j}(r,R)=B_r(x_{\alpha_j})\setminus B_{Rr_{\alpha_j}}(x_{\alpha_j}).
\]
Then the energy identity is equivalent to
\[
 \lim_{r\searrow0}\lim_{R\to\infty}\limsup_{j\to\infty}
 E(u_{\alpha_j};A_{\alpha_j}(r,R))=0,
\]
for sufficiently small $r> 0$ and large enough $R > 0$.
The \emph{no-neck property}, also called \emph{necklessness},
asserts that
\begin{equation}\label{eq:no neck}
\lim_{r\searrow0}\lim_{R\to\infty}\limsup_{j\to\infty}
 \operatorname{osc}_{A_{\alpha_j}(r,R)}u_{\alpha_j}=0,
 \qquad
 \operatorname{osc}_A u=\sup_{x,y\in A}|u(x)-u(y)|.
\end{equation}
No neck property ensures that adjacent components meet at their attachment points, so the union of the images of the base map and all bubbles is connected.

For sequences of $\alpha$-harmonic maps from a fixed closed Riemann
surface into a compact Riemannian manifold, Jost
\cite[Theorem 4.2.1]{Jost1991} studied the energy identity for a
min-max sequence, while Chen-Tian \cite{ChenTian1999} proved that
the necks connecting the limiting harmonic components converge to
geodesics of finite length for a minimizing sequence. For a general
sequence, however, the ordinary energy identity may fail, as shown by
the counterexample of Li-Wang \cite{LiWang2015}. In their earlier
work \cite{LiWang2010}, Li-Wang established a generalized energy
identity and analyzed convergence of the necks to geodesics, with
explicit formulas for their lengths. Energy identities under additional
hypotheses were also studied in \cite{Lamm2010,LiLiuWang2017,Moore2017}.
When the target is the standard round sphere $\Sph^n$, Li-Zhu
\cite[Theorem~1.1]{LiZhu2019} proved the following bubble tree
compactness theorem, including the energy identity and the no-neck
property.

\begin{theorem}\label{thm:compactness-input}
Let $(\Sigma,g)$ be a fixed closed Riemannian surface and let
$u_{\alpha_j}:\Sigma\to\Sph^n$ be a sequence of smooth
$\alpha_j$-harmonic maps, where $\alpha_j\searrow1$ and
$\sup_jE_{\alpha_j,\kappa}(u_{\alpha_j})<\infty$. After passing to a
subsequence, still denoted by $u_{\alpha_j}$, there exist a smooth
harmonic map $u_\infty:\Sigma\to\Sph^n$, a nonnegative integer $L$,
and, for each $\ell=1,\ldots,L$, a sequence of  points $x_{\alpha_j}^\ell\in\Sigma$, a sequence of positive scales $r_{\alpha_j}^\ell$, and a nonconstant harmonic sphere $\omega_\ell:\Sph^2\to\Sph^n$ with the following properties:
\begin{enumerate}
\item The maps $u_{\alpha_j}$ converge weakly to $u_\infty$ in $W^{1,2}(\Sigma)$ and smoothly on compact subsets of
$\Sigma\setminus\mathcal S$, where
$\mathcal S=\{p_1,\ldots,p_L\}$ is the energy concentration set and $p_\ell=\lim_{j\to\infty}x_{\alpha_j}^\ell$.
\item In local conformal coordinates of $p_\ell$, the convergence in \eqref{eq:bubble-rescaling} holds for the rescaled map $u_{\alpha_j}(x_{\alpha_j}^\ell + r_{\alpha_j}^\ell z)$. Distinct bubbles represented in the same coordinate chart satisfy the separation
condition \eqref{eq:scale separation}.
\item The energy identity \eqref{eq:energy-identity} holds.
\item The no neck property \eqref{eq:no neck} holds on every neck joining adjacent harmonic components. In particular, the union of the images of $u_\infty$ and $\omega_1,\ldots,\omega_L$ is connected.
\end{enumerate}
Every point of $\mathcal S$ gives rise to at least 1 nonconstant bubble.
If $L=0$, then $\mathcal S=\varnothing$ and the convergence is smooth
on all of $\Sigma$.
\end{theorem}

For later comparisons with the standard normalization of $\alpha$-energy in \cite{SacksUhlenbeck1981}, the rescaled metric
$\widehat g=\kappa g$ satisfies
\begin{equation}\label{eq:kappa-normalization}
 |du|_{\widehat g}^2=\kappa^{-1}|du|_g^2,\qquad
 E_{\alpha,1}^{\widehat g}(u)
 =\kappa^{1-\alpha}E_{\alpha,\kappa}^g(u).
\end{equation}
Recall the domain balance identity \eqref{eq:Omega-balance}, by \eqref{eq:Omega-derivative}, we have
\begin{equation}\label{eq:Omega-integral}
 \Omega_{\alpha,\kappa}(s)=\alpha(\alpha-1)\int_0^s t(\kappa+t)^{\alpha-2}\dd t.
\end{equation}
For $1<\alpha\leq2$, the integral representation gives
\begin{equation}\label{eq:Omega-two-sided}
 \frac\alpha2(\alpha-1)(\kappa+s)^{\alpha-2}s^2
 \leq\Omega_{\alpha,\kappa}(s)
 \leq(\alpha-1)(\kappa+s)^{\alpha-2}s^2.
\end{equation}
The lower bound follows from the monotonicity of $(\kappa+t)^{\alpha-2}$. For the upper bound, compute
\[
 \frac{d}{dt}\bigl[t^2(\kappa+t)^{\alpha-2}\bigr]
 =t(\kappa+t)^{\alpha-3}(2\kappa+\alpha t)
 \geq\alpha t(\kappa+t)^{\alpha-2},
\]
and integrate from 0 to $s$.

\begin{lemma}\label{lem:noncollapse}
Let $u_{\alpha_j}:\Sph^2\to\Sph^n$ be nonconstant $\alpha_j$-harmonic maps, where $\alpha_j\searrow1$. Suppose
\begin{equation}\label{eq:noncollapse-hyp}
 \sup_jE_{\alpha_j,\kappa}(u_{\alpha_j})<\infty,
 \qquad \inf_jE(u_{\alpha_j})>0.
\end{equation}
There is no point $p\in\Sph^2$ such that $u_{\alpha_j}$ converges to a constant in $C^1_{\rm loc}(\Sph^2\setminus\{p\})$.
\end{lemma}
\begin{proof}
Discarding finitely many terms,  we can assume $1<\alpha_j\leq2$. Define
\[
 d\mu_{\alpha_j}=\Omega_{\alpha_j,\kappa}({|du_{\alpha_j}|_g^2})\dd A.
\]
The measures $\mu_{\alpha_j}$ are positive, and \eqref{eq:Omega-balance} gives $\int x\dd\mu_{\alpha_j}=0$.

We first show that $\mu_{\alpha_j}(\Sph^2)\geq c(\alpha_j-1)$ for a fixed $c>0$. By Cauchy-Schwarz inequality, there holds
\[
 \left(\int {|du_{\alpha_j}|_g^2}\dd A\right)^2
 \leq\left(\int {|du_{\alpha_j}|_g^4}(\kappa+{|du_{\alpha_j}|_g^2})^{\alpha_j-2}\dd A\right)
       \left(\int(\kappa+{|du_{\alpha_j}|_g^2})^{2-\alpha_j}\dd A\right).
\]
The second factor in above inequality is uniformly bounded, since $0\leq2-\alpha_j\leq1$ and $s^{2-\alpha_j}\leq1+s^{\alpha_j}$ for $s\geq0$. The  assumptions in \eqref{eq:noncollapse-hyp} therefore give a uniform positive lower bound for the first factor in above inequality. The lower inequality in \eqref{eq:Omega-two-sided} proves  $\mu_{\alpha_j}(\Sph^2)\geq c(\alpha_j-1)$.

Suppose, contrary to the conclusion, that the maps converge to a constant away from $p$. For each fixed $\delta>0$, ${|du_{\alpha_j}|_g^2}\to0$ uniformly on $\Sph^2\setminus B_\delta(p)$. The upper bound in \eqref{eq:Omega-two-sided}, with $\kappa$ fixed, yields
\[
 \mu_{\alpha_j}(\Sph^2\setminus B_\delta(p))=o((\alpha_j-1)).
\]
Thus the probability measures $\mu_{\alpha_j}/\mu_{\alpha_j}(\Sph^2)$ converge weakly
to $\delta_p$. Dividing the domain balance identity \eqref{eq:Omega-balance} by $\mu_{\alpha_j}(\Sph^2)$ and passing
to the limit, the resulting identity $x(p)=0$ contradicts $|x(p)|=1$. Thus, it completes the proof of Lemma \ref{lem:noncollapse}. 
\end{proof}

\begin{proposition}\label{prop:no-bubbling-diffeo}
Let $\alpha_j\searrow1$ and $u_{\alpha_j}:\Sph^2\to\Sph^2$ be a sequence of orientation preserving $\alpha_j$-harmonic diffeomorphisms. After passing to a subsequence, $u_{\alpha_j}$ converges smoothly to a harmonic 2-sphere of degree 1.
\end{proposition}
\begin{proof}
Theorem~\ref{thm:energy-bound} supplies the uniform $\alpha_j$-energy
bound, so Theorem~\ref{thm:compactness-input} gives a harmonic base map
$u_\infty$ and bubbles $\omega_1,\ldots,\omega_L$ as $\alpha_j \searrow 1$. Since each $u_{\alpha_j}$
preserves orientation, local smooth convergence on the base map and under
orientation preserving rescaling map gives a nonnegative Jacobian for every limiting component. The Hopf differential of a harmonic map $\Sph^2\to\Sph^2$ is holomorphic and hence vanishes. Consequently, each nonconstant component is holomorphic, has positive integer degree, and satisfies $E(v)=4\pi\deg v$.

The energy identity implies that the neck regions carry no Jacobian
mass. Integrating the pullback of the
target area form over the base and bubble regions, we get
\[
 1=\deg u_{\alpha_j}=\deg u_\infty+\sum_{\ell=1}^L\deg\omega_\ell.
\]
All terms on the right are nonnegative integers, and every nonconstant
component contributes at least 1. Thus there is exactly 1
nonconstant component, of degree 1, in particular,
\eqref{eq:energy-identity} also gives $E(u_{\alpha_j})\to4\pi$.

If the base map were constant, there would be exactly 1 bubble and
hence a single concentration point $p$. The convergence would then be to a constant in $C^1_{\rm loc}(\Sph^2\setminus\{p\})$, contrary to Lemma~\ref{lem:noncollapse}, whose energy lower bound assumption follows from $E(u_{\alpha_j})\geq4\pi\deg u_{\alpha_j}=4\pi$. Therefore $u_\infty$ is the
nonconstant degree 1 component and $L=0$. The final assertion of Theorem~\ref{thm:compactness-input} gives smooth convergence on all of $\Sph^2$.
\end{proof}

\subsection{The rotational limit and rigidity for \texorpdfstring{$\alpha$}{alpha} near 1}\label{subsec:diffeo-near-one}
In this subsection, we pass the domain balance identity to the limit and show that the harmonic limits of $\alpha_j$-harmonic diffeomorphisms are orthogonal when $\alpha_j\searrow1$ (Proposition~\ref{prop:rotational-limit}). We then exclude nonorthogonal maps converging smoothly to a rotation by Lemma~\ref{lem:sequential-local}, which yields rigidity for exponents sufficiently close to 1 in Theorem~\ref{thm:near-one-rigidity}.

For $s$ in a bounded interval of $(0,\infty)$, we have
\begin{equation}\label{eq:Xi}
 \frac{\Omega_{\alpha,\kappa}(s)}{\alpha-1}
 \longrightarrow\Xi_\kappa(s):=s-\kappa\log\left(1+\frac s\kappa\right)
\end{equation}
in every fixed $C^k$ norm as $\alpha\searrow1$. In addition, $\Xi_\kappa$ satisfies
\begin{equation}\label{eq:Xi-monotone}
 \Xi_\kappa(0)=0,\qquad \Xi_\kappa'(s)=\frac{s}{\kappa+s}>0\quad(s>0).
\end{equation}

\begin{lemma}\label{lem:Mobius-balance}
Let $u:\Sph^2\to\Sph^2$ be an orientation preserving M\"obius transformation. If
\begin{equation}\label{eq:limit-balance}
 \int_{\Sph^2}\Xi_\kappa(|du|^2)x\dd A=0,
\end{equation}
then $u$ is a rotation.
\end{lemma}
\begin{proof}
Since the balance identity \eqref{eq:limit-balance} is invariant under domain and target
rotations, for $R_1,R_2\in SO(3)$, a change of variables gives
\[
 \int_{\Sph^2}\Xi_\kappa(|d(R_2\circ u\circ R_1)|^2)x\dd A
 =R_1^{-1}\int_{\Sph^2}\Xi_\kappa(|du|^2)x\dd A =0.
\]
We represent the M\"obius transformation by a determinant 1 complex matrix. The singular value decomposition of the representing matrix has the form
\[
 U\begin{pmatrix}s&0\\0&s^{-1}\end{pmatrix}V,
 \qquad U,V\in SU(2),\quad s\geq1.
\]
The unitary factors induce rotations, and the diagonal factor induces
the stereographic dilation $z\mapsto s^2z$. Up to domain and target
rotations, the Jacobian therefore has the form
\[
 J_q(x)=\left(\frac{1-|q|^2}{1+|q|^2-2q\cdot x}\right)^2,
 \qquad q\in\R^3,\quad |q|<1.
\]
The case $q=0$ is a rotation. If $q\ne0$, let $\omega=q/|q|$ and $t=\omega\cdot x$. The function $J_q(t)$ is strictly increasing. Since a M\"obius map is conformal, $|du|^2=2J_q$, and by \eqref{eq:Xi-monotone},
\[
 \omega\cdot\int\Xi_\kappa(2J_q)x\dd A
 =2\pi\int_0^1t\bigl[\Xi_\kappa(2J_q(t))-\Xi_\kappa(2J_q(-t))\bigr]\dd t>0,
\]
which contradicts \eqref{eq:limit-balance}. Therefore, $q=0$ and $u$ is a rotation.
\end{proof}

\begin{proposition}\label{prop:rotational-limit}
Every sequence of $\alpha_j$-harmonic  diffeomorphisms with $\alpha_j\searrow1$ has a subsequence converging smoothly to an element of $O(3)$.
\end{proposition}
\begin{proof}
Passing to a subsequence with fixed orientation and using a reflection if
necessary, Proposition~\ref{prop:no-bubbling-diffeo} gives smooth
convergence to a harmonic map of degree 1. Hence, it is a M\"obius
transformation. The energy densities also converge smoothly. Dividing
\eqref{eq:Omega-balance} by $\alpha_j-1$ and applying \eqref{eq:Xi}
gives \eqref{eq:limit-balance},  Lemma~\ref{lem:Mobius-balance} then
identifies the limit as a rotation.
\end{proof}

As $\alpha\searrow1$, the Hessian of $E_{\alpha,\kappa}$ degenerates along the 3 gradient fields in \eqref{eq:harmonic-kernel}. The following lemma uses the domain balance identity \eqref{eq:Omega-balance} to exclude these directions and establish local rigidity uniformly as the exponent approaches 1.

\begin{lemma}\label{lem:sequential-local}
Let $u_{\alpha_j}:(\Sph^2,g)\to(\Sph^2,\bar g)$ be a sequence of smooth $\alpha_j$-harmonic maps where $\alpha_j>1$ and $\alpha_j\searrow1$. If $u_{\alpha_j}$ converges smoothly to a rotation, then $u_{\alpha_j}$ is a rotation for every sufficiently large $j$.
\end{lemma}
\begin{proof}
Arguing by contradiction, we pass to a subsequence for which no $u_{\alpha_j}$ is a rotation and postcompose with the inverse of the limiting rotation, so that the maps converge smoothly to the identity map. The rotation slice $R(t)$ constructed in the proof of Proposition~\ref{prop:local-rigidity} depends only on the round metrics and is independent of the exponent $\alpha \searrow 1$. Since target rotations preserve $E_{\alpha_j,\kappa}$, further postcomposition by rotations tending to the identity gives, with the same notation,
\[
 u_{\alpha_j}(x)=\exp_xV_{\alpha_j}(x),\qquad
 V_{\alpha_j}\perp\Killing\text{ in }L^2,\qquad
 V_{\alpha_j}\longrightarrow0\text{ smoothly}.
\]
The assumption that none of the maps is a rotation ensures that $\delta_{\alpha_j}:=\|V_{\alpha_j}\|_{L^2}>0$, which allows us to normalize the displacement fields $V_{\alpha_j}$ once we have controlled their stronger norms.

For a fixed $0<\gamma<1$, consider the $L^2$ gradient of $V\mapsto E_{\alpha,\kappa}(x\mapsto\exp_xV(x))$ in exponential coordinates. This operator depends smoothly on $(\alpha,V)$ near $(1,0)$ as a map from $\R\times C^{2,\gamma}(T\Sph^2)$ to $C^{0,\gamma}(T\Sph^2)$, since $\kappa>0$ keeps the coefficients smooth. Its value at $V=0$ vanishes for every $\alpha$, while its derivative in $V$ at $(1,0)$ is the harmonic Jacobi operator $\mathcal J_0=\nabla^*\nabla-\Id$. Taylor expansion of the equation satisfied by $V_{\alpha_j}$ and the global Schauder estimate therefore give constants $C_1,C_2$, independent of $j$, such that
\begin{align}
 \norm{\mathcal J_0V_{\alpha_j}}_{C^{0,\gamma}}
 &\leq C_1\bigl((\alpha_j-1)+\norm{V_{\alpha_j}}_{C^{2,\gamma}}\bigr)
                 \norm{V_{\alpha_j}}_{C^{2,\gamma}},\label{eq:Jacobi-remainder}\\
 \norm{V_{\alpha_j}}_{C^{2,\gamma}}
 &\leq C_2\bigl(\norm{\mathcal J_0V_{\alpha_j}}_{C^{0,\gamma}}+\delta_{\alpha_j}\bigr).
 \label{eq:Jacobi-estimate}
\end{align}
Here the $L^2$ term in \eqref{eq:Jacobi-estimate} controls the finite dimensional kernel described in \eqref{eq:harmonic-kernel}, while the restriction of $\mathcal J_0$ to its orthogonal complement is invertible. Because $\alpha_j-1$ and $\|V_{\alpha_j}\|_{C^{2,\gamma}}$ both tend to 0, substitution of \eqref{eq:Jacobi-remainder} into \eqref{eq:Jacobi-estimate} allows the term containing the stronger norm to be absorbed, yielding
\[
 \norm{V_{\alpha_j}}_{C^{2,\gamma}}\leq2C_2\delta_{\alpha_j}
\]
for all sufficiently large $j$. Consequently, the normalized fields $W_{\alpha_j}:=V_{\alpha_j}/\delta_{\alpha_j}$ are bounded in $C^{2,\gamma}$, and division of \eqref{eq:Jacobi-remainder} by $\delta_{\alpha_j}$ shows that $\mathcal J_0W_{\alpha_j}\to0$ in $C^{0,\gamma}$. For any fixed $0<\gamma'<\gamma$, compactness therefore gives a subsequence converging in $C^{2,\gamma'}$ to a field $W$ satisfying
\[
 \mathcal J_0W=0,\qquad W\perp\Killing,\qquad \norm{W}_{L^2}=1.
\]
Lemma~\ref{lem:harmonic-Jacobi} then identifies $W$ as $\nabla(\xi\cdot x)$ for some $\xi\in\R^3$, so it remains to exclude these gradient fields by the domain balance identity \eqref{eq:source-balance}.

The first variation of the Dirichlet energy density at the identity map, together with the bound just obtained, gives the uniform expansion
\[
 |du_{\alpha_j}|_g^2-2
 =2\Div V_{\alpha_j}+\mathcal R_{\alpha_j},\qquad
 |\mathcal R_{\alpha_j}|\leq C\norm{V_{\alpha_j}}_{C^1}^2
 \leq C\delta_{\alpha_j}^2.
\]
By \eqref{eq:Omega-derivative}, the weight in the domain balance identity satisfies
\[
 \frac{\Omega'_{\alpha_j,\kappa}(2)}{\alpha_j-1}
 =2\alpha_j(\kappa+2)^{\alpha_j-2}
 \longrightarrow\frac{2}{\kappa+2}>0,\qquad
 \sup_{|s-2|\leq c}|\Omega''_{\alpha_j,\kappa}(s)|
 \leq C(\alpha_j-1)
\]
for a fixed $0<c<1$. Since $\int_{\Sph^2}x\,\dd A_g=0$, we may subtract $\Omega_{\alpha_j,\kappa}(2)$ from the integrand in \eqref{eq:Omega-balance}. Taylor expansion and the preceding estimates of the energy density then yield
\begin{align*}
 0 =\int_{\Sph^2}
 \bigl[\Omega_{\alpha_j,\kappa}(|du_{\alpha_j}|_g^2)
       -\Omega_{\alpha_j,\kappa}(2)\bigr]x\,\dd A_g =2\Omega'_{\alpha_j,\kappa}(2)
       \int_{\Sph^2}(\Div V_{\alpha_j})x\,\dd A_g
       +O\bigl((\alpha_j-1)\delta_{\alpha_j}^2\bigr).
\end{align*}
Dividing by $\Omega'_{\alpha_j,\kappa}(2)\delta_{\alpha_j}$ and using its positive comparison with $(\alpha_j-1)\delta_{\alpha_j}$ gives
\[
 \int_{\Sph^2}(\Div W_{\alpha_j})x\,\dd A_g
 =O(\delta_{\alpha_j}),
\]
so the limiting field $W$ satisfies $\int_{\Sph^2}(\Div W)x\,\dd A_g=0$. On the other hand, the representation $W=\nabla(\xi\cdot x)$ implies
\[
 \int_{\Sph^2}(\Div W)x\,\dd A_g
 =-2\int_{\Sph^2}(\xi\cdot x)x\,\dd A_g
 =-\frac{8\pi}{3}\xi,
\]
which forces $\xi=0$ and contradicts $\norm{W}_{L^2}=1$. This contradiction proves that the original $\alpha_j$-harmonic maps $u_{\alpha_j}$ are rotations for all sufficiently large $j$.
\end{proof}

\begin{theorem}\label{thm:near-one-rigidity}
For every $\kappa>0$ there exists $\alpha_\kappa>1$ such that every smooth $\alpha$-harmonic diffeomorphism of the unit round 2-sphere with $1<\alpha<\alpha_\kappa$ is orthogonal. Equivalently, every nonconstant $\alpha$-harmonic map $u:\Sph^2\to\Sph^n$ with $I_\alpha(u)=1$ and $1 < \alpha < \alpha_\kappa$ is an isometric equatorial embedding.
\end{theorem}
\begin{proof}
If the first assertion were false, there would be $\alpha_j\searrow1$ and nonorthogonal $\alpha_j$-harmonic diffeomorphisms $u_{\alpha_j}$. Proposition~\ref{prop:rotational-limit} gives smooth convergence along a subsequence of $u_{\alpha_j}$ to an orthogonal map. After a reflection, Lemma~\ref{lem:sequential-local} implies that the maps are rotations for all sufficiently large indices, a contradiction. The equivalence follows from Propositions~\ref{prop:index-one-forward} and \ref{prop:automatic-index}.
\end{proof}

\subsection{Morse index and rotation symmetries}\label{sec:full-index-rigidity}
In this subsection, we compare the second variation with the scalar form to bound the sum of the Morse index and nullity by 6 (Proposition~\ref{prop:six-dimensional-bound}), extending the harmonic map estimate in \cite{Karpukhin2021}. For a nonorthogonal diffeomorphism, we then identify the Jacobi kernel with the fields generated by domain and target rotations, obtaining Morse index 0 and nullity 6 in Proposition~\ref{prop:nonround-nondegenerate}.

\begin{proposition}\label{prop:six-dimensional-bound}
Every smooth $\alpha$-harmonic diffeomorphism $u:(\Sph^2,g)\to(\Sph^2,\bar g)$ satisfies
\begin{equation}\label{eq:six-dimensional-bound}
 \ind_{\Sph^2}E_{\alpha,\kappa}(u)
 +\nul_{\Sph^2}E_{\alpha,\kappa}(u)\leq6.
\end{equation}
\end{proposition}
\begin{proof}
By Proposition~\ref{prop:automatic-index}, the scalar operator $J_{\alpha,u}$ has exactly 1 negative eigenvalue with multiplicity 1 and
\[
 \ker J_{\alpha,u}=\Span\{u^1,u^2,u^3\}.
\]
Let $h>0$ be the eigenfunction corresponding to the negative eigenvalue, normalized in $L^2(\Sph^2,(\kappa+|du|_g^2)^{\alpha-1}\dd A_g)$. On the orthogonal complement of $h$ in this weighted space, the scalar spectral decomposition gives $Q_{\alpha,u}(f,f)\geq0$ for $f\in W^{1,2}(\Sph^2) \cap h^{\perp}$, with equality precisely when $f\in\ker J_{\alpha,u}$.

Let $\mathcal Z$ be the direct sum of the negative and zero eigenspaces of the Jacobi operator of $E_{\alpha,\kappa}$ on sections of $u^{-1}T\Sph^2$, realized with the same weighted $L^2$ inner product. Since these eigenspaces represent the negative and zero directions of the second variation, the Hessian is nonpositive on $\mathcal Z$ and
\[
 \dim\mathcal Z=\ind_{\Sph^2}E_{\alpha,\kappa}(u)
                  +\nul_{\Sph^2}E_{\alpha,\kappa}(u).
\]
To compare $\mathcal Z$ with the scalar negative eigenspace, define
\begin{equation}\label{eq:ground-state-projection}
 \mathcal P:\mathcal Z\longrightarrow\R^3,
 \qquad
 \mathcal P(V)=\int_{\Sph^2}(\kappa+|du|_g^2)^{\alpha-1}hV\,\dd A_g.
\end{equation}
For $V\in\ker\mathcal P$, each ambient component $V^i$ is orthogonal to $h$, $1 \leq i \leq 3$, so the scalar forms  $Q_{\alpha,u}$ are nonnegative and \eqref{eq:hessian} yields
\begin{equation}\label{eq:saturation-calculation}
\begin{aligned}
 0\geq D^2E_{\alpha,\kappa}(u)[V,V]
 =\sum_{i=1}^3Q_{\alpha,u}(V^i,V^i) +2\alpha(\alpha-1)\int_{\Sph^2}
 (\kappa+|du|_g^2)^{\alpha-2}
 \ip{du}{dV}_g^2\,\dd A_g\geq0.
\end{aligned}
\end{equation}
All the nonnegative terms must therefore vanish, and the scalar equality characterization above implies that every $V^i$ belongs to $\Span\{u^1,u^2,u^3\}$. Consequently, $V=Bu$ for a constant real $3\times3$ matrix $B$, while the tangency condition $\langle V,u\rangle=0$ gives $u^TBu=0$. Surjectivity of $u$ then implies $y^TBy=0$ for every $y\in\Sph^2$, which forces the symmetric part of $B$ to vanish and hence $B\in\mathfrak{so}(3)$.

We have thus shown that $\ker\mathcal P\subset\{Bu:B\in\mathfrak{so}(3)\}$, and this latter space has dimension 3 as surjectivity of $u$ makes $B\mapsto Bu$ injective. Since $\mathcal P$ takes values in $\R^3$, the rank-nullity formula in linear algebra now gives
\[
 \dim\mathcal Z
 =\dim\ker\mathcal P+\operatorname{rank}\mathcal P
 \leq3+3=6,
\]
which proves \eqref{eq:six-dimensional-bound}.
\end{proof}

\begin{lemma}\label{lem:continuous-isotropy}
Let $u:\Sph^2\to\Sph^2$ be a smooth $\alpha$-harmonic diffeomorphism. If there are $A,B\in\mathfrak{so}(3)$, not both 0, such that
\begin{equation}\label{eq:isotropy-equation}
 du_x(Ax)=Bu(x)\qquad \forall x\in\Sph^2,
\end{equation}
then $u\in O(3)$.
\end{lemma}
\begin{proof}
Since \eqref{eq:isotropy-equation} and surjectivity of $u$ force $B=0$ whenever $A=0$, while invertibility of $du$ forces $A=0$ whenever $B=0$, the hypothesis that the matrices are not both 0 implies that both are nonzero. For each $x\in\Sph^2$, the chain rule and \eqref{eq:isotropy-equation} show that the curve $t\mapsto u(e^{tA}x)$ satisfies
\[
 \frac{d}{dt}u(e^{tA}x)
 =du_{e^{tA}x}(Ae^{tA}x)
 =Bu(e^{tA}x).
\]
Uniqueness for this linear ordinary differential equation, with initial value $u(x)$, gives
\[
 u(e^{tA}x)=e^{tB}u(x)
 \qquad\text{for every }t\in\R\text{ and }x\in\Sph^2.
\]

This identity conjugates the two rotation flows by the diffeomorphism $u$, so $e^{tA}=\Id$ if and only if $e^{tB}=\Id$. The flows thus have the same least positive period and consequently the same angular speed. After independent rotations of the domain and target $2$-sphere coordinates, which align their oriented rotation axes with the vertical axis, and a common rescaling of $t$, both flows become the standard rotations $R_\varphi$ about the vertical axis. Retaining the notation $u$ for the transformed map, we obtain $ u(R_\varphi x)=R_\varphi u(x)$ for every $\varphi\in\R$ and $x\in\Sph^2$. This equivariance relation implies that $u$ maps the two poles onto the two poles of 2-spheres. If $u$ interchanges them, we postcompose with the target reflection $(y^1,y^2,y^3)\mapsto(y^1,y^2,-y^3)$,  so that after this normalization $u$ fixes both poles. On the complement of the poles, the polar representation and the smooth real lift of the longitude described before Proposition~\ref{prop:equivariant} give the form \eqref{eq:twisted-map}, with $0<f(r)<\pi$ for $0<r<\pi$, $f(0)=0$ and $f(\pi) = \pi$. In these polar coordinates, the determinant for the differential of $(r,\theta)\mapsto(f(r),\theta+\gamma(r))$ cannot vanish as $u$ is a diffeomorphism. Hence $f'$ has a constant nonzero sign on $(0,\pi)$, and the endpoint values force $f'$ to be positive. The normalization of original $u$ is therefore a smooth twisted $1$-equivariant $\alpha$-harmonic diffeomorphism to which Proposition~\ref{prop:equivariant} applies, yielding $u\in O(3)$. Since all normalizations were orthogonal changes of domain or target coordinates, the original $\alpha$-harmonic diffeomorphism also belongs to $O(3)$.
\end{proof}

\begin{proposition}\label{prop:nonround-nondegenerate}
Let $u:\Sph^2\to\Sph^2$ be a smooth $\alpha$-harmonic diffeomorphism with $u\notin O(3)$. Then its Morse index and nullity satisfy
\begin{equation}\label{eq:nonround-nondegenerate}
 \ind_{\Sph^2}E_{\alpha,\kappa}(u)=0,
 \qquad \nul_{\Sph^2}E_{\alpha,\kappa}(u)=6.
\end{equation}
Moreover, its Jacobi kernel is
\begin{equation}\label{eq:symmetry-kernel}
 \{Bu+du(Ax):A,B\in\mathfrak{so}(3)\},
\end{equation}
which is the tangent space at $u$ to its orbit under domain and target rotations.
\end{proposition}
\begin{proof}
Since both metrics on $\S^2$ are round, domain and target rotations preserve $E_{\alpha,\kappa}$, so the maps $x\mapsto e^{tB}u(e^{tA}x)$ are $\alpha$-harmonic for every $A,B\in\mathfrak{so}(3)$ and $t\in\R$. Differentiating their Euler-Lagrange equation at $t=0$ shows that $Bu+du(Ax)$ is a Jacobi field, and hence the space in \eqref{eq:symmetry-kernel} is contained in the Jacobi kernel of $u$.

If $Bu+du(Ax)=0$, then $du_x(Ax)=(-B)u(x)$ for every $x\in\Sph^2$, so Lemma~\ref{lem:continuous-isotropy} and the assumption $u\notin O(3)$ imply $A=B=0$. Consequently, the linear map $(A,B)\mapsto Bu+du(Ax)$ is injective on $\mathfrak{so}(3)\oplus\mathfrak{so}(3)$, and its image has dimension 6. Combining this resulting lower bound on the nullity with Proposition~\ref{prop:six-dimensional-bound}, we obtain
\[
 6\leq\nul_{\Sph^2}E_{\alpha,\kappa}(u)
 \leq\ind_{\Sph^2}E_{\alpha,\kappa}(u)
       +\nul_{\Sph^2}E_{\alpha,\kappa}(u)\leq6.
\]
This proves \eqref{eq:nonround-nondegenerate}, while equality of dimensions shows that the space in \eqref{eq:symmetry-kernel} exhausts the Jacobi kernel of $u$.
\end{proof}

\subsection{Uniform compactness on closed exponent intervals}\label{subsec:uniform-exponent-compactness}
In this subsection, we establish a uniform gradient bound independent of the target dimension for $\alpha$-harmonic maps with bounded $\alpha$-energy and exponents in a compact subinterval of $(1,\infty)$ (Lemma~\ref{lem:uniform-gradient}). Higher regularity gives smooth compactness for a fixed target sphere; for diffeomorphisms of the unit round $2$-sphere, we use preservation of degree and lower semicontinuity of the scalar index to show that the limit remains a diffeomorphism, as stated in Corollary~\ref{cor:compact-diffeo-branch}.

\begin{lemma}\label{lem:uniform-gradient}
Let $(\Sigma,g)$ be a closed Riemannian surface, and fix $\kappa>0$,
$1<\alpha_-\leq\alpha_+<\infty$ and $\Lambda<\infty$. There is a constant
$C=C(\Sigma,g,\kappa,\alpha_-,\alpha_+,\Lambda)$, independent of $n$, such
that every smooth $\alpha$-harmonic map $u:\Sigma\to\Sph^n$ satisfying $ \alpha\in[\alpha_-,\alpha_+]$ and $E_{\alpha,\kappa}(u)\leq\Lambda$ admits $\norm{du}_{L^\infty(\Sigma)}\leq C$. For fixed $n$, these maps
also have uniform bounds in every $C^k$ norm, with constants that may
depend additionally on $n$ and $k$. Consequently, every sequence
$u_{\alpha_j}$ in this family with $\alpha_j\to\alpha_*$ has a
subsequence converging smoothly to an $\alpha_*$-harmonic map $u_*:\Sigma\to\Sph^n$.
\end{lemma}
\begin{proof}
Writing $e=|du|_g^2$, we first derive a scalar differential inequality
whose coefficients can be controlled independently of $n$. With
$\tau(u)=\tr_g\nabla du$ denoting the tension field of $u$, \eqref{eq:drift}
gives
\[
 \tau(u)=-\frac{\alpha-1}{\kappa+e}du(\nabla e).
\]
Differentiating the pullback metric and using the symmetry of
$\nabla du$, we obtain
\begin{equation}\label{eq:div-pullback-metric}
 \Div_g(u^*\bar g)
 =\frac12\,de-\frac{\alpha-1}{\kappa+e}
       \ip{du(\nabla e)}{du},
\end{equation}
where the last term is the 1-form
$X\mapsto\ip{du(\nabla e)}{du(X)}$. With the convention
$\Delta_g=-\Div_g\nabla$ fixed earlier, the Bochner identity for the
unit round target sphere reads
\[
 -\frac12\Delta_g e
 =|\nabla du|_g^2+\ip{\operatorname{Ric}_g}{u^*\bar g}_g
  -e^2+|u^*\bar g|_g^2+\ip{\nabla\tau(u)}{du}_g.
\]
By the previous formula for tension field and $\frac12\,de(X)=\ip{(\nabla du)(X,\cdot)}{du}_g$, we get
\begin{align*}
 \ip{\nabla\tau(u)}{du}_g
 ={}-\frac{\alpha-1}{\kappa+e}
       \ip{\Hess_g e}{u^*\bar g}_g
     -\frac{\alpha-1}{2(\kappa+e)}|\nabla e|_g^2 +\frac{\alpha-1}{(\kappa+e)^2}|du(\nabla e)|^2.
\end{align*}
To combine these terms in divergence form, let $\mathcal B$ be the
symmetric endomorphism of $T\Sigma$ determined by
\[
 g(\mathcal BX,Y)
 =g(X,Y)+\frac{2(\alpha-1)}{\kappa+e}
                 \ip{du(X)}{du(Y)}.
\]
Substituting the preceding formulas and
\eqref{eq:div-pullback-metric} into the divergence product rule yields
\begin{equation}\label{eq:alpha-Bochner-divergence}
\begin{aligned}
 \frac12\Div_g\bigl((\kappa+e)^{\alpha-1}\mathcal B\nabla e\bigr)&=(\kappa+e)^{\alpha-1}
  \bigl(|\nabla du|_g^2+\ip{\operatorname{Ric}_g}{u^*\bar g}_g
        -e^2+|u^*\bar g|_g^2\bigr)\\
 &\quad+\frac{\alpha-1}{2}(\kappa+e)^{\alpha-2}|\nabla e|_g^2.
\end{aligned}
\end{equation}
Since $|du(X)|^2\leq e|X|_g^2$ and $\alpha\leq\alpha_+$, the
endomorphism $\mathcal{B}$ satisfies
\begin{equation}\label{eq:Bochner-ellipticity}
 \Id\leq\mathcal B\leq(2\alpha_+-1)\Id.
\end{equation}
Choose $C_g\geq0$ such that $\operatorname{Ric}_g\geq-C_gg$.
Using $\nabla(\kappa+e)^\alpha
=\alpha(\kappa+e)^{\alpha-1}\nabla e$
in \eqref{eq:alpha-Bochner-divergence}, together with the  Ricci
curvature lower bound, we obtain
\begin{equation}\label{eq:F-subsolution}
 -\Div_g\bigl(\mathcal B\nabla(\kappa+e)^\alpha\bigr)
 \leq2\alpha(e+C_g)(\kappa+e)^\alpha.
\end{equation}
Both the ellipticity constants and the curvature bound used here
depend only on the fixed exponent interval $[\alpha_-,\alpha_+]$ and the domain metric $g$.

We next apply Moser iteration to \eqref{eq:F-subsolution}, setting
$F=(\kappa+e)^\alpha$ and letting $C$ denote constants depending only
on the data in the statement of Lemma \ref{lem:uniform-gradient}, independently of $n$. The uniform  $\alpha$-energy bound implies
\[
 \norm{F}_{L^1(\Sigma)}\leq2\Lambda,
 \qquad
 \norm{2\alpha(e+C_g)}_{L^{\alpha_-}(\Sigma)}\leq C.
\]
Testing \eqref{eq:F-subsolution} against $F^{s-1}$, for $s\geq2$,
and applying \eqref{eq:Bochner-ellipticity} and H\"older's
inequality therefore gives
\begin{equation}\label{eq:Moser-first}
\begin{aligned}
 \frac{4(s-1)}{s^2}\norm{\nabla F^{s/2}}_{L^2}^2
 \leq\int_\Sigma2\alpha(e+C_g)F^s\dd A_g \leq C\norm{F^{s/2}}_{L^{2\alpha_-/(\alpha_-1)}}^2.
\end{aligned}
\end{equation}
By the Gagliardo–Nirenberg interpolation inequality and $\varepsilon$-Young inequality on the fixed closed surface, for
$0<\varepsilon\leq1$, we have 
\[
 \norm{F^{s/2}}_{L^{2\alpha_-/(\alpha_-1)}}^2
 \leq\varepsilon\norm{\nabla F^{s/2}}_{L^2}^2
 +C\bigl(1+\varepsilon^{-1/(\alpha_-1)}\bigr)
       \norm{F^{s/2}}_{L^2}^2.
\]
After increasing the constant $C$ in \eqref{eq:Moser-first} to be at
least 1, we choose $\varepsilon=2(s-1)/(Cs^2)$ and absorb the
gradient term into the left-hand side, obtaining
\[
 \norm{\nabla F^{s/2}}_{L^2}^2
 \leq Cs^{\alpha_-/(\alpha_-1)}\norm{F^{s/2}}_{L^2}^2.
\]
Combining this estimate with
$W^{1,2}(\Sigma)\hookrightarrow L^4(\Sigma)$ yields
\[
 \norm{F}_{L^{2s}}
 \leq\bigl(Cs^{\alpha_-/(\alpha_-1)}\bigr)^{1/s}
       \norm{F}_{L^s}.
\]
Iteration over $s=2^i$, $i\geq1$, gives
$\norm{F}_{L^\infty}\leq C\norm{F}_{L^2}$, as the product
of the iteration constants converges. As each map $u$ is smooth on the closed
surface, $\norm{F}_{L^\infty}$ is finite, so interpolation with
the energy bound gives
\[
 \norm{F}_{L^\infty}
 \leq C\norm{F}_{L^2}
 \leq C\norm{F}_{L^\infty}^{1/2}\norm{F}_{L^1}^{1/2},
 \qquad
 \norm{F}_{L^\infty}\leq C\norm{F}_{L^1}\leq C\Lambda.
\]
Since $\alpha$ lies in the prescribed compact interval, this proves
the required bound for $\norm{du}_{L^\infty}$, with a constant $C$ independent of $n$.

With the uniform gradient bound established, we fix $n$ and obtain higher
derivative estimates for $u$ by differentiating \eqref{eq:sphere} in local
coordinates. We use $i,j,k,l,\ell$ for
domain indices and $a,b$ for ambient target components, ranging
from 1 to 2 and from 1 to $n+1$, respectively, and sum repeated
contracted indices. For a fixed coordinate direction $\ell$, the
map $w=\partial_\ell u$ satisfies
\begin{equation}\label{eq:differentiated-system}
 \partial_i\bigl(A_{ab}^{ij}\partial_jw^b\bigr)=\partial_iG_a^i,
\end{equation}
where
\[
 A_{ab}^{ij}
 =\sqrt{\det g}\,(\kappa+e)^{\alpha-1}
 \left(g^{ij}\delta_{ab}
 +\frac{2(\alpha-1)}{\kappa+e}
    g^{ik}g^{jl}\partial_ku^a\partial_lu^b\right)
\]
and differentiation of the metric coefficients gives
\begin{align*}
 G_a^i
 ={}&-(\kappa+e)^{\alpha-1}
 \biggl[\delta_\ell^i\sqrt{\det g}\,e\,u^a
       +\partial_\ell\bigl(\sqrt{\det g}\,g^{ij}\bigr)
          \partial_ju^a + \frac{\alpha-1}{\kappa+e}\sqrt{\det g}\,
          g^{ij}(\partial_\ell g^{kl})
          \ip{\partial_ku}{\partial_lu}\partial_ju^a
 \biggr].
\end{align*}
Thus $A$ and $G$ depend only on $u$, $Du$, $\alpha$ and the fixed
metric, and the gradient bound makes both uniformly bounded.
Moreover, for every matrix $\xi=(\xi_i^a)$,
\begin{equation}\label{eq:system-coercivity}
\begin{aligned}
 A_{ab}^{ij}\xi_i^a\xi_j^b
 =\sqrt{\det g}\,(\kappa+e)^{\alpha-1}
 \left(g^{ij}\xi_i\cdot\xi_j
 +\frac{2(\alpha-1)}{\kappa+e}
       \bigl(g^{ij}\partial_i u\cdot\xi_j\bigr)^2\right) \geq c|\xi|^2,
\end{aligned}
\end{equation}
where $c>0$ is uniform on each fixed chart and exponent interval.

To obtain a H\"older semi-norm bound for $w$, let $B_r\subset B_{2r}$ be
concentric coordinate balls and test
\eqref{eq:differentiated-system} with
$\zeta^2(w-w_{B_{2r}\setminus B_r})$, where the subscript denotes
the average over the annulus and $\zeta$ is a smooth cutoff equal to 1
on $B_r$, supported in $B_{2r}$, with $|D\zeta|\leq C/r$.
By Young's inequality and the Poincar\'e inequality on the
annulus, we have 
\begin{align*}
 \int_{B_r}|Dw|^2\dd x
 &\leq\frac{C}{r^2}\int_{B_{2r}\setminus B_r}
               |w-w_{B_{2r}\setminus B_r}|^2\dd x+Cr^2 \leq C\int_{B_{2r}\setminus B_r}|Dw|^2\dd x+Cr^2.
\end{align*}
The usual Caccioppoli estimate and the uniform bound for $w$
control the initial integral on each fixed interior ball.
The Widman's hole-filling iteration argument applied to the last inequality yields
\[
 \int_{B_r}|Dw|^2\dd x\leq Cr^{2\beta}
\]
for some $\beta\in(0,1)$, uniformly for centers in a smaller chart.
Since the domain has dimension 2, Morrey's lemma gives a uniform
$C^{0,\beta}$ bound for $w$, and hence a uniform $C^{1,\beta}$
bound for $u$.

It follows from their displayed formulas that $A$ and $G$ are now
uniformly bounded in $C^{0,\beta}$, so interior Schauder estimates for
\eqref{eq:differentiated-system} give a uniform $C^{1,\beta}$
bound for $w$ on smaller charts. Repeated differentiation and the
same estimates yield bounds for every derivative of $u$, and a
finite cover of $\Sigma$ gives the asserted global $C^k$ norm bounds.
For a sequence with $\alpha_j\to\alpha_*$, by Arzel\`a-Ascoli theorem and
a diagonal choice of subsequence, we get the smooth convergence of $u_{\alpha_j}$, while passage to \eqref{eq:sphere} shows that the limit $u_*$  is $\alpha_*$-harmonic.
\end{proof}

\begin{corollary}\label{cor:compact-diffeo-branch}
Fix $\kappa>0$ and $1<\alpha_-\leq\alpha_+<\infty$. If $\alpha_j\in[\alpha_-,\alpha_+]$ converges to $\alpha_*$ and $u_{\alpha_j}:\Sph^2\to\Sph^2$ are $\alpha_j$-harmonic diffeomorphisms, after passing to a subsequence, $u_{\alpha_j}$ converges smoothly to an $\alpha_*$-harmonic diffeomorphism $u_{\alpha_*}$.
\end{corollary}
\begin{proof}
Theorem~\ref{thm:energy-bound} gives a uniform bound for
$E_{\alpha_j,\kappa}(u_{\alpha_j})$ for $\alpha_j \in [\alpha_-,\alpha_+]$, so Lemma~\ref{lem:uniform-gradient} yields, after passing to a subsequence, $u_{\alpha_j}$ converges smoothly to an $\alpha_*$-harmonic map $u_*:\Sph^2\to\Sph^2$. Since each $u_{\alpha_j}$ has degree
$1$ or $-1$, we may choose the subsequence to have constant degree,
which is preserved by smooth convergence. In particular, $u_*$
is nonconstant.

Smooth convergence of the coefficients of the scalar forms
$Q_{\alpha_j,u_{\alpha_j}}$ implies that their restrictions to any
fixed finite dimensional space of smooth functions converge
uniformly on its $L^2$ unit sphere. Hence, whenever
$Q_{\alpha_*,u_*}$ is negative definite on such a subspace,
$Q_{\alpha_j,u_{\alpha_j}}$ is also negative definite on that
subspace for all sufficiently large $j$.
This proves lower semicontinuity of the scalar index, which,
together with Proposition~\ref{prop:automatic-index} and 
nonconstancy of the bound \eqref{eq:basic-index}, gives
\[
 1\leq I_{\alpha_*}(u_*)
 \leq\liminf_{j\to\infty}I_{\alpha_j}(u_{\alpha_j})=1.
\]
Thus $I_{\alpha_*}(u_*)=1$, and
Proposition~\ref{prop:index-one-forward} shows that $u_*$ is a
diffeomorphism.
\end{proof}

\subsection{Local dependence on the exponent and global rigidity}\label{subsec:exponent-global-rigidity}
In this subsection, for fixed $\kappa>0$, we first show that every nonorthogonal $\alpha_*$-harmonic diffeomorphism, with $\alpha_*>1$, admits nearby nonorthogonal $\alpha$-harmonic diffeomorphisms for all $\alpha$ sufficiently close to $\alpha_*$, by applying the implicit function theorem to variations that are $L^2$-orthogonal to the rotational Jacobi fields. The same argument near $O(3)$ shows that every $\alpha$-harmonic map sufficiently close to an orthogonal map is orthogonal, uniformly for exponents near $\alpha_*$, as stated in Lemma~\ref{lem:parameter-continuation}. Together with Corollary~\ref{cor:compact-diffeo-branch}, these results show that the set of exponents admitting nonorthogonal $\alpha$-harmonic diffeomorphisms is both open and relatively closed in $(1,\infty)$. The rigidity near 1 established in Theorem~\ref{thm:near-one-rigidity}, combined with connectedness of $(1,\infty)$, then forces this set to be empty, proving Theorem~\ref{thm:global-diffeo-rigidity}.

\begin{lemma}\label{lem:parameter-continuation}
Fix $\alpha_*>1$ and $\kappa>0$.
\begin{enumerate}
\item If $u_*$ is a nonorthogonal $\alpha_*$-harmonic diffeomorphism, then for every $\alpha$ sufficiently close to $\alpha_*$ there is a nonorthogonal $\alpha$-harmonic diffeomorphism close to $u_*$.
\item There is a neighborhood of $(\alpha_*,O(3))$ in the product of the exponent parameter space $(1,\infty)$ and $C^{2,\gamma}(\Sph^2,\Sph^2)$, $0<\gamma<1$, in which every $\alpha$-harmonic map is orthogonal.
\end{enumerate}
\end{lemma}
\begin{proof}
Fix $0<\gamma<1$. For (1), let $\mathcal K_*$ denote the Jacobi
kernel at $u_*$ identified in Proposition~\ref{prop:nonround-nondegenerate}
and \eqref{eq:symmetry-kernel}. In the exponential chart
$u_V=\exp_{u_*}V$, let $\mathcal F(\alpha,V)$ be the unprojected
$L^2$ gradient of $V\mapsto E_{\alpha,\kappa}(u_V)$, using the
construction in the proof of Proposition~\ref{prop:local-rigidity}.
Here the inner product is fixed on sections of $u_*^{-1}T\Sph^2$.
Since $\kappa>0$, $\mathcal{F}(\alpha,V)$ depends smoothly on $(\alpha,V)$
as a map from $C^{2,\gamma}$ to $C^{0,\gamma}$. Its derivative in
$V$ at $(\alpha_*,0)$ is the self adjoint strongly elliptic Hessian
operator with kernel $\mathcal K_*$. The Fredholm alternative and
Schauder estimates therefore make the derivative of the projected
gradient an isomorphism from $C^{2,\gamma}\cap\mathcal K_*^\perp$
onto $C^{0,\gamma}\cap\mathcal K_*^\perp$. The implicit function
theorem with parameter $\alpha$ gives a smooth family $V(\alpha)$
for $\alpha$ near $\alpha_*$ such that
\[
 V(\alpha_*)=0,\qquad V(\alpha)\perp\mathcal K_*,\qquad
 \operatorname{proj}_{\mathcal K_*^\perp}
 \mathcal F(\alpha,V(\alpha))=0.
\]

To show that these maps satisfy the $\alpha$-harmonic map equation,
choose a basis $K_1,\ldots,K_6$ of $\mathcal K_*$ consisting of
fields generated by domain and target rotations, and let $W_i(V)$
be the velocities of the corresponding rotation curves in exponential
coordinates, so that $W_i(0)=K_i$. Rotation invariance of the  $\alpha$-energy
gives
\[
 \ip{\mathcal F(\alpha,V)}{W_i(V)}_{L^2}=0,
 \qquad 1\leq i\leq6.
\]
The projected equation places
$\mathcal F(\alpha,V(\alpha))$ in $\mathcal K_*$, while the matrix
$(\ip{K_j}{W_i(V)})_{i,j=1}^6$ is the Gram matrix of the chosen
basis at $V=0$ and remains invertible nearby. These identities thus
force $\mathcal F(\alpha,V(\alpha))=0$. Elliptic regularity makes
$u_{V(\alpha)}$ smooth, and, after shrinking the parameter interval,
these maps remain diffeomorphisms outside $O(3)$, since being a
diffeomorphism is open in $C^1$ and $O(3)$ is closed.

For (2), use the target rotation slice and the spaces $X,Y$ constructed
in the proof of Proposition~\ref{prop:local-rigidity}. At $(\alpha_*,0)$, Proposition \ref{prop:local-rigidity} gives an isomorphism from $X$ onto
$Y$, so the implicit function theorem yields a unique small solution $V$ of $\mathcal{F}(\alpha,V) = 0$
in the slice for every $\alpha$ nearby $\alpha_*$. Since the identity map is
$\alpha$-harmonic for every $\alpha > 1$, this solution is always $V=0$.
The rotation normalization used there is independent of the exponent $\alpha$,
hence every $\alpha$-harmonic map sufficiently close to the identity,
with $\alpha$ sufficiently close to $\alpha_*$, is a rotation.
Finally, postcomposition with $R^{-1}$ reduces a map close to any
$R\in O(3)$ to this case while preserving both the $\alpha$-energy and the ambient $C^{2,\gamma}$ norm. The same neighborhood therefore works
along all of $O(3)$, proving (2).
\end{proof}

\begin{theorem}\label{thm:global-diffeo-rigidity}
For every $\alpha>1$ and $\kappa>0$, every smooth $\alpha$-harmonic diffeomorphism $u:\Sph^2\to\Sph^2$ between unit round $2$-spheres belongs to $O(3)$.
\end{theorem}
\begin{proof}
Fix $\kappa>0$ and define
\[
 \mathcal D_\kappa=
 \{\alpha>1:\text{there is a nonorthogonal $\alpha$-harmonic diffeomorphism}
                    \ \Sph^2\to\Sph^2\}.
\]
Part (1) of Lemma~\ref{lem:parameter-continuation} shows that $\mathcal D_\kappa$ is open in $(1,\infty)$. We next prove that the $\mathcal D_\kappa$ is relatively closed in $(1,\infty)$. Let $\alpha_j\in\mathcal D_\kappa$ converge to $\alpha_*>1$, and choose corresponding nonorthogonal $\alpha_j$-harmonic diffeomorphisms $u_{\alpha_j}$. The exponents lie in a compact subinterval of $(1,\infty)$ after finitely many terms are discarded. By Corollary~\ref{cor:compact-diffeo-branch}, there is a subsequence of $u_{\alpha_j}$ converges smoothly to an $\alpha_*$-harmonic diffeomorphism $u_*$. If $u_*$ were orthogonal, part (2) of Lemma~\ref{lem:parameter-continuation} would force $u_{\alpha_j}$ to be orthogonal for all sufficiently large $j$, a contradiction. Thus $u_*$ is nonorthogonal and $\alpha_*\in\mathcal D_\kappa$. Since the interval $(1,\infty)$ is connected, the set $\mathcal D_\kappa$ is both open and relatively closed, whereas Theorem~\ref{thm:near-one-rigidity} excludes every exponent in a nonempty interval adjacent to 1. Therefore $\mathcal D_\kappa$ is empty, completing the proof the Theorem \ref{thm:global-diffeo-rigidity}.
\end{proof}

\begin{corollary}\label{cor:global-scalar-rigidity}
Let $u:\Sph^2\to\Sph^n$, $n\geq2$, be a smooth nonconstant $\alpha$-harmonic map, where $\alpha>1$ and $\kappa>0$. Then
\begin{equation}\label{eq:global-scalar-classification}
 I_\alpha(u)=1
 \quad\Longleftrightarrow\quad
 u(x)=Lx\quad\text{for a linear map }L:\R^3\to\R^{n+1}
 \text{ with }L^TL=\Id.
\end{equation}
In this case
\[
 |du|^2=2,\quad E_{\alpha,\kappa}(u)=2\pi(\kappa+2)^\alpha,
 \quad N_\alpha(u)=3,\quad
 \ind_{\Sph^n}E_{\alpha,\kappa}(u)=n-2,\quad
 \nul_{\Sph^n}E_{\alpha,\kappa}(u)=3n-3.
\]
Moreover, the only positive Codazzi tensor $C$ satisfying ${\bigl((\det C)^{1/(2(\alpha-1))}-\kappa\bigr)}>0$ and $K_{G_C}=1$ in \eqref{eq:GC-intro} is $C=(\kappa+2)^{\alpha-1}\bar g$.
\end{corollary}
\begin{proof}
If $I_\alpha(u)=1$, Proposition~\ref{prop:index-one-forward}
shows that $u$ is an $\alpha$-harmonic diffeomorphism onto an
equatorial $2$-sphere. After identifying this $2$-sphere isometrically
with $\Sph^2$, Theorem~\ref{thm:global-diffeo-rigidity} shows
that $u$ is an isometry onto its image and therefore has the form
in \eqref{eq:global-scalar-classification}. Conversely, every map
of this form is an isometric equatorial embedding, so
Corollary~\ref{cor:equatorial-index} gives scalar index 1 and the
asserted index and nullity values. Since $L^TL=\Id$ implies
$|du|^2=2$, the stated $\alpha$-energy formula follows as well. Finally, let $C$ be the Codazzi tensor satisfying the stated determinant and curvature conditions. Theorem~\ref{thm:reconstruction} associates to $C$
an orientation preserving $\alpha$-harmonic diffeomorphism of the
unit round $2$-sphere. By Theorem~\ref{thm:global-diffeo-rigidity},
this diffeomorphism is a rotation, whose corresponding tensor is
$C=(\kappa+2)^{\alpha-1}\bar g$ via the same correspondence.
\end{proof}

\vskip1cm
\section{Harmonic index estimates and reduction to 2 bubbles}\label{sec:bubble-index}

In this section, we classify the bubble trees arising under the energy and Morse index bounds of Theorem~\ref{thm:main-full-gap}. The Morse index estimates for harmonic maps in Subsection~\ref{subsec:harmonic-index-gap} show that at most 2 nonconstant components can occur, each has harmonic degree 1 and image being an equatorial $2$-sphere. In Subsection~\ref{sec:two-bubble-reduction}, we use the limiting domain balance identity \eqref{eq:limit-balance} to prove that, after passing to a subsequence, either the limiting harmonic maps are isometric equatorial parametrizations for all sufficiently large $j$, or the base map is constant and exactly 2 bubbles form at antipodal points of the domain, as stated in Proposition~\ref{prop:two-bubble-reduction}.

Throughout this section, $(\Sph^2,g)$ is the unit round sphere and $\kappa>0$ is fixed. We consider a sequence of smooth $\alpha_j$-harmonic maps $u_{\alpha_j}:\Sph^2\to\Sph^n$ with $\alpha_j\searrow1$ and uniformly bounded $E_{\alpha_j,\kappa}(u_{\alpha_j})$. After passing to a subsequence as in Theorem~\ref{thm:compactness-input}, we retain the notation $u_\infty$ for the base map and $\omega_1,\ldots,\omega_L$ for the bubbles from Subsection \ref{subsec:diffeo-compactness}. For each harmonic component $v$, we use the scalar form \eqref{eq:scalar-form} and operator \eqref{eq:operator-intro} at $\alpha=1$, with the both scalar and Morse index and nullity conventions of Subsection~\ref{subsec:main-index}. All derivatives and integrals are taken with respect to the metric on the closed domain of $v$. For a harmonic bubble, we may choose any smooth metric in its conformal class on $\Sph^2$, since the scalar form at $\alpha=1$, and hence its scalar and Morse index and nullity, are conformally invariant.

\subsection{Morse index lower semicontinuity and the harmonic spectral gap}\label{subsec:harmonic-index-gap}
In this subsection, we combine Morse index lower semicontinuity under bubble tree convergence, recalled in Theorem~\ref{thm:full-index-input}, with scalar index estimates for harmonic spheres. Lemmas~\ref{lem:harmonic-sphere-index-one} and~\ref{lem:harmonic-scalar-gap} separate the equatorial conformal parametrizations of scalar index 1 from all other nonconstant harmonic spheres, whose scalar index is at least 3.

\begin{theorem}[Morse index lower semicontinuity {\cite{Moore-Ream2019}}]\label{thm:full-index-input}
Fix $n\geq2$ and $\kappa>0$, and let
$u_{\alpha_j}:\Sph^2\to\Sph^n$ be a sequence of smooth
$\alpha_j$-harmonic maps with $\alpha_j\searrow1$ and
$\sup_jE_{\alpha_j,\kappa}(u_{\alpha_j})<\infty$.
After passing to a subsequence as in
Theorem~\ref{thm:compactness-input}, using the same notation
$u_\infty$ and $\omega_1,\ldots,\omega_L$ for the harmonic base
map and the nonconstant harmonic bubbles, there holds
\begin{equation}\label{eq:full-index-lower}
 \ind_{\Sph^n}E(u_\infty)+\sum_{\ell=1}^L\ind_{\Sph^n}E(\omega_\ell)
 \leq\liminf_{j\to\infty}\ind_{\Sph^n}E_{\alpha_j,\kappa}(u_{\alpha_j}).
\end{equation}
\end{theorem}

Let $v:\Sph^2\to\Sph^n$ be a nonconstant harmonic map, and let $\Sph^r\subset\Sph^n$ be the smallest equatorial sphere containing the image $v(\S^2)$. The lifting argument in the proof of Proposition~\ref{prop:index-one-forward} excludes the case $r=1$, so $r\geq2$ and $v$ is linearly full as a map into $\Sph^r$. For $r>2$,
Karpukhin's comparison \cite{Karpukhin2021} gives
$\ind_{\Sph^r}E(v)\geq(r-2)I_1(v)$, while for $r=2$ the same bound follows from $\ind_{\Sph^2}E(v)\geq0$. The normal splitting formula
\eqref{eq:index-splitting}, applied at $\alpha=1$, yields
\begin{equation}\label{eq:harmonic-index-comparison}
 \ind_{\Sph^n}E(v)
 =\ind_{\Sph^r}E(v)+(n-r)I_1(v)
 \geq(n-2)I_1(v).
\end{equation}

\begin{lemma}\label{lem:harmonic-sphere-index-one}
A nonconstant harmonic map $\omega:\Sph^2\to\Sph^n$ has
$I_1(\omega)=1$ if and only if it is a conformal diffeomorphism
onto an equatorial $2$-sphere, in which case $N_1(\omega)=3$.
\end{lemma}
\begin{proof}
Since Proposition~\ref{prop:index-one-forward} is applicable at $\alpha=1$,
the assumption $I_1(\omega)=1$ implies that $\omega$ is a
diffeomorphism onto an equatorial $2$-sphere and that
$N_1(\omega)=3$. The Hopf differential of $\omega$ is holomorphic
and vanishes on $\Sph^2$, so $\omega$ is
conformal and has degree $1$ or $-1$.

Conversely, suppose that $\omega$ is a conformal diffeomorphism
onto an equatorial $2$-sphere, which we identify isometrically
with $(\Sph^2,\bar g)$. Changing variables by $\omega$ and using
conformal invariance gives
\[
 Q_{1,\omega}(f\circ\omega,f\circ\omega)
 =\int_{\Sph^2}\bigl(|\nabla f|_{\bar g}^2-2f^2\bigr)\,\dd A_{\bar g},
 \qquad f\in W^{1,2}(\Sph^2).
\]
The spectrum of the round Laplacian shows that this form has
exactly 1 negative direction, given by the constant functions,
and a kernel of dimension 3 spanned by the coordinate functions.
Thus $I_1(\omega)=1$ and $N_1(\omega)=3$. This completes the proof of Lemma \ref{lem:harmonic-sphere-index-one}.
\end{proof}

\begin{lemma}\label{lem:harmonic-scalar-gap}
Every nonconstant harmonic sphere $v:\Sph^2\to\Sph^n$ satisfies
\begin{equation}\label{eq:harmonic-scalar-gap}
 I_1(v)=1\quad\text{or}\quad I_1(v)\geq3.
\end{equation}
If the harmonic degree $d=E(v)/(4\pi)$ of $v$ is 2, then $I_1(v)=N_1(v)=3$.
\end{lemma}
\begin{proof}
By the results of Calabi \cite{Calabi1967} and Barbosa
\cite{Barbosa1975}, the smallest equatorial target of $v$ is
$\Sph^{2m}$ for some $m\geq1$, and its harmonic degree $d$ is a
positive integer satisfying $d\geq m(m+1)/2$. Viewing $v$ as a
linearly full map into this target, Karpukhin's scalar index
estimate \cite{Karpukhin2021} and the nullity results of
Ejiri-Kotani \cite{EjiriKotani1992} and Kotani
\cite{Kotani1997} give
\begin{equation}\label{eq:harmonic-degree-estimates}
 I_1(v)+N_1(v)\geq2d+2,\qquad
 N_1(v)\text{ is odd},\qquad N_1(v)^2\leq8d+1.
\end{equation}
The coordinate functions of $v$ in the smallest equatorial target $\S^{2m}$ are
linearly independent elements of the scalar kernel, so
$N_1(v)\geq2m+1\geq3$.

If $d=1$, the degree bound forces $m=1$, and $v$ is a conformal
diffeomorphism onto an equatorial $2$-sphere. Hence
Lemma~\ref{lem:harmonic-sphere-index-one} gives $I_1(v)=1$.
If $d=2$, the same degree bound again gives $m=1$, so, after a
target reflection if necessary, $v$ is a holomorphic map of degree
2 onto this sphere. The degree 2 case of Montiel-Ros's index
theorem \cite{MontielRos1991} gives $I_1(v)=3$, while
\eqref{eq:harmonic-degree-estimates} and $N_1(v)\geq3$ force
$N_1(v)=3$. Finally, if $d\geq3$, then
\eqref{eq:harmonic-degree-estimates} yields
\[
 I_1(v)\geq2d+2-\sqrt{8d+1}\geq3.
\]
These cases prove both assertions of Lemma \ref{lem:harmonic-scalar-gap}.
\end{proof}

\subsection{Balance measures and the reduction to 2 bubbles}\label{sec:two-bubble-reduction}
In this subsection, we pass the domain balance identity to the limit and identify the resulting measure in Lemma~\ref{lem:balance-atoms}. Combined with the Morse and scalar index estimates for harmonic $2$-spheres in Subsection \ref{subsec:harmonic-index-gap}, this yields Proposition~\ref{prop:two-bubble-reduction}, which gives the alternative between eventual isometric equatorial parametrizations and exactly 2 equatorial bubbles of degree 1 at antipodal points, when passing to the limit as $\alpha \searrow 1$, for the $\alpha$-harmonic $2$-spheres satisfying the uniform $\alpha$-energy and sharp Morse index bounds in Theorem \ref{thm:main-full-gap}.

\begin{lemma}\label{lem:balance-atoms}
Let $u_{\alpha_j}:\Sph^2\to\Sph^n$ be a sequence of smooth
$\alpha_j$-harmonic maps satisfying $\alpha_j\searrow1$ and
$E_{\alpha_j,\kappa}(u_{\alpha_j})\leq\Lambda$. With the bubble tree
notation of Theorem~\ref{thm:compactness-input}, the nonnegative measures
\begin{equation}\label{eq:neck-balance-measure}
 \dd\mu_{\alpha_j}
 =\frac{\Omega_{\alpha_j,\kappa}(|du_{\alpha_j}|_g^2)}
        {\alpha_j-1}\,\dd A
\end{equation}
have total mass at most $2\Lambda$ and, after passing to a further
subsequence, satisfy
\begin{equation}\label{eq:balance-measure-limit}
\begin{gathered}
 \mu_{\alpha_j}\rightharpoonup
 \Xi_\kappa(|du_\infty|_g^2)\,\dd A
       +\sum_{p\in\mathcal S}b_p\delta_p,\qquad
 \int_{\Sph^2}x\,\Xi_\kappa(|du_\infty|_g^2)\,\dd A
       +\sum_{p\in\mathcal S}b_px(p)=0,
\end{gathered}
\end{equation}
where the convergence is weak convergence of finite Radon measures
and $\Xi_\kappa$ is defined in \eqref{eq:Xi}. The atomic masses obey
\[
 b_p\geq2\sum_{p_\ell=p}E(\omega_\ell)>0,
 \qquad p\in\mathcal S.
\]
In particular, a bubble of harmonic degree 1 at $p$ implies
$b_p\geq8\pi$.
\end{lemma}
\begin{proof}
The definition of $\Omega_{\alpha,\kappa}$ in \eqref{eq:Omega} and its integral representation in \eqref{eq:Omega-integral}, together with \eqref{eq:Xi}, give
\begin{equation}\label{eq:balance-density-bounds}
\begin{aligned}
 0\leq\frac{\Omega_{\alpha_j,\kappa}(s)}{\alpha_j-1}
 \leq s(\kappa+s)^{\alpha_j-1},\qquad
 \frac{\Omega_{\alpha_j,\kappa}(s)}{\alpha_j-1}
 \geq c_{\alpha_j}\Xi_\kappa(s),
 \qquad \text{for }\,\,  s\geq0,
\end{aligned}
\end{equation}
where $c_{\alpha_j}=\alpha_j\min\{1,\kappa^{\alpha_j-1}\}\to1$. The upper bound in \eqref{eq:balance-density-bounds} yields
\[
 \mu_{\alpha_j}(\Sph^2)
 \leq\int_{\Sph^2}(\kappa+|du_{\alpha_j}|_g^2)^{\alpha_j}\,\dd A
 =2E_{\alpha_j,\kappa}(u_{\alpha_j})\leq2\Lambda,
\]
so after passing to a  subsequence $\mu_{\alpha_j}$ converges weakly as finite Radon measures.
Smooth convergence of $u_{\alpha_j}$ on compact subsets of $\Sph^2\setminus\mathcal S$
and \eqref{eq:Xi} identify the diffuse part of the limit, giving the first weak convergence of \eqref{eq:balance-measure-limit} with $b_p\geq0$.

To compare the atomic masses with the bubble energies, observe that
the $\alpha$-energy bound gives
$\int_{\Sph^2}|du_{\alpha_j}|_g^2\,\dd A\leq4\pi+2\Lambda$.
For every measurable set $A\subset\Sph^2$ of positive area $|A|$,
Jensen's inequality therefore yields
\begin{equation}\label{eq:log-nonconcentration}
 \int_A\log\left(1+\frac{|du_{\alpha_j}|_g^2}{\kappa}\right)\,\dd A
 \leq |A|\log\left(1+\frac{4\pi+2\Lambda}{\kappa|A|}\right).
\end{equation}
Since the right-hand side tends to 0 as $|A|\to0$, uniformly in $j$, these
logarithmic densities are uniformly integrable. Together with smooth
convergence away from the finite set $\mathcal S$, this implies their
$L^1$ convergence to $\log(1+|du_\infty|_g^2/\kappa)$, so they contribute
no atomic mass. On the other hand, the energy identity
\eqref{eq:energy-identity} gives
\[
 |du_{\alpha_j}|_g^2\,\dd A
 \rightharpoonup |du_\infty|_g^2\,\dd A
       +2\sum_{\ell=1}^L E(\omega_\ell)\delta_{p_\ell}.
\]
Using the expression for $\Xi_\kappa$ in \eqref{eq:Xi}, we may therefore
pass to the weak limit in the lower bound of
\eqref{eq:balance-density-bounds} and compare the atomic masses.
Since $c_{\alpha_j}\to1$, this gives
$b_p\geq2\sum_{p_\ell=p}E(\omega_\ell)$ for every $p\in\mathcal S$.
Each such point carries a nonconstant bubble, so $b_p>0$, and a bubble of
harmonic degree 1 has energy $4\pi$ and hence gives $b_p\geq8\pi$. Finally, the domain balance identity \eqref{eq:Omega-balance} gives
$\int_{\Sph^2}x\,\dd\mu_{\alpha_j}=0$ for every $j$, and passage to the
weak limit against the position vector $x$ yields the
second identity in \eqref{eq:balance-measure-limit}.
\end{proof}

\begin{proposition}\label{prop:two-bubble-reduction}
Let $n\geq3$, $\alpha_j\searrow1$, and let $u_{\alpha_j}:\Sph^2\to\Sph^n$ be a sequence of nonconstant smooth $\alpha_j$-harmonic maps with
\begin{equation}\label{eq:two-bubble-budget}
 E_{\alpha_j,\kappa}(u_{\alpha_j})\leq\Lambda,\qquad
 \ind_{\Sph^n}E_{\alpha_j,\kappa}(u_{\alpha_j})<3(n-2).
\end{equation}
After passing to a subsequence, either $u_{\alpha_j}$ is an isometric equatorial parametrization for all sufficiently large $j$, or the base $u_\infty$ is constant and there are exactly two equatorial bubbles of degree 1 at distinct antipodal domain points. In the alternative with two bubbles, there holds
\begin{equation}\label{eq:antipodal-energy-measure}
 \frac12|du_{\alpha_j}|^2\dd A\rightharpoonup
 4\pi\delta_p+4\pi\delta_{-p}.
\end{equation}
\end{proposition}
\begin{proof}
After passing to a subsequence as in
Theorem~\ref{thm:compactness-input}, we combine the Morse index
lower semicontinuity of Theorem~\ref{thm:full-index-input} with
\eqref{eq:harmonic-index-comparison} and the assumed bound
\eqref{eq:two-bubble-budget} to obtain
\begin{equation}\label{eq:harmonic-tree-budget}
 (n-2)\left(I_1(u_\infty)+\sum_{\ell=1}^L I_1(\omega_\ell)\right)
 \leq\liminf_{j\to\infty}
       \ind_{\Sph^n}E_{\alpha_j,\kappa}(u_{\alpha_j})
 \leq3(n-2)-1,
\end{equation}
where $I_1(u_\infty)=0$ if the base map is constant. There must be at least 1 nonconstant component in the bubble tree limit, since otherwise
Theorem~\ref{thm:compactness-input} would give smooth convergence of $u_{\alpha_j}$
to a constant on all of $\Sph^2$. The images of $u_{\alpha_j}$
would then lie in a fixed open hemisphere for large $j$, whereas
integration of \eqref{eq:sphere} gives
\begin{equation}\label{eq:hemisphere-balance}
 \int_{\Sph^2}
 (\kappa+|du_{\alpha_j}|_g^2)^{\alpha_j-1}
 |du_{\alpha_j}|_g^2u_{\alpha_j}\,\dd A=0.
\end{equation}
Taking the inner product with the center of that hemisphere contradicts this identity, as $u_{\alpha_j}$ is nonconstant.
Since $n\geq3$, \eqref{eq:harmonic-tree-budget} and
Lemmas~\ref{lem:harmonic-scalar-gap}
and~\ref{lem:harmonic-sphere-index-one} therefore show that the bubble tree limit has either 1 or 2 nonconstant components, each a conformal diffeomorphism onto an equatorial $2$-sphere with harmonic degree 1 and energy $4\pi$.

Suppose first that $u_\infty$ is nonconstant, so the index bound
allows at most 1 bubble. If a bubble formed at $p$, then
Lemma~\ref{lem:balance-atoms} would give
\[
 \int_{\Sph^2}x\,\Xi_\kappa(|du_\infty|_g^2)\,\dd A
       +b_px(p)=0,
 \qquad b_p\geq8\pi.
\]
However, \eqref{eq:Xi} and $E(u_\infty)=4\pi$ imply
\[
 \left|\int_{\Sph^2}x\,\Xi_\kappa(|du_\infty|_g^2)\,\dd A\right|
 \leq\int_{\Sph^2}\Xi_\kappa(|du_\infty|_g^2)\,\dd A
 <\int_{\Sph^2}|du_\infty|_g^2\,\dd A=8\pi,
\]
where the strict inequality follows from $\kappa>0$ and the
nonconstancy of $u_\infty$. This contradicts the above balance
identity, so $L=0$ and the convergence is smooth on all of
$\Sph^2$. The potentials of the conjugated operators in
\eqref{eq:conjugation} consequently converge uniformly to that
of $J_{1,u_\infty}$. Since
Lemma~\ref{lem:harmonic-sphere-index-one} gives
$I_1(u_\infty)=1$ and $N_1(u_\infty)=3$, the min-max principle
yields
\[
 I_{\alpha_j}(u_{\alpha_j})
       +N_{\alpha_j}(u_{\alpha_j})\leq4
\]
for all sufficiently large $j$. Moreover, $u_\infty$ has 3
linearly independent ambient coordinate functions, and smooth
convergence of $u_{\alpha_j}$ preserves the independence of the corresponding coordinate functions of $u_{\alpha_j}$. Thus
\eqref{eq:sphere} and \eqref{eq:basic-index} give
$N_{\alpha_j}(u_{\alpha_j})\geq3$ and
$I_{\alpha_j}(u_{\alpha_j})\geq1$, forcing
$I_{\alpha_j}(u_{\alpha_j})=1$ for large enough $j$.
Corollary~\ref{cor:global-scalar-rigidity} then shows that
$u_{\alpha_j}$ is an isometric equatorial parametrization of $\S^2$ for
all sufficiently large $j$.

If $u_\infty$ is constant, Lemma~\ref{lem:balance-atoms} gives a
nonzero sum of positive atoms as the limiting balance measure.
The balance identity in \eqref{eq:balance-measure-limit} excludes
support at 1 point, while \eqref{eq:harmonic-tree-budget} allows at most
2 bubbles. Consequently, $L=2$ and the concentration points
$p_1$ and $p_2$ are distinct, so the limiting domain balance
identity becomes
\[
 b_{p_1}x(p_1)+b_{p_2}x(p_2)=0,
 \qquad b_{p_1},b_{p_2}>0.
\]
Since $|x(p_1)|=|x(p_2)|=1$, this identity implies
$b_{p_1}=b_{p_2}$ and $x(p_2)=-x(p_1)$.
Finally, the base map $u_\infty$ is constant and each bubble has energy $4\pi$, so the energy identity \eqref{eq:energy-identity}
gives \eqref{eq:antipodal-energy-measure}.
\end{proof}

\vskip1cm

\section{Sharp Morse index rigidity}\label{sec:two-bubble-necks}

In this section, we complete the proof of
Theorem~\ref{thm:main-full-gap} by establishing a Morse index
lower bound for the configuration of 2 bubbles in
Proposition~\ref{prop:two-bubble-reduction}. The conservation laws
and decay estimates in Subsections~\ref{subsec:neck-conservation}
and~\ref{subsec:neck-decay-hamiltonian} yield the relative expansion
in Subsection~\ref{subsec:relative-neck-expansion}, which we use in
Subsection~\ref{subsec:negative-variations} to construct an additional
$(n-2)$-dimensional subspace on which the second variation is
negative definite.
Finally, Subsection~\ref{subsec:degree-two-sharpness} constructs
degree 2 minimizers whose equatorial inclusions have Morse index
exactly $3(n-2)$ and completes the proof of
Theorem~\ref{thm:main-full-gap}.
Throughout this section, the domain $\Sph^2$ carries the unit round
metric, expressed in the coordinates used in each argument.

\begin{theorem}\label{thm:two-bubble-index}
Fix $n\geq3$, $\kappa>0$, and $\Lambda<\infty$, and let
$u_{\alpha_j}:\Sph^2\to\Sph^n$ be smooth $\alpha_j$-harmonic maps
with $\alpha_j\searrow1$ and
$E_{\alpha_j,\kappa}(u_{\alpha_j})\leq\Lambda$.
Suppose that, with the notation of
Theorem~\ref{thm:compactness-input}, the base map $u_\infty$ is
constant and there are exactly 2 nonconstant bubbles
$\omega_1,\omega_2$ at distinct concentration points
$p_1,p_2$, each an equatorial harmonic sphere of harmonic degree 1.
Then
\begin{equation}\label{eq:two-bubble-index-bound}
 \ind_{\Sph^n}E_{\alpha_j,\kappa}(u_{\alpha_j})\geq3(n-2)
\end{equation}
for all sufficiently large $j$.
\end{theorem}

All constants below are uniform in $j$ and the neck length, once
the parameters, coordinate choices, and cutoffs have been fixed.

\subsection{Cylinder coordinates and rotational conservation laws}\label{subsec:neck-conservation}
In this subsection, we introduce cylinder coordinates on the neck domain and use Lemma~\ref{lem:neck-conductivity} to control the scalar weight $(\kappa + |d u_{\alpha}|^2)^{\alpha -1}$ in  the $\alpha$-harmonic equation and the image of the neck domain.
We then use gnomonic coordinates to derive the rotational
conservation laws in Lemma~\ref{lem:gnomonic-flux}, which provide the zero-flux identities used in the subsequent neck estimates.

Let $x_{\alpha_j}^\ell\to p_\ell$, $\ell=1,2$, be the bubble centers from \eqref{eq:bubble-rescaling}, and fix
$p_0\in\Sph^2\setminus\{p_1,p_2\}$.
Write $\widehat{\mathbb C}=\mathbb C\cup\{\infty\}$ for the
Riemann sphere. For all sufficiently large $j$, let
$\Phi_{\alpha_j}:\widehat{\mathbb C}\to\Sph^2$ be the unique
orientation preserving conformal diffeomorphism satisfying
\[
 \Phi_{\alpha_j}(0)=x_{\alpha_j}^1,\qquad
 \Phi_{\alpha_j}(\infty)=x_{\alpha_j}^2,\qquad
 \Phi_{\alpha_j}(1)=p_0.
\]
Such a parametrization can be obtained by composing an inverse
stereographic projection with a M\"obius transformation of
$\widehat{\mathbb C}$. Since the limiting marked points are
distinct, $\Phi_{\alpha_j}$ converges smoothly to a conformal
diffeomorphism. We use $\zeta$ for the complex coordinate on $\widehat{\mathbb{C}}$, where the local coordinates at the 2 bubble centers are
$\zeta$ and $1/\zeta$, respectively.

We measure the bubble scales in these local coordinates and
retain the notation $r_{\alpha_j}^1,r_{\alpha_j}^2$ from
\eqref{eq:bubble-rescaling}. The transition maps from the earlier
fixed coordinates have derivatives converging to nonzero limits,
so the new bubble scales under $\zeta$ differ from the earlier scales in \eqref{eq:bubble-rescaling} by positive
factors converging to nonzero limits. After absorbing the
resulting rotations of the rescaled variables into the bubble parametrizations, we have
\[
\begin{aligned}
 u_{\alpha_j}\bigl(\Phi_{\alpha_j}(r_{\alpha_j}^1z)\bigr)
     \longrightarrow\omega_1(z),\qquad
 u_{\alpha_j}\bigl(\Phi_{\alpha_j}((r_{\alpha_j}^2z)^{-1})\bigr)
     \longrightarrow\omega_2(z)
\end{aligned}
\]
smoothly on compact subsets of $\mathbb C$, since neither bubble
has a further bubble attached to it. 

On $\mathbb C\setminus\{0\}$, set $\zeta=e^{t+i\theta}$,
where $t\in\R$ and $\theta\in S^1$.
We continue to write $u_{\alpha_j}$ for
$u_{\alpha_j}\circ\Phi_{\alpha_j}$ on the Riemann sphere $\widehat{\mathbb{C}}$,
identified with $\Sph^2$, and set
$g_{\alpha_j}=\Phi_{\alpha_j}^*g$ there. On the cylinder $\R\times S^1$,
both the map $u_{\alpha_j}$ and the metric $g_{\alpha_j}$ are further pulled back by
$(t,\theta)\mapsto e^{t+i\theta}$.
Thus both the map and the domain metric are expressed in the
same cylindrical coordinates, and the equation \eqref{eq:sphere}, the
Hessian \eqref{eq:hessian}, and its Morse index are unchanged.
Under this cylindrical coordinates, the round metric takes the form
\begin{equation}\label{eq:round-cylinder-metric}
 g_{\alpha_j}=\frac{dt^2+d\theta^2}
 {\bigl(A_{\alpha_j}e^t+B_{\alpha_j}e^{-t}
       +C_{\alpha_j}\cos\theta+D_{\alpha_j}\sin\theta\bigr)^2}.
\end{equation}
The coefficients converge as $j\to\infty$, with
$A_{\alpha_j},B_{\alpha_j}$ tending to positive limits, and the
pullback of the round stereographic metric gives $4A_{\alpha_j}B_{\alpha_j} -C_{\alpha_j}^2-D_{\alpha_j}^2=1$.
In particular, the first derivatives of
$\log(g_{\alpha_j})_{tt}$ with respect to $t$ and $\theta$ are
uniformly bounded. Since
$(g_{\alpha_j})_{tt}=(g_{\alpha_j})_{\theta\theta}$ and
$(g_{\alpha_j})_{t\theta}=0$, we have
\[
 |du_{\alpha_j}|_{g_{\alpha_j}}^2
 =\frac{|\partial_tu_{\alpha_j}|^2+
              |\partial_\theta u_{\alpha_j}|^2}
        {(g_{\alpha_j})_{tt}}.
\]
For every $\theta$, the metric coefficient
$(g_{\alpha_j})_{tt}$ increases in $t$ up to
\begin{equation}\label{eq:metric-sign-change}
 t_{c,\alpha_j}=\frac12\log(B_{\alpha_j}/A_{\alpha_j})
\end{equation}
and decreases thereafter.

To centre the connecting annulus in the logarithmic coordinate,
define
\begin{equation}\label{eq:neck-scales}
 \delta_{\alpha_j}=\sqrt{r_{\alpha_j}^1r_{\alpha_j}^2},\qquad
 t_{\alpha_j}^*=\frac12\log\frac{r_{\alpha_j}^1}{r_{\alpha_j}^2},
 \qquad s=t-t_{\alpha_j}^*.
\end{equation}
The rescaled radial variables of the 2 bubbles are
$|\zeta|/r_{\alpha_j}^1$ and
$(r_{\alpha_j}^2|\zeta|)^{-1}$, respectively, and their
reciprocals satisfy
\[
 \frac{r_{\alpha_j}^1}{|\zeta|}
       =\delta_{\alpha_j}e^{-s},\qquad
 r_{\alpha_j}^2|\zeta|=\delta_{\alpha_j}e^s,
\]
The angular variables in the corresponding rescaled coordinates
are $\theta$ and $-\theta$, respectively.

For a fixed $R_0>1$, to be chosen sufficiently large, the annulus
$R_0r_{\alpha_j}^1\leq|\zeta|\leq(R_0r_{\alpha_j}^2)^{-1}$
is represented in $(s,\theta)$  coordinates by
\begin{equation}\label{eq:neck-cylinder-weight}
 N_{\alpha_j}=[-L_{\alpha_j},L_{\alpha_j}]\times S^1,\qquad
 L_{\alpha_j}=-\log(R_0\delta_{\alpha_j}).
\end{equation}
To apply interior elliptic estimates near both ends of $N_{\alpha_j}$,
we first work on the larger cylinder
\[
 \widehat N_{\alpha_j}
 =[-L_{\alpha_j}-10,L_{\alpha_j}+10]\times S^1.
\]
This enlargement amounts to replacing $R_0$ by $R_0e^{-10}$
in the annulus above, and we choose $R_0$ sufficiently large such that the following bound for the oscillation of $u_{\alpha_j}$ to hold on
$\widehat N_{\alpha_j}$.

\begin{lemma}\label{lem:neck-conductivity}
The sequence in Theorem~\ref{thm:two-bubble-index} satisfies
\begin{equation}\label{eq:uniform-neck-conductivity}
 \left\|(\kappa+|du_{\alpha_j}|_{g_{\alpha_j}}^2)^{\alpha_j-1}
       -1\right\|_{L^\infty(\Sph^2)}\longrightarrow0.
\end{equation}
Moreover, for every $\nu>0$, one can fix $R_0$ sufficiently
large, independently of $j$, so that
$u_{\alpha_j}(\widehat N_{\alpha_j})$ lies in the geodesic ball
of radius $\nu$ centred at the constant value of $u_\infty$,
for all sufficiently large $j$.
\end{lemma}
\begin{proof}
By the normalization \eqref{eq:kappa-normalization}, a fixed
$\kappa>0$ does not change Li-Zhu's argument
\cite[Lemma~4.1]{LiZhu2019}, applied separately near the
2 concentration points, and \eqref{eq:uniform-neck-conductivity}
follows directly. The assertion about the neck image follows from the no-neck property
in Theorem~\ref{thm:compactness-input}, after choosing $R_0$
sufficiently large to account for the fixed enlarged neck region $\widehat{N}_{\alpha_j}$.
\end{proof}

Let $q\in\Sph^n$ be the constant value of $u_\infty$.
The no-neck property implies
$u_{\alpha_j}(t_{\alpha_j}^*,0)\to q$, so, after taking suitable target rotations
converging to the identity, we may assume that
$u_{\alpha_j}(t_{\alpha_j}^*,0)=q$ for every $j$.
On the hemisphere $\langle u_{\alpha_j},q\rangle>0$, define
the gnomonic coordinate by
$z_{\alpha_j}=u_{\alpha_j}/\langle u_{\alpha_j},q\rangle-q$.
The central projection of $u_{\alpha_j}$ onto the affine tangent
plane $q+q^\perp$ is then $q+z_{\alpha_j}$, and the inverse
coordinate map is
\begin{equation}\label{eq:gnomonic-coordinate}
 u_{\alpha_j}=\frac{q+z_{\alpha_j}}{\sqrt{1+|z_{\alpha_j}|^2}},\qquad
 z_{\alpha_j}\in q^\perp,\qquad z_{\alpha_j}(0,0)=0.
\end{equation}
By Lemma~\ref{lem:neck-conductivity}, the image of the neck domain lies
in this hemisphere for $R_0$ sufficiently large. Increasing
$R_0$ further if necessary gives
\begin{equation}\label{eq:gnomonic-smallness}
 |z_{\alpha_j}|\leq\nu\quad\hbox{on }\widehat N_{\alpha_j}.
\end{equation}
Under the coordinates on $\widehat{N}_{\alpha_j}$, we use $D=(\partial_s,\partial_\theta)$,
$\operatorname{div}_0$, and
$\Delta_0=\partial_s^2+\partial_\theta^2$ to denote the flat
derivative, divergence, and Laplacian, respectively. For $s\in\R$ and $r>0$, write
$Q_r(s)=(s-r,s+r)\times S^1$.

\begin{lemma}\label{lem:gnomonic-flux}
The gnomonic coordinate function $z_{\alpha_j}$ satisfies the  equation
\begin{equation}\label{eq:gnomonic-conductivity}
 \operatorname{div}_0\left({\frac{\bigl(\kappa+|du_{\alpha_j}|_{g_{\alpha_j}}^2\bigr)^{\alpha_j-1}}{1+|z_{\alpha_j}|^2}}Dz_{\alpha_j}\right)=0.
\end{equation}
For every circle in the cylinder $\widehat{N}_{\alpha_j}$, there holds
\begin{equation}\label{eq:gnomonic-fluxes}
\begin{aligned}
 \int_{S^1}\frac{(\kappa+|du_{\alpha_j}|_{g_{\alpha_j}}^2)^{\alpha_j-1}}{1+|z_{\alpha_j}|^2}(z_{\alpha_j})_s\dd\theta=0,\qquad \int_{S^1}\frac{(\kappa+|du_{\alpha_j}|_{g_{\alpha_j}}^2)^{\alpha_j-1}}{1+|z_{\alpha_j}|^2}z_{\alpha_j}\wedge(z_{\alpha_j})_s\dd\theta=0.
\end{aligned}
\end{equation}
Here $v\wedge w:=vw^T-wv^T$, which is an antisymmetric endomorphism of the ambient Euclidean space $\R^{n+1}$.
\end{lemma}
\begin{proof}
For $1\leq i,l\leq n+1$, multiply the equation for
$u_{\alpha_j}^l$ in \eqref{eq:sphere} by $u_{\alpha_j}^i$
and subtract the equation for $u_{\alpha_j}^i$ multiplied by
$u_{\alpha_j}^l$. The terms involving
$|du_{\alpha_j}|_{g_{\alpha_j}}^2u_{\alpha_j}$ cancel, as do
the mixed products of first derivatives, so the resulting
identities can be written as
\[
 \operatorname{div}_{g_{\alpha_j}}\!\left(
 (\kappa+|du_{\alpha_j}|_{g_{\alpha_j}}^2)^{\alpha_j-1}
 u_{\alpha_j}\wedge\nabla_{g_{\alpha_j}}u_{\alpha_j}
 \right)=0.
\]
Each circle $\{s\}\times S^1$ bounds a disk in $\Sph^2$
on which both $u_{\alpha_j}$ and the domain metric are smooth.
Applying the divergence theorem  on
this disk gives
\[
 \int_{S^1}
 (\kappa+|du_{\alpha_j}|_{g_{\alpha_j}}^2)^{\alpha_j-1}
 u_{\alpha_j}\wedge\partial_su_{\alpha_j}\,\dd\theta=0,
\]
since the factors from the unit normal and the boundary
length cancel in conformal coordinates.

Since $g_{\alpha_j}$ is conformal to $\dd s^2+\dd\theta^2$
on the neck $\widehat{N}_{\alpha_j}$, the divergence identity above is equivalent to
$\operatorname{div}_0\bigl(
(\kappa+|du_{\alpha_j}|_{g_{\alpha_j}}^2)^{\alpha_j-1}
u_{\alpha_j}\wedge Du_{\alpha_j}\bigr)=0$.
Substituting \eqref{eq:gnomonic-coordinate} gives
\[
\begin{aligned}
 (\kappa+|du_{\alpha_j}|_{g_{\alpha_j}}^2)^{\alpha_j-1}
 u_{\alpha_j}\wedge Du_{\alpha_j} =
 \frac{(\kappa+|du_{\alpha_j}|_{g_{\alpha_j}}^2)^{\alpha_j-1}}
      {1+|z_{\alpha_j}|^2}
 \bigl(q\wedge Dz_{\alpha_j}
          +z_{\alpha_j}\wedge Dz_{\alpha_j}\bigr).
\end{aligned}
\]
Since $q$ is constant and $z_{\alpha_j},Dz_{\alpha_j}$ take
values in $q^\perp$, we have
$(q\wedge Dz_{\alpha_j})q=-Dz_{\alpha_j}$ and
$(z_{\alpha_j}\wedge Dz_{\alpha_j})q=0$.
Applying the flat divergence equation to $q$ therefore gives
\eqref{eq:gnomonic-conductivity}. Applying the same calculation
to the above vanishing circle integral gives the first equality
in \eqref{eq:gnomonic-fluxes}, and substituting this equality back
into that integral eliminates the term containing $q$
gives the second equality of \eqref{eq:gnomonic-fluxes}.
\end{proof}

\subsection{Neck estimates and bubble scales}\label{subsec:neck-decay-hamiltonian}
In this subsection, we use the first identity in
\eqref{eq:gnomonic-fluxes} to prove decay of
$\|Dz_{\alpha_j}\|_{L^2(Q_2(s))}$ in
Lemma~\ref{lem:neck-decay}, and then obtain pointwise bounds
for $z_{\alpha_j}$ and $Dz_{\alpha_j}$ and an $L^4$ bound
for $D^2z_{\alpha_j}$ in Lemma~\ref{lem:gnomonic-regularity}.
We also prove an estimate \eqref{eq:epsilon-delta}, relating $\alpha_j-1$ to the product of the 2 bubble scales, in
Lemma~\ref{lem:neck-Hamiltonian}. This comparison result will be used
to estimate $\Delta_0z_{\alpha_j}$ in \eqref{eq:neck-forcing} in Subsection \ref{subsec:relative-neck-expansion}.

\begin{lemma}\label{lem:neck-decay}
If $\nu$ in \eqref{eq:gnomonic-smallness} is sufficiently
small and $R_0$ is sufficiently large, there exists $C>0$,
independent of $j$ and $s$, such that for all sufficiently
large $j$ and every $s\in\R$ satisfying
$Q_2(s)\subset\widehat N_{\alpha_j}$, we have
\begin{equation}\label{eq:neck-L2-decay}
 \|Dz_{\alpha_j}\|_{L^2(Q_2(s))}\leq C\bigl[\delta_{\alpha_j}(e^s+e^{-s})\bigr]^{3/5}.
\end{equation}
\end{lemma}
\begin{proof}
In this proof, to simplify the notation, we write
$u=u_{\alpha_j}$, $z=z_{\alpha_j}$, and $\alpha=\alpha_j$, and set
\[
 a=\frac{(\kappa+|du|_{g_{\alpha_j}}^2)^{\alpha-1}}{1+|z|^2},
 \qquad \widehat L=L_{\alpha_j}+10.
\]
By \eqref{eq:uniform-neck-conductivity} and
\eqref{eq:gnomonic-smallness}, choosing $\nu$ sufficiently small
and then $j$ sufficiently large ensures that
$3/4\leq a\leq5/4$ on $\widehat N_{\alpha_j}$. Let $\bar z(s)=(2\pi)^{-1}\int_{S^1}z(s,\theta)\,\dd\theta$
be the average of $z(s,\theta)$ on the slice $s\times S^1$, and define
\[
 P(s)=\int_{S^1}a\,\langle z,z_s\rangle\,\dd\theta.
\]
Using \eqref{eq:gnomonic-conductivity} and integrating by parts
in $\theta$, we obtain
\[
 P'(s)=\int_{S^1}a|Dz|^2\,\dd\theta
 \geq\frac34\int_{S^1}|Dz|^2\,\dd\theta.
\]
The first identity in \eqref{eq:gnomonic-fluxes} gives
$P(s)=\int_{S^1}a\,\langle z-\bar z(s),z_s\rangle\,\dd\theta$.
Hence the Poincar\'e inequality on $S^1$ yields
\begin{equation}\label{eq:neck-flux-energy}
\begin{aligned}
 |P(s)|
 \leq\frac54\|z-\bar z(s)\|_{L^2(S^1)}
                 \|z_s\|_{L^2(S^1)} \leq\frac58\int_{S^1}|Dz|^2\,\dd\theta
 \leq\frac56P'(s).
\end{aligned}
\end{equation}

For $0\leq r<\widehat L$, define
\[
\begin{aligned}
 F(r)=\int_{-\widehat L+r}^{\widehat L-r} \int_{S^1}a|Dz|^2\,\dd\theta\,\dd s
      =P(\widehat L-r)-P(-\widehat L+r).
\end{aligned}
\]
Since $F'(r)=-P'(\widehat L-r)-P'(-\widehat L+r)$,
applying \eqref{eq:neck-flux-energy} at the 2 boundary circles gives
\[
 F(r)\leq\frac56\bigl(P'(\widehat L-r)+P'(-\widehat L+r)\bigr)
       =-\frac56F'(r).
\]
Integrating this differential inequality gives
$F(r)\leq F(0)e^{-6r/5}$.Differentiating
\eqref{eq:gnomonic-coordinate} gives
\begin{equation}\label{eq:gnomonic-flat-density}
 |Du|^2=(g_{\alpha_j})_{tt}|du|_{g_{\alpha_j}}^2
 =\frac{|Dz|^2}{1+|z|^2}
 -\frac{\sum_{i=s,\theta}\langle z,\partial_i z\rangle^2}
        {(1+|z|^2)^2},
\end{equation}
so $|Du|$ and $|Dz|$ are uniformly mutually controlled on $\widehat N_{\alpha_j}$.

conformal form of $g_{\alpha_j}$ and the uniform $\alpha$-energy bound yields
\[
\begin{aligned}
 F(0)\leq C\int_{\Sph^2}
 (\kappa+|du|_{g_{\alpha_j}}^2)^{\alpha-1}
 |du|_{g_{\alpha_j}}^2\,\dd A_{g_{\alpha_j}} \leq2C E_{\alpha,\kappa}(u)\leq2C\Lambda.
\end{aligned}
\]

For any $s$ with $Q_2(s)\subset\widehat N_{\alpha_j}$,
take $r=\widehat L-|s|-2\geq0$. Then
$Q_2(s)\subset[-\widehat L+r,\widehat L-r]\times S^1$,
and we obtain
\[
 \|Dz\|_{L^2(Q_2(s))}^2
 \leq\frac43F(r)
 \leq C e^{-6(\widehat L-|s|-2)/5}.
\]
Taking square roots and using
$e^{-\widehat L}=e^{-10}R_0\delta_{\alpha_j}$ from
\eqref{eq:neck-cylinder-weight}, together with
$e^{|s|}\leq e^s+e^{-s}$, proves \eqref{eq:neck-L2-decay}, which completes the proof of Lemma \ref{lem:neck-decay}.
\end{proof}

\begin{lemma}\label{lem:gnomonic-regularity}
If $\nu$ is sufficiently small and $R_0$ is sufficiently large
as in Lemma~\ref{lem:neck-decay}, there exists $C>0$, independent
of $j$ and $s$, such that for all sufficiently large $j$ and
every $(s,\theta)\in N_{\alpha_j}$, there holds
\begin{equation}\label{eq:neck-C1-W24}
\begin{aligned}
 |z_{\alpha_j}(s,\theta)|+|Dz_{\alpha_j}(s,\theta)|
 \leq C\bigl[\delta_{\alpha_j}(e^s+e^{-s})\bigr]^{3/5},\qquad
 \|D^2z_{\alpha_j}\|_{L^4(Q_1(s))}
 \leq C\bigl[\delta_{\alpha_j}(e^s+e^{-s})\bigr]^{3/5}.
\end{aligned}
\end{equation}
\end{lemma}
\begin{proof}
In this proof, we write $u=u_{\alpha_j}$, $z=z_{\alpha_j}$,
and $\alpha=\alpha_j$, with all metric coefficients evaluated
at $t=s+t_{\alpha_j}^*$. 
To expand
\eqref{eq:gnomonic-conductivity},  we first observe that 
\[
 \frac{D|du|_{g_{\alpha_j}}^2}
      {\kappa+|du|_{g_{\alpha_j}}^2}
 =\frac{D|Du|^2-|Du|^2D\log(g_{\alpha_j})_{tt}}
        {\kappa(g_{\alpha_j})_{tt}+|Du|^2}.
\]
The derivative of logarithmic of the metric is uniformly bounded by \eqref{eq:round-cylinder-metric}, while
\eqref{eq:gnomonic-flat-density} implies
$|D|Du|^2|\leq C|Dz||D^2z|+C|z||Dz|^3$.
Since the denominator above is bounded below by $\frac14|Dz|^2$
when $\nu\leq1$,
the expanded equation \eqref{eq:gnomonic-conductivity} gives
\begin{equation}\label{eq:gnomonic-nondivergence}
 |\Delta_0z|
 \leq C(\alpha-1)|D^2z|
       +C(\alpha-1)|Dz|+C|z||Dz|^2.
\end{equation}

Fix $s_0$ with $|s_0|\leq L_{\alpha_j}$, and let $c$ be the
average of $z$ on $Q_4(s_0)$. By
\eqref{eq:uniform-neck-conductivity} and
\eqref{eq:gnomonic-smallness}, the coefficient
\[
 a=\frac{(\kappa+|du|_{g_{\alpha_j}}^2)^{\alpha-1}}{1+|z|^2}
\]
in \eqref{eq:gnomonic-conductivity} can be made arbitrarily
close to 1 in $L^\infty(\widehat N_{\alpha_j})$ by choosing
$\nu$ small and then $j$ large. The Poincar\'e inequality and the Sobolev embedding
$W^{1,2}(Q_4(s_0))\subset L^8(Q_4(s_0))$ give
\begin{equation}\label{eq:local-oscillation}
 \|z-c\|_{L^8(Q_4(s_0))}
 \leq C\|Dz\|_{L^2(Q_4(s_0))}.
\end{equation}
The pointwise bound \eqref{eq:gnomonic-smallness} also implies
$|c|\leq\nu$ and hence $|z-c|\leq2\nu$ on
$\widehat N_{\alpha_j}$. Testing
\eqref{eq:gnomonic-conductivity} with $\chi^2(z-c)$, where
$\chi$ is supported in $Q_5(s_0)$ and equals 1 on
$Q_4(s_0)$, therefore yields
$\|Dz\|_{L^2(Q_4(s_0))}\leq C\nu$, since
$Q_5(s_0)\subset\widehat N_{\alpha_j}$.

Writing the equation as
$\Delta_0z=\operatorname{div}_0((1-a)Dz)$, we apply the
interior $W^{1,8}$ estimate for the flat Laplacian with a cutoff
supported in $Q_3(s_0)$. The term containing $(1-a)Dz$
is absorbed by choosing $\|a-1\|_{L^\infty}$ sufficiently
small, and \eqref{eq:local-oscillation} then gives
\[
 \|Dz\|_{L^8(Q_2(s_0))}
 \leq C\bigl(\|z-c\|_{L^8(Q_3(s_0))}
             +\|Dz\|_{L^2(Q_3(s_0))}\bigr)
 \leq C\|Dz\|_{L^2(Q_4(s_0))}.
\]
In the expanded inequality \eqref{eq:gnomonic-nondivergence}, the coefficients of the terms
containing $D^2z$ are bounded by $C(\alpha-1)$.
The interior $W^{2,4}$ estimate therefore allows these terms
to be absorbed for large $j$, giving
\[
\begin{aligned}
 \|D^2z\|_{L^4(Q_{3/2}(s_0))}
 \leq C\bigl(\|z-c\|_{L^4(Q_2(s_0))}
       +(\alpha-1)\|Dz\|_{L^4(Q_2(s_0))}\bigr)
       +C\nu\|Dz\|_{L^8(Q_2(s_0))}^2\leq C\|Dz\|_{L^2(Q_4(s_0))}.
\end{aligned}
\]
Applying the Sobolev embedding $W^{1,4}\subset L^\infty$
to $Dz$ on $Q_{3/2}(s_0)$ now yields
\begin{equation}\label{eq:local-absorption}
 \|Dz\|_{L^\infty(Q_1(s_0))}
 +\|D^2z\|_{L^4(Q_1(s_0))}
 \leq C\|Dz\|_{L^2(Q_4(s_0))}.
\end{equation}

Since the strips $Q_2(s_0-2)$, $Q_2(s_0)$, and $Q_2(s_0+2)$
cover $Q_4(s_0)$ and lie in $\widehat N_{\alpha_j}$, applying \eqref{eq:local-absorption} and Lemma~\ref{lem:neck-decay} on these strips proves the second derivative estimates in \eqref{eq:neck-C1-W24}.
To estimate $z(s,\theta)$, integrate from $(0,0)$ to $(s,0)$
and then along the circle $\{s\}\times S^1$.
Since $z(0,0)=0$, the gradient bound gives
\[
\begin{aligned}
 |z(s,\theta)| \leq C\delta_{\alpha_j}^{3/5}
       \int_0^{|s|}(e^r+e^{-r})^{3/5}\,\dd r
       +C\bigl[\delta_{\alpha_j}(e^s+e^{-s})\bigr]^{3/5} \leq C\bigl[\delta_{\alpha_j}(e^s+e^{-s})\bigr]^{3/5},
\end{aligned}
\]
which completes the proof of Lemma \ref{lem:gnomonic-regularity}.
\end{proof}

We next compare $\alpha_j-1$ with the product of the bubble
scales $\delta_{\alpha_j}$ in \eqref{eq:neck-scales}.
\begin{lemma}\label{lem:neck-Hamiltonian}
For the sequence $u_{\alpha_j}$ given in
Theorem~\ref{thm:two-bubble-index}, the quantity
$\delta_{\alpha_j}$ defined in \eqref{eq:neck-scales} satisfies
\begin{equation}\label{eq:epsilon-delta}
 \alpha_j-1\leq C\delta_{\alpha_j}^{6/5}
\end{equation}
for all sufficiently large $j$, where $C>0$ is independent of $j$.
\end{lemma}
\begin{proof}
Let $S_{\alpha_j}$ denote the stress tensor \eqref{eq:stress}
of $u_{\alpha_j}$ with domain metric $g_{\alpha_j}$. In the
coordinate $t$ introduced in Subsection~\ref{subsec:neck-conservation}, set
\begin{equation}\label{eq:neck-Hamiltonian}
 \mathcal H_{\alpha_j}(t)
 =-\int_{S^1}\left[
       (S_{\alpha_j})_{tt}
       -\frac{\kappa^{\alpha_j}}2(g_{\alpha_j})_{tt}
                    \right]\,\dd\theta.
\end{equation}
The divergence free identity established in
Subsection~\ref{subsec:domain-balance-identity}, expressed in
the conformal coordinates $(t,\theta)$, reads
\[
 \partial_t(S_{\alpha_j})_{tt}
 +\partial_\theta(S_{\alpha_j})_{\theta t}
 =\frac12\partial_t(g_{\alpha_j})_{tt}\,
          \tr_{g_{\alpha_j}}S_{\alpha_j}.
\]
Integrating over $S^1$ and using \eqref{eq:stress-trace} and
\eqref{eq:Omega} therefore gives
\begin{equation}\label{eq:neck-Hamiltonian-derivative}
 \mathcal H_{\alpha_j}'(t)
 =\frac12\int_{S^1}\partial_t(g_{\alpha_j})_{tt}\,
       \Omega_{\alpha_j,\kappa}
            (|du_{\alpha_j}|_{g_{\alpha_j}}^2)\,\dd\theta.
\end{equation}
Smoothness at the 2 domain points corresponding to
$t=\pm\infty$ implies
$\lim_{t\to\pm\infty}\mathcal H_{\alpha_j}(t)=0$.
By \eqref{eq:metric-sign-change} and the nonnegativity of
$\Omega_{\alpha_j,\kappa}$, the function
$\mathcal H_{\alpha_j}$ is nondecreasing for
$t\leq t_{c,\alpha_j}$ and nonincreasing for
$t\geq t_{c,\alpha_j}$, so
\begin{equation}\label{eq:neck-Hamiltonian-positive}
 \mathcal H_{\alpha_j}(t)\geq0\qquad(t\in\R).
\end{equation}

Suppose first that $t_{\alpha_j}^*\leq t_{c,\alpha_j}$.
By \eqref{eq:neck-scales},
\[
 t_{\alpha_j}^*-\log r_{\alpha_j}^1
 =-\log\delta_{\alpha_j}\longrightarrow\infty.
\]
Choose a fixed compact disk in $\mathbb C\setminus\{0\}$
on which $|d\omega_1|$ is bounded below by a positive
constant. The smooth rescaled convergence in
Subsection~\ref{subsec:neck-conservation} shows that the
corresponding disks in the coordinate $\zeta$ lie in
$\{t<t_{\alpha_j}^*\}$ and carry a uniformly positive
Dirichlet energy, while
$|du_{\alpha_j}|_{g_{\alpha_j}}^2\to\infty$ uniformly there.
On these disks, \eqref{eq:round-cylinder-metric} gives
$\frac12\partial_t\log(g_{\alpha_j})_{tt}\to1$,
and the lower bound in \eqref{eq:Omega-two-sided} gives
\[
 \Omega_{\alpha_j,\kappa}
       (|du_{\alpha_j}|_{g_{\alpha_j}}^2)
 \geq\frac14(\alpha_j-1)|du_{\alpha_j}|_{g_{\alpha_j}}^2
\]
for all sufficiently large $j$. Since the integrand in
\eqref{eq:neck-Hamiltonian-derivative} is nonnegative for
$t\leq t_{\alpha_j}^*$, integration from $-\infty$ yields
\begin{equation}\label{eq:neck-Hamiltonian-lower}
\begin{aligned}
 \mathcal H_{\alpha_j}(t_{\alpha_j}^*)
 =\frac12\int_{(-\infty,t_{\alpha_j}^*)\times S^1}
       \partial_t\log(g_{\alpha_j})_{tt}\,
       \Omega_{\alpha_j,\kappa}
       (|du_{\alpha_j}|_{g_{\alpha_j}}^2)\,
                         \dd A_{g_{\alpha_j}}\geq c(\alpha_j-1),
\end{aligned}
\end{equation}
where $c>0$ is independent of $j$.
If $t_{\alpha_j}^*\geq t_{c,\alpha_j}$, integrate instead
from $t_{\alpha_j}^*$ to $+\infty$ and use the second bubble.
Indeed, $-\log r_{\alpha_j}^2-t_{\alpha_j}^*
=-\log\delta_{\alpha_j}\to\infty$, and the compact disks
selected as above for $\omega_2$ satisfy
$\frac12\partial_t\log(g_{\alpha_j})_{tt}\to-1$.
The same argument,
with the sign reversed in the integral, gives the lower
bound in \eqref{eq:neck-Hamiltonian-lower}.

Finally, \eqref{eq:stress} and \eqref{eq:neck-Hamiltonian}
give
\[
 \mathcal H_{\alpha_j}(t_{\alpha_j}^*)
 \leq\alpha_j\int_{S^1}
 (\kappa+|du_{\alpha_j}|_{g_{\alpha_j}}^2)^{\alpha_j-1}
 |(u_{\alpha_j})_t(t_{\alpha_j}^*,\theta)|^2\,\dd\theta,
\]
because $(\kappa+|du_{\alpha_j}|_{g_{\alpha_j}}^2)^{\alpha_j}
-\kappa^{\alpha_j}\geq0$. Using
\eqref{eq:uniform-neck-conductivity},
\eqref{eq:gnomonic-flat-density}, and
\eqref{eq:neck-C1-W24} at $s=0$, we obtain
\[
 \mathcal H_{\alpha_j}(t_{\alpha_j}^*)
 \leq C\int_{S^1}|Dz_{\alpha_j}(0,\theta)|^2\,\dd\theta
 \leq C\delta_{\alpha_j}^{6/5}.
\]
Together with \eqref{eq:neck-Hamiltonian-lower}, this proves
\eqref{eq:epsilon-delta}.
\end{proof}

\subsection{Asymptotic expansion on the neck and rotational fields}\label{subsec:relative-neck-expansion}
In this subsection, we prove the expansion
\eqref{eq:first-neck-expansion} of $z_{\alpha_j}(s,\theta)$ in
Proposition~\ref{prop:first-neck-expansion}, with coefficient
maps $D_\pm$ determined by the expansions of the gnomonic coordinates for two bubbles $\omega_1$ and $\omega_2$ in  \eqref{eq:bubble-expansion-at-infinity}.
The second identity in \eqref{eq:gnomonic-fluxes}, together
with \eqref{eq:neck-Hamiltonian-positive}, then yields
Lemma~\ref{lem:odd-rotational-mode}, which identifies the
infinitesimal rotations will be used to construct the additional
negative variation vector fields in Subsection~\ref{subsec:negative-variations}.

For $Q_1(s)\subset N_{\alpha_j}$, applying
\eqref{eq:neck-C1-W24} to the terms on the right-hand side of
\eqref{eq:gnomonic-nondivergence} gives
\begin{equation}\label{eq:neck-forcing}
\begin{aligned}
 \|\Delta_0z_{\alpha_j}\|_{L^4(Q_1(s))}
 \leq C(\alpha_j-1)
       [\delta_{\alpha_j}(e^s+e^{-s})]^{3/5}
       +C[\delta_{\alpha_j}(e^s+e^{-s})]^{9/5}\leq C[\delta_{\alpha_j}(e^s+e^{-s})]^{9/5}.
\end{aligned}
\end{equation}
Here the last inequality follows
from \eqref{eq:epsilon-delta}.
We next prove the following Poisson estimate needed to approximate
$z_{\alpha_j}$ by a harmonic function, inspired by the approaches
of \cite{Qing-tian,Yin2022}.

\begin{lemma}\label{lem:weighted-cylinder-Poisson}
Let $L>4$, $\delta>0$, and $1<\sigma<2$, and let $f$ be a
scalar or vector valued $L^4$ function on $[-L,L]\times S^1$
satisfying
\[
 \|f\|_{L^4(Q_1(s))}
 \leq C[\delta(e^s+e^{-s})]^\sigma
 \qquad\text{for }|s|\leq L-1,
\]
for some $C>0$. Then there exists
$v\in W^{2,4}_{\mathrm{loc}}((-L,L)\times S^1)$
such that $\Delta_0v=f$  and
\begin{equation}\label{eq:weighted-Poisson-estimate}
 \|v\|_{W^{2,4}(Q_1(s))}
 \leq C_\sigma C[\delta(e^s+e^{-s})]^\sigma
 \qquad\text{for }|s|\leq L-2,
\end{equation}
where $C_\sigma$ depends only on $\sigma$. 
\end{lemma}
\begin{proof}
Extend $f$ by 0 outside the $[-L,L]\times \S^1$ to obtain a function on $\R\times S^1$, which is still denoted by $f$. Near the 2 ends of $[-L,L]\times \S^1$, the
$L^4$ norm of $f$ on a strip in $\R\times S^1$ is controlled by the hypothesis on $Q_1(L-1)$ or $Q_1(-L+1)$. Since
$e^{|s|}\leq e^s+e^{-s}\leq2e^{|s|}$, the extension of $f$ satisfies
\[
 \|f\|_{L^4(Q_1(s))}
 \leq C_\sigma C\delta^\sigma e^{\sigma|s|},
 \qquad \text{for } s\in\R.
\]
Here and below, constants denoted by $C_\sigma$ depend only
on $\sigma$. Write
\[
 f_k(s)=\frac1{2\pi}\int_{S^1}
                  f(s,\theta)e^{-ik\theta}\,\dd\theta,
 \qquad k\in\mathbb Z.
\]
By H\"older's inequality, the same $L^4$ bound as $f$ also holds for $f_k$
on $(s-1,s+1)$ uniformly in $k$.

For $k=0,\pm1$, choose the solutions of
$v_k''-k^2v_k=f_k$ with $v_k(0)=v_k'(0)=0$ as following
\[
 v_0(s):=\int_0^s(s-y)f_0(y)\,\dd y,\qquad
 v_{\pm1}(s):=\int_0^s\sinh(s-y)f_{\pm1}(y)\,\dd y.
\]
Since the convolution kernels in $v_0$ and $v_{\pm1 }$ are bounded by a constant times $e^{|s-y|}$, summing over unit intervals between 0 and $s$  gives
\[
 \sum_{|k|\leq1}|v_k(s)|
 \leq C_\sigma C\delta^\sigma e^{\sigma|s|}.
\]

For $|k|\geq2$, define
\[
 v_k(s):=-\frac1{2|k|}
       \int_{\R}e^{-|k||s-y|}f_k(y)\,\dd y.
\]
Applying H\"older's inequality on the intervals
$[s+\ell,s+\ell+1]$ and $[s-\ell-1,s-\ell]$,
where $\ell=0,1,\ldots$ and using $L^4$ norm bounds for $f_k$, we get
\[
\begin{aligned}
 |v_k(s)| \leq\frac{C_\sigma C\delta^\sigma e^{\sigma|s|}}{|k|}
       \|e^{-|k|y}\|_{L^{4/3}(0,1)}
       \sum_{\ell=0}^{\infty}e^{-(|k|-\sigma)\ell} \leq C_\sigma C\delta^\sigma e^{\sigma|s|}|k|^{-7/4}.
\end{aligned}
\]
Here, the last inequality uses $\sigma<2$ and
$\|e^{-|k|y}\|_{L^{4/3}(0,1)}\leq C_\sigma|k|^{-3/4}$.
Since $\sum_{|k|\geq2}|k|^{-7/4}<\infty$, the Fourier series
\[
 v(s,\theta):=\sum_{k\in\mathbb Z}v_k(s)e^{ik\theta}
\]
converges uniformly on $[-L,L]\times S^1$ and satisfies
$|v(s,\theta)|\leq C_\sigma C\delta^\sigma e^{\sigma|s|}$.
Since each coefficient $v_k$ solves $v_k''-k^2v_k=f_k$, so
$\Delta_0v=f$ in the sense of distributions.

Interior elliptic regularity gives the asserted regularity of
$v$ and, since $Q_2(s)$ lies in the original cylinder whenever
$|s|\leq L-2$, there holds
\[
\begin{aligned}
 \|v\|_{W^{2,4}(Q_1(s))}
 \leq C_\sigma\bigl(\|v\|_{L^4(Q_2(s))}
              +\|f\|_{L^4(Q_2(s))}\bigr)
 \leq C_\sigma C[\delta(e^s+e^{-s})]^\sigma.
\end{aligned}
\]
This proves
\eqref{eq:weighted-Poisson-estimate} and Lemma \ref{lem:weighted-cylinder-Poisson}.
\end{proof}

\begin{proposition}\label{prop:first-neck-expansion}
For the sequence and coordinates introduced in
Subsection~\ref{subsec:neck-conservation}, after further taking 
rotations on target $\S^n$ through angles $O(\delta_{\alpha_j})$ and passage to a subsequence, the gnomonic coordinate functions relative to $q$,
still denoted by $z_{\alpha_j}$, satisfy
\begin{equation}\label{eq:first-neck-expansion}
\begin{gathered}
 \delta_{\alpha_j}^{-1}z_{\alpha_j}(s,\theta)
       \longrightarrow w(s,\theta),\qquad \text{with }\,\,
 w(s,\theta):=e^sD_+\vartheta(\theta)
       +e^{-s}D_-\vartheta(\theta),\quad
 \vartheta(\theta):=
       \begin{pmatrix}\cos\theta\\\sin\theta\end{pmatrix},
\end{gathered}
\end{equation}
in $C^1_{\mathrm{loc}}(\R\times S^1)$, where
$D_\pm:\R^2\to q^\perp$ are linear maps of rank 2.
These rotations on $\S^n$ leave the limiting bubbles unchanged, and
$D_\pm$ are determined by their expansions at the attaching
points. More precisely, since
$\omega_1(\infty)=\omega_2(\infty)=q$, the gnomonic coordinate
functions
\[
 z_{\omega_\ell}(z)
 =\frac{\omega_\ell(z)}{\langle\omega_\ell(z),q\rangle}-q,
 \qquad \ell=1,2,
\]
are defined for $|z|$ sufficiently large and satisfy
\begin{equation}\label{eq:bubble-expansion-at-infinity}
\begin{aligned}
 z_{\omega_1}(Re^{i\theta})
     =R^{-1}D_-\vartheta(\theta)+O(R^{-2}),\qquad
 z_{\omega_2}(Re^{-i\theta}) =R^{-1}D_+\vartheta(\theta)+O(R^{-2})
\end{aligned}
\end{equation}
as $R\to\infty$, with remainders $O(R^{-2})$ independent of $\theta$.
\end{proposition}
\begin{proof}
By \eqref{eq:neck-forcing} and
Lemma~\ref{lem:weighted-cylinder-Poisson} with $\sigma=9/5$,
there exists a solution $v_{\alpha_j}$ of
$\Delta_0v_{\alpha_j}=\Delta_0z_{\alpha_j}$ on $N_{\alpha_j}$
satisfying \eqref{eq:weighted-Poisson-estimate}.
Thus $h_{\alpha_j}:=z_{\alpha_j}-v_{\alpha_j}$ is harmonic,
and utilizing \eqref{eq:gnomonic-smallness} together with the estimate \eqref{eq:weighted-Poisson-estimate} 
and the Sobolev embedding, we can obtain a uniform bound for
$h_{\alpha_j}$ on the circles $s=\pm L$, where
$L=L_{\alpha_j}-3$.

The zero angular Fourier mode of $h_{\alpha_j}$ has the form
$p_{\alpha_j}+q_{\alpha_j}s$, while each first angular Fourier mode
is a linear combination of $e^s$ and $e^{-s}$.
Solving for these coefficients from the values at $s=\pm L$
bounds them by $Ce^{-L}=O(\delta_{\alpha_j})$.
If $h_{\alpha_j,\geq2}$ denotes the summation of the Fourier modes with
$|k|\geq2$, the formulas for the harmonic extension of the
endpoint Fourier coefficients give
\[
 |D^\ell h_{\alpha_j,\geq2}(s,\theta)|
 \leq C\sum_{k=2}^{\infty}k^\ell e^{-k(L-|s|)}
 \leq Ce^{-2(L-|s|)},
 \qquad \ell=0,1,
\]
whenever $|s|\leq L-1$.
By \eqref{eq:neck-cylinder-weight}, this last bound is at most
$C[\delta_{\alpha_j}(e^s+e^{-s})]^{9/5}$, since
$\delta_{\alpha_j}e^{|s|}$ is uniformly bounded on the neck.
Combining the Fourier decomposition with
\eqref{eq:weighted-Poisson-estimate} and the Sobolev embedding
$W^{2,4}\hookrightarrow C^{1,1/2}$, for $|s|\leq L_{\alpha_j}-5$, we obtain
\begin{equation}\label{eq:full-neck-expansion}
\begin{gathered}
 z_{\alpha_j}=p_{\alpha_j}+q_{\alpha_j}s
 +\delta_{\alpha_j}
       (e^sD_{+,\alpha_j}+e^{-s}D_{-,\alpha_j})\vartheta
 +R_{\alpha_j},\quad
 \|R_{\alpha_j}\|_{C^1(Q_1(s))}
 \leq C[\delta_{\alpha_j}(e^s+e^{-s})]^{9/5},
\end{gathered}
\end{equation}
where $p_{\alpha_j},q_{\alpha_j}\in q^\perp$ and
$D_{\pm,\alpha_j}:\R^2\to q^\perp$ are uniformly bounded
linear maps.

At $s=0$, the coefficient in the first identity of
\eqref{eq:gnomonic-fluxes} converges uniformly to 1 by
\eqref{eq:uniform-neck-conductivity} and
\eqref{eq:neck-C1-W24}.
Substituting the derivative of \eqref{eq:full-neck-expansion}
into this identity, and using the zero angular average of
$\vartheta$, gives
\[
 (2\pi+o(1))q_{\alpha_j}
 =o(\delta_{\alpha_j})+O(\delta_{\alpha_j}^{9/5}).
\]
The normalization $z_{\alpha_j}(0,0)=0$ in
\eqref{eq:gnomonic-coordinate} then also controls the constant
term, so
\begin{equation}\label{eq:neck-affine-terms}
 q_{\alpha_j}=o(\delta_{\alpha_j}),\qquad
 p_{\alpha_j}=O(\delta_{\alpha_j}).
\end{equation}

After passing to a subsequence, let $D_\pm$ be the limits of
$D_{\pm,\alpha_j}$.
To identify $D_-$, fix $R>e^5R_0$ and evaluate
\eqref{eq:full-neck-expansion} at
$s=\log(R\delta_{\alpha_j})$, which corresponds to the first
rescaled bubble circle $\zeta=r_{\alpha_j}^1Re^{i\theta}$.
For all sufficiently large $j$, the circle $\log(R\delta_{\alpha_j})\times S^1$ lies in the range
of the remainder estimate of $R_{\alpha_j}$, and the expansion \eqref{eq:full-neck-expansion} gives
\[
\begin{aligned}
 z_{\alpha_j}(\log(R\delta_{\alpha_j}),\theta)
 =R^{-1}D_{-,\alpha_j}\vartheta
   +\delta_{\alpha_j}^2R D_{+,\alpha_j}\vartheta
   +p_{\alpha_j} +q_{\alpha_j}\log(R\delta_{\alpha_j})
   +O\bigl((R^{-1}+R\delta_{\alpha_j}^2)^{9/5}\bigr).
\end{aligned}
\]
At fixed $R$, the affine terms tend to 0 by
\eqref{eq:neck-affine-terms} and
$\delta_{\alpha_j}|\log\delta_{\alpha_j}|\to0$.
The rescaled convergence established in
Subsection~\ref{subsec:neck-conservation} therefore yields
\[
 z_{\omega_1}(Re^{i\theta})
 =R^{-1}D_-\vartheta(\theta)+O(R^{-9/5}).
\]
Since $\omega_1$ is a conformal diffeomorphism onto an
equatorial 2-sphere, its Taylor expansion in the coordinate
$1/z$ at infinity has the form
$z_{\omega_1}(Re^{i\theta})
=R^{-1}B_-\vartheta(\theta)+O(R^{-2})$, with
$\operatorname{rank}B_-=2$.
Comparing these 2 expansions, multiplying by $R$, and then
letting $R\to\infty$ gives $D_-=B_-$.
At the other end, use $s=-\log(R\delta_{\alpha_j})$,
which corresponds to the second rescaled bubble coordinate
$Re^{-i\theta}$, to obtain the same conclusion for $D_+$.
This proves \eqref{eq:bubble-expansion-at-infinity} and the
rank $2$ assertions of $D_{\pm}$.

To remove the constant term $p_{\alpha_j}$, compose
$u_{\alpha_j}$ with the rotation in
$\operatorname{span}\{q,p_{\alpha_j}\}$ that sends
$(q+p_{\alpha_j})/\sqrt{1+|p_{\alpha_j}|^2}$ to $q$ and
fixes the orthogonal complement, choosing the identity when
$p_{\alpha_j}=0$.
Its angle is $\arctan|p_{\alpha_j}|=O(\delta_{\alpha_j})$,
and it sends $q$ to
$(q-p_{\alpha_j})/\sqrt{1+|p_{\alpha_j}|^2}$.
On every fixed compact subcylinder of $\R\times S^1$, \eqref{eq:full-neck-expansion}
and \eqref{eq:neck-affine-terms} give
$z_{\alpha_j},Dz_{\alpha_j}=O(\delta_{\alpha_j})$.
The rotation therefore sends $q+z_{\alpha_j}$ to
$q+z_{\alpha_j}-p_{\alpha_j}+O(\delta_{\alpha_j}^2)$ in
$C^1$, and taking the inner product of the rotated vector $q+z_{\alpha_j}-p_{\alpha_j}+O(\delta_{\alpha_j}^2)$ with \(q\) gives
$1+O(\delta_{\alpha_j}^2)$.
Dividing by this component as in
\eqref{eq:gnomonic-coordinate}, we obtain the gnomonic
coordinates of the rotated map
\[
 \widetilde z_{\alpha_j}
 =z_{\alpha_j}-p_{\alpha_j}+O(\delta_{\alpha_j}^2)
 \qquad\text{in }C^1
\]
on that subcylinder.
Substituting this relation into \eqref{eq:full-neck-expansion}
and using \eqref{eq:neck-affine-terms} gives
$\delta_{\alpha_j}^{-1}\widetilde z_{\alpha_j}\to w$ in
$C^1_{\mathrm{loc}}(\R\times S^1)$.
The rotations converge to the identity, so the limiting
bubbles and their coefficient linear maps remain unchanged.
We henceforth write $z_{\alpha_j}$ for
$\widetilde z_{\alpha_j}$, without imposing the earlier
condition $z_{\alpha_j}(0,0)=0$, which proves
\eqref{eq:first-neck-expansion} and Proposition \ref{prop:first-neck-expansion}.
\end{proof}

The second identity in \eqref{eq:gnomonic-fluxes} and the
nonnegativity of $\mathcal H_{\alpha_j}$ in
\eqref{eq:neck-Hamiltonian-positive} impose the following
constraints on $D_\pm$, which allow us to construct the
rotational fields used in the Morse index estimate in Theorem \ref{thm:two-bubble-index}.
\begin{lemma}\label{lem:odd-rotational-mode}
Let $w$ and $D_\pm$ be as in
Proposition~\ref{prop:first-neck-expansion}. The images of
$D_+$ and $D_-$ coincide, and we write
\[
 P:=\operatorname{range}D_+=\operatorname{range}D_-,
 \qquad \dim P=2.
\]
Define
\begin{equation}\label{eq:neck-normal-space}
 N:=(\R q\oplus P)^\perp,\qquad \dim N=n-2.
\end{equation}
After a rotation or reflection of $\theta$ and a translation
of $s$, both independent of $j$, there is a family of skew
symmetric endomorphisms $K_\xi$ of $\R^{n+1}$, depending
linearly on $\xi\in N$, such that
\begin{equation}\label{eq:odd-rotational-mode}
 K_\xi q=0,\qquad
 K_\xi w(s,\theta)=\sinh s\cos\theta\,\xi
 \qquad\text{for all }(s,\theta)\in\R\times S^1.
\end{equation}
\end{lemma}
\begin{proof}
Dividing the second identity in \eqref{eq:gnomonic-fluxes}
by $\delta_{\alpha_j}^2$ and passing to the limit using
\eqref{eq:uniform-neck-conductivity} and
\eqref{eq:first-neck-expansion}, we obtain
$\int_{S^1}w\wedge w_s\,\dd\theta=0$, which, upon
substitution of the formula for $w$, gives
\begin{equation}\label{eq:tail-torque-matrix}
 D_+D_-^T=D_-D_+^T.
\end{equation}
Since $D_\pm$ have rank 2, their transposes map onto
$\R^2$, so the ranges of the 2 sides are respectively
$\operatorname{range}D_+$ and $\operatorname{range}D_-$.
Their common image is therefore a 2-dimensional subspace
$P\subset q^\perp$, and there is an invertible $2\times2$
matrix $M$ such that $D_-=D_+M$.
Equation~\eqref{eq:tail-torque-matrix} becomes
$D_+(M^T-M)D_+^T=0$. Multiplying by a left inverse of
$D_+$ and its transpose yields
\begin{equation}\label{eq:relative-tail-symmetric}
 M=M^T.
\end{equation}

To determine the sign of an eigenvalue of $M$, we use the
quantity $\mathcal H_{\alpha_j}$ defined in
\eqref{eq:neck-Hamiltonian}.
For $x\geq0$ and $\alpha>1$, the identity
\begin{equation}\label{eq:energy-increment-comparison}
 (\kappa+x)^\alpha-\kappa^\alpha
 -(\kappa+x)^{\alpha-1}x
 =\kappa\bigl[(\kappa+x)^{\alpha-1}-\kappa^{\alpha-1}\bigr]
 \geq0
\end{equation}
and the stress tensor formula \eqref{eq:stress}, together
with \eqref{eq:neck-Hamiltonian-positive}, imply
\[
 0\leq\int_{S^1}
 (\kappa+|du_{\alpha_j}|_{g_{\alpha_j}}^2)^{\alpha_j-1}
 \left[(2\alpha_j-1)|(u_{\alpha_j})_s|^2
       -|(u_{\alpha_j})_\theta|^2\right]\,\dd\theta.
\]
Here, the conformal factor cancels as $(g_{\alpha_j})_{tt}|du_{\alpha_j}|_{g_{\alpha_j}}^2
=|(u_{\alpha_j})_s|^2+|(u_{\alpha_j})_\theta|^2$.
By \eqref{eq:gnomonic-coordinate} and
\eqref{eq:first-neck-expansion},
$\delta_{\alpha_j}^{-1}Du_{\alpha_j}\to Dw$ uniformly
on each fixed circle $s \times S^1$, while
\eqref{eq:uniform-neck-conductivity} gives uniform convergence
of the conductivity to 1. Dividing the preceding inequality by
$2\delta_{\alpha_j}^2$ and passing to the limit therefore yields
\begin{equation}\label{eq:tail-Hamiltonian-sign}
 0\leq\frac12\int_{S^1}
       (|w_s|^2-|w_\theta|^2)\,\dd\theta
 =-2\pi\tr(D_+^TD_-).
\end{equation}
Thus $\tr(D_+^TD_+M)\leq0$.
Since $D_+^TD_+$ is positive definite and $M$ is symmetric
and invertible, $M$ must have a negative eigenvalue.

A rotation or reflection of $\theta$ diagonalizes $M$, so
we may write $M=\operatorname{diag}(\mu_1,\mu_2)$ with
$\mu_1<0$ and $D_+=(d_1,d_2)$.
Set $c=\frac12\log(-\mu_1)$, so that
\[
 e^{s+c}+\mu_1e^{-s-c}=2\sqrt{-\mu_1}\sinh s.
\]
Because $d_1,d_2$ form a basis of $P$, there is a unique
$b\in P$ satisfying
\[
 b\cdot d_1=(2\sqrt{-\mu_1})^{-1},\qquad b\cdot d_2=0.
\]
It follows that $b\cdot w(s+c,\theta)=\sinh s\cos\theta$.
For $\xi\in N$, define
\begin{equation}\label{eq:neck-Killing-generator}
 K_\xi y=(b\cdot y)\xi-(\xi\cdot y)b,
 \qquad y\in\R^{n+1}.
\end{equation}
This formula is linear in $\xi$ and defines a skew symmetric
endomorphism. Since $b\in P$, $\xi\perp\R q\oplus P$,
and $w$ takes values in $P$, it gives
$K_\xi q=0$ and
$K_\xi w(s+c,\theta)=\sinh s\cos\theta\,\xi$.
In the translated coordinate where the limiting function is
$w(s+c,\theta)$, retaining the notation $w$ for this function
proves \eqref{eq:odd-rotational-mode}.
\end{proof}

\subsection{Global negative variations and exact orthogonality}\label{subsec:negative-variations}
In this subsection, we construct global variation vector fields that equal
$-K_\xi u_{\alpha_j}$ on the disk containing the first energy concentration
point and $K_\xi u_{\alpha_j}$ on the disk containing the second energy
concentration point. Across a fixed cylinder in the neck, we join
their boundary values by harmonic extension followed by orthogonal
projection onto $u_{\alpha_j}^{-1}T\Sph^n$.
Lemma~\ref{lem:negative-bridge} shows that, for all sufficiently
large $j$, these vector fields span a space of dimension $n-2$ on which
$D^2E_{\alpha_j,\kappa}(u_{\alpha_j})$ is negative definite.
We then construct a space of dimension $2(n-2)$ on which the
same Hessian is negative definite, using variation vector fields
supported near the 2 concentration points. The orthogonality of
the 2 spaces with respect to the Hessian gives the Morse index
lower bound in Theorem~\ref{thm:two-bubble-index}.

The skew symmetric endomorphisms $K_\xi$, for $\xi\in N$, defined
in \eqref{eq:neck-Killing-generator} satisfy $\langle K_\xi u_{\alpha_j},u_{\alpha_j}\rangle=0$, so $K_\xi u_{\alpha_j}$ is a smooth section of $u_{\alpha_j}^{-1}T\Sph^n$. These sections are the infinitesimal
variations of $u_{\alpha_j}$ under the target rotations
$e^{\tau K_\xi}$ for $-\varepsilon < \tau < \varepsilon$.

\begin{lemma}\label{lem:Killing-kernel}
For every $\xi\in N$, the section $K_\xi u_{\alpha_j}$ is a
Jacobi field along $u_{\alpha_j}$, that is,
\begin{equation}\label{eq:Killing-kernel-pairing}
 D^2E_{\alpha_j,\kappa}(u_{\alpha_j})
 [K_\xi u_{\alpha_j},Y]=0
 \qquad\text{for every }Y\in W^{1,2}(\Sph^2,u_{\alpha_j}^{-1}T\Sph^n).
\end{equation}
In the cylinder coordinates obtained after the rotation or
reflection of $\theta$ and the translation of $s$ in
Lemma~\ref{lem:odd-rotational-mode}, these vector fields satisfy
\begin{equation}\label{eq:rotational-profile-limit}
 \delta_{\alpha_j}^{-1}K_\xi u_{\alpha_j}
 \longrightarrow \sinh s\cos\theta\,\xi
\end{equation}
in $C^1([-T,T]\times S^1)$ for every fixed $T>0$, uniformly
for $\xi$ in bounded subsets of $N$.
\end{lemma}
\begin{proof}
Fix $\xi\in N$ and first take $Y$ to be a smooth variation vector field along $u_{\alpha_j}$. Since $K_\xi$ is skew symmetric,
$e^{\tau K_\xi}$ preserves the round target metric and hence
$E_{\alpha_j,\kappa}(e^{\tau K_\xi}v)
=E_{\alpha_j,\kappa}(v)$ for every smooth map $v$.
Applying this identity to the variation $\exp_{u_{\alpha_j}}(\rho Y)$, the resulting family
$e^{\tau K_\xi}\exp_{u_{\alpha_j}}(\rho Y)$ has initial
variation fields $K_\xi u_{\alpha_j}$ and $Y$.
Since $u_{\alpha_j}$ is $\alpha_j$-harmonic, the rotational invariance of $E_{\alpha_j,\kappa}$ gives
\[
\begin{aligned}
 0=\left.\frac{\partial^2}{\partial\tau\,\partial\rho}
       \right|_{(\tau,\rho)=(0,0)}
 E_{\alpha_j,\kappa}\left(
 e^{\tau K_\xi}\exp_{u_{\alpha_j}}(\rho Y)\right) =D^2E_{\alpha_j,\kappa}(u_{\alpha_j})
   [K_\xi u_{\alpha_j},Y].
\end{aligned}
\]
The continuity of the Hessian on $W^{1,2}$ sections, established
in Subsection~\ref{subsec:normal-splitting}, extends this
identity to every $Y\in W^{1,2}(\Sph^2,u_{\alpha_j}^{-1}T\Sph^n)$, proving \eqref{eq:Killing-kernel-pairing}.

For the convergence in \eqref{eq:rotational-profile-limit}, fix $T>0$, using
\eqref{eq:gnomonic-coordinate} and $K_\xi q=0$, we have 
\[
 \delta_{\alpha_j}^{-1}K_\xi u_{\alpha_j}
 =\frac{K_\xi(\delta_{\alpha_j}^{-1}z_{\alpha_j})}
        {\sqrt{1+|z_{\alpha_j}|^2}}
 \longrightarrow K_\xi w
 =\sinh s\cos\theta\,\xi
 \quad\text{in }C^1([-T,T]\times S^1),
\]
where the convergence follows from
\eqref{eq:first-neck-expansion} and the last equality is
\eqref{eq:odd-rotational-mode}.
The bound $\|K_\xi\|\leq C|\xi|$ from \eqref{eq:neck-Killing-generator} makes the convergence
uniform for $\xi$ in bounded subsets of $N$.
\end{proof}
Fix $T>0$ independently of $j$. Since $L_{\alpha_j}\to\infty$,
the cylinder $[-T,T]\times S^1$ lies in the neck region for
all sufficiently large $j$, also after the fixed coordinate
changes in Lemma~\ref{lem:odd-rotational-mode}.
Using the identification in Subsection~\ref{subsec:neck-conservation},
we regard this cylinder as a closed annulus in the domain sphere.
The closures of the 2 components of its complement are
smooth disks bounded by $s=-T$ and $s=T$, respectively.
We call these the left and right caps, more precisely, the left cap consists
of $\{s\leq-T\}$ together with the first bubble centre
$\zeta=0$, whereas the right cap consists of $\{s\geq T\}$
together with the second bubble centre $\zeta=\infty$.

For each $\xi\in N$, define the variation vector field $V_{\alpha_j,\xi}$ on the 2 caps by
\begin{equation}\label{eq:signed-fields-caps}
\begin{aligned}
 V_{\alpha_j,\xi}:=-K_\xi u_{\alpha_j}\quad
 \text{on the left cap},\qquad
 V_{\alpha_j,\xi}:=K_\xi u_{\alpha_j}\quad
 \text{on the right cap}.
\end{aligned}
\end{equation}
To join the vector fields on left and right caps, let
$H_{\alpha_j,\xi}:[-T,T]\times S^1\to\R^{n+1}$ solve the following
Dirichlet problem
\[
\begin{cases}
 \Delta_0H_{\alpha_j,\xi}=0
 &\text{in }(-T,T)\times S^1,\\
 H_{\alpha_j,\xi}(-T,\theta)
 =-K_\xi u_{\alpha_j}(-T,\theta)
 &\text{for }\theta\in S^1,\\
 H_{\alpha_j,\xi}(T,\theta)
 =K_\xi u_{\alpha_j}(T,\theta)
 &\text{for }\theta\in S^1.
\end{cases}
\]
We then define the variation vector field $V_{\alpha_j,\xi}$ on this cylinder by pointwise orthogonal
projection onto the tangent space of the target $\S^n$
\begin{equation}\label{eq:projected-bridge}
 V_{\alpha_j,\xi}:=P_{u_{\alpha_j}}H_{\alpha_j,\xi}
 \quad\text{on }[-T,T]\times S^1,\qquad
 P_uY:=Y-\langle Y,u\rangle u.
\end{equation}

Since $\langle K_\xi u_{\alpha_j},u_{\alpha_j}\rangle=0$, this projection fixes the boundary values of $H_{\alpha_j,\xi}$, so the vector field in \eqref{eq:projected-bridge} agrees with \eqref{eq:signed-fields-caps} on both boundary circles $-T \times S^1$ and $T\times S^1$. Hence, $V_{\alpha_j,\xi}$ is a global
section $V_{\alpha_j,\xi}\in W^{1,2}(\Sph^2,u_{\alpha_j}^{-1}T\Sph^n)$.
Finally, the linear dependence of $K_\xi$ on $\xi$, uniqueness
of harmonic extension, and linearity of the projection $P_{u_{\alpha_j}}$ imply
that $V_{\alpha_j,\xi}$ depends linearly on $\xi\in N$.

\begin{lemma}\label{lem:negative-bridge}
For all sufficiently large $j$, the space
$\mathcal V_{\alpha_j}=\{V_{\alpha_j,\xi}:\xi\in N\}$ has
dimension $n-2$, and $D^2E_{\alpha_j,\kappa}(u_{\alpha_j})$
is negative definite on $\mathcal V_{\alpha_j}$.
More precisely, there holds
\begin{equation}\label{eq:bridge-negative-limit}
 \frac{D^2E_{\alpha_j,\kappa}(u_{\alpha_j})
       [V_{\alpha_j,\xi},V_{\alpha_j,\eta}]}
      {\alpha_j\delta_{\alpha_j}^2}
 \longrightarrow -2\pi\tanh T\,\langle\xi,\eta\rangle,
 \qquad \xi,\eta\in N.
\end{equation}
\end{lemma}
\begin{proof}
The choice of signs in \eqref{eq:signed-fields-caps} and the convergence
\eqref{eq:rotational-profile-limit} imply that the boundary values of
$H_{\alpha_j,\xi}/\delta_{\alpha_j}$ converge in $C^1$ to
$\sinh T\cos\theta\,\xi$ on both circles $s=\pm T$.
The harmonic extension of these limiting boundary values is
\begin{equation}\label{eq:limiting-harmonic-bridge}
 W_\xi(s,\theta)=\tanh T\,\cosh s\cos\theta\,\xi.
\end{equation}
Continuity of harmonic extension from
$W^{1/2,2}(\{-T,T\}\times S^1)$ to
$W^{1,2}([-T,T]\times S^1)$ therefore gives
$H_{\alpha_j,\xi}/\delta_{\alpha_j}\to W_\xi$ strongly in $W^{1,2}([-T,T]\times S^1)$.
Moreover, \eqref{eq:gnomonic-coordinate} and
\eqref{eq:first-neck-expansion} give
$u_{\alpha_j}=q+O(\delta_{\alpha_j})$ in
$C^1([-T,T]\times S^1)$, so $P_{u_{\alpha_j}}\to P_q$ in $C^1([-T,T]\times S^1)$.
Since $\xi\in N\subset q^\perp$, we have $P_qW_\xi=W_\xi$,
and \eqref{eq:projected-bridge} yields
\begin{equation}\label{eq:bridge-strong-convergence}
 V_{\alpha_j,\xi}/\delta_{\alpha_j}
 \longrightarrow W_\xi
 \quad\text{strongly in }W^{1,2}([-T,T]\times S^1).
\end{equation}
The maximum principle applied to $|H_{\alpha_j,\xi}|^2$ and
$|P_{u_{\alpha_j}}H_{\alpha_j,\xi}|\leq|H_{\alpha_j,\xi}|$ give
an $L^\infty$ bound for $V_{\alpha_j,\xi}$, while
\eqref{eq:bridge-strong-convergence} controls its first derivatives.
Consequently,
\[
 \|V_{\alpha_j,\xi}\|_{L^\infty([-T,T]\times S^1)}
 +\|DV_{\alpha_j,\xi}\|_{L^2([-T,T]\times S^1)}
 \leq C(T)\delta_{\alpha_j}|\xi|.
\]
By \eqref{eq:rotational-profile-limit}, the same bound holds with
$V_{\alpha_j,\xi}$ replaced by $K_\xi u_{\alpha_j}$, and $C(T)$
is independent of $j$ and $\xi$.

Since $K_\xi$ is constant and skew symmetric,
$d(K_\xi u_{\alpha_j})=K_\xi du_{\alpha_j}$ and
$\langle du_{\alpha_j},d(K_\xi u_{\alpha_j})\rangle_{g_{\alpha_j}}=0$ pointwise on the domain. 
Thus the second integral in \eqref{eq:hessian} vanishes identically
for $K_\xi u_{\alpha_j}$.
On each cap, \eqref{eq:signed-fields-caps} also gives
$|V_{\alpha_j,\xi}|^2=|K_\xi u_{\alpha_j}|^2$ and
$|dV_{\alpha_j,\xi}|_{g_{\alpha_j}}^2
 =|d(K_\xi u_{\alpha_j})|_{g_{\alpha_j}}^2$, as well as
$\langle du_{\alpha_j},dV_{\alpha_j,\xi}\rangle_{g_{\alpha_j}}=0$.
Subtracting
$D^2E_{\alpha_j,\kappa}(u_{\alpha_j})
[K_\xi u_{\alpha_j},K_\xi u_{\alpha_j}]=0$ in $D^2E_{\alpha_j,\kappa}(u_{\alpha_j}) [V_{\alpha_j,\xi},V_{\alpha_j,\xi}]$, which follows from
\eqref{eq:Killing-kernel-pairing}, therefore cancels all
contributions from the 2 caps in \eqref{eq:hessian}. To write the remaining integrals in cylinder coordinates, observe that
\[
|dV_{\alpha_j,\xi}|_{g_{\alpha_j}}^2
 =\frac{|(V_{\alpha_j,\xi})_s|^2
             +|(V_{\alpha_j,\xi})_\theta|^2}{(g_{\alpha_j})_{ss}},
 \qquad
 \dd A_{g_{\alpha_j}}=(g_{\alpha_j})_{ss}\,\dd\theta\,\dd s,
\]
and the same derivative formula holds for $u_{\alpha_j}$ and
$K_\xi u_{\alpha_j}$.
Hence \eqref{eq:hessian} gives the  identity
\begin{equation}\label{eq:bridge-hessian-cylinder}
\begin{aligned}
 D^2E_{\alpha_j,\kappa}(u_{\alpha_j})
       [V_{\alpha_j,\xi},V_{\alpha_j,\xi}]
 &=\alpha_j\int_{-T}^{T}\int_{S^1}
       (\kappa+|du_{\alpha_j}|_{g_{\alpha_j}}^2)^{\alpha_j-1}
       \bigl(|DV_{\alpha_j,\xi}|^2
             -|D(K_\xi u_{\alpha_j})|^2\bigr)\,\dd\theta\,\dd s\\
 &\quad-\alpha_j\int_{-T}^{T}\int_{S^1}
       (\kappa+|du_{\alpha_j}|_{g_{\alpha_j}}^2)^{\alpha_j-1}
       |Du_{\alpha_j}|^2
       \bigl(|V_{\alpha_j,\xi}|^2
             -|K_\xi u_{\alpha_j}|^2\bigr)\,\dd\theta\,\dd s\\
 &\quad+2\alpha_j(\alpha_j-1)
       \int_{[-T,T]\times S^1}
       (\kappa+|du_{\alpha_j}|_{g_{\alpha_j}}^2)^{\alpha_j-2}
       \langle du_{\alpha_j},dV_{\alpha_j,\xi}\rangle_{g_{\alpha_j}}^2
       \,\dd A_{g_{\alpha_j}}.
\end{aligned}
\end{equation}

The absolute value of the second term in
\eqref{eq:bridge-hessian-cylinder} is bounded by
\[
\begin{aligned}
 &\alpha_j\int_{-T}^{T}\int_{S^1}
 (\kappa+|du_{\alpha_j}|_{g_{\alpha_j}}^2)^{\alpha_j-1}
 |Du_{\alpha_j}|^2
 \bigl(|V_{\alpha_j,\xi}|^2+|K_\xi u_{\alpha_j}|^2\bigr)
 \,\dd\theta\,\dd s\\
 &\qquad\leq C(T)\delta_{\alpha_j}^2
 \int_{-T}^{T}\int_{S^1}
 \bigl(|V_{\alpha_j,\xi}|^2+|K_\xi u_{\alpha_j}|^2\bigr)
 \,\dd\theta\,\dd s
 \leq C(T)\delta_{\alpha_j}^4|\xi|^2,
\end{aligned}
\]
where we used \eqref{eq:uniform-neck-conductivity},
$|Du_{\alpha_j}|=O(\delta_{\alpha_j})$, and the preceding
$L^\infty$ bounds of $V_{\alpha_j,\xi}$ and $K_\xi u_{\alpha_j}$.
For the third term in
\eqref{eq:bridge-hessian-cylinder}, by Cauchy-Schwarz inequality, we have
\begin{equation}\label{eq:alpha-correction-bridge}
\begin{aligned}
 0&\leq2\alpha_j(\alpha_j-1)
 \int_{[-T,T]\times S^1}
 (\kappa+|du_{\alpha_j}|_{g_{\alpha_j}}^2)^{\alpha_j-2}
 \langle du_{\alpha_j},dV_{\alpha_j,\xi}\rangle_{g_{\alpha_j}}^2
 \,\dd A_{g_{\alpha_j}}\\
 &\leq2\alpha_j(\alpha_j-1)
 \int_{[-T,T]\times S^1}
 (\kappa+|du_{\alpha_j}|_{g_{\alpha_j}}^2)^{\alpha_j-1}
 |dV_{\alpha_j,\xi}|_{g_{\alpha_j}}^2\,\dd A_{g_{\alpha_j}}\\
 &=2\alpha_j(\alpha_j-1)
 \int_{-T}^{T}\int_{S^1}
 (\kappa+|du_{\alpha_j}|_{g_{\alpha_j}}^2)^{\alpha_j-1}
 |DV_{\alpha_j,\xi}|^2\,\dd\theta\,\dd s\\
 &\leq C(T)(\alpha_j-1)\delta_{\alpha_j}^2|\xi|^2
 =o(\delta_{\alpha_j}^2)|\xi|^2.
\end{aligned}
\end{equation}
Finally, replacing the weight by 1 in the first integral of
\eqref{eq:bridge-hessian-cylinder} produces an error, after division
by $\alpha_j$, bounded by
\[
\begin{aligned}
 &\left\|(\kappa+|du_{\alpha_j}|_{g_{\alpha_j}}^2)^{\alpha_j-1}
                  -1\right\|_{L^\infty(\Sph^2)}\int_{-T}^{T}\int_{S^1}
 \bigl(|DV_{\alpha_j,\xi}|^2+|D(K_\xi u_{\alpha_j})|^2\bigr)
 \,\dd\theta\,\dd s
 =o(\delta_{\alpha_j}^2)|\xi|^2,
\end{aligned}
\]
by \eqref{eq:uniform-neck-conductivity} and the $L^2$ bounds of $D V_{\alpha_j,\xi}$ and $D( K_\xi u_{\alpha_j})$ above.
Combining these estimates with \eqref{eq:bridge-hessian-cylinder},
we obtain
\[
\begin{aligned}
 \frac{D^2E_{\alpha_j,\kappa}(u_{\alpha_j})
       [V_{\alpha_j,\xi},V_{\alpha_j,\xi}]}
      {\alpha_j\delta_{\alpha_j}^2}=\int_{-T}^{T}\int_{S^1}
 \left(\left|D\frac{V_{\alpha_j,\xi}}{\delta_{\alpha_j}}\right|^2
       -\left|D\frac{K_\xi u_{\alpha_j}}{\delta_{\alpha_j}}\right|^2\right)
 \,\dd\theta\,\dd s+o(1)|\xi|^2.
\end{aligned}
\]
The convergences \eqref{eq:bridge-strong-convergence} and
\eqref{eq:rotational-profile-limit} allow passage to the limit in
these integrals. It follows that
\[
\begin{aligned}
 \lim_{j\to\infty}
 \frac{D^2E_{\alpha_j,\kappa}(u_{\alpha_j})
       [V_{\alpha_j,\xi},V_{\alpha_j,\xi}]}
      {\alpha_j\delta_{\alpha_j}^2}
 =\pi|\xi|^2(\tanh^2T-1)
              \int_{-T}^{T}\cosh(2s)\,\dd s=-2\pi\tanh T\,|\xi|^2.
\end{aligned}
\]
Polarization gives \eqref{eq:bridge-negative-limit}.
Applying this bilinear convergence to an orthonormal basis of $N$
gives uniform convergence of the corresponding quadratic forms
on the unit sphere of $N$, where the limiting value is
$-2\pi\tanh T<0$. For all sufficiently large $j$, the strict
negativity of these quadratic forms for nonzero $\xi$ implies
that $\xi\mapsto V_{\alpha_j,\xi}$ is injective.
Consequently, $\dim\mathcal V_{\alpha_j}=\dim N=n-2$, and
$D^2E_{\alpha_j,\kappa}(u_{\alpha_j})$ is negative definite on
$\mathcal V_{\alpha_j}$.
\end{proof}

\begin{proof}[Proof of Theorem~\ref{thm:two-bubble-index}]
We first construct negative variations supported near the 2 bubbles.
For each $\ell=1,2$, the equatorial image of $\omega_\ell$ contains
$q$, and its ambient 3-dimensional linear span is generated by $q$
and its tangent plane at $q$.
The rank 2 coefficient linear maps in \eqref{eq:bubble-expansion-at-infinity}
span these tangent planes, which both equal $P$ by
Lemma~\ref{lem:odd-rotational-mode}.
Consequently, both bubbles take values in
$\Sph^n\cap(\R q\oplus P)$, and every constant vector $\nu\in N$
is perpendicular to $\omega_\ell$ and $d\omega_\ell$.

For $R>1$, define the logarithmic cutoff
\[
 \chi_R(z)=
 \begin{cases}
  1,&|z|\leq R,\\
  \dfrac{\log(R^2/|z|)}{\log R},&R<|z|<R^2,\\
  0,&|z|\geq R^2.
 \end{cases}
\]
Utilizing Euclidean norms and area measure on $\mathbb C$, we have
\[
 \int_{\mathbb C}|d\chi_R|^2
 =\frac{2\pi}{(\log R)^2}\int_R^{R^2}\frac{\dd r}{r}
 =\frac{2\pi}{\log R}\longrightarrow0 \qquad \text{as } R \to 
 \infty.
\]
Since $\nu \in N$ is normal to the image of $\omega_\ell$, using
\eqref{eq:hessian} at $\alpha=1$, we have 
\begin{equation}\label{eq:bubble-core-negative}
\begin{aligned}
 D^2E(\omega_\ell)[\chi_R\nu,\chi_R\nu]
 =|\nu|^2\left(\frac{2\pi}{\log R}
       -\int_{\mathbb C}|d\omega_\ell|^2\chi_R^2\right) \longrightarrow-8\pi|\nu|^2
 \qquad\text{as }R\to\infty.
\end{aligned}
\end{equation}
We now fix $R$ sufficiently large such that the Hessian in
\eqref{eq:bubble-core-negative} is at most $-4\pi|\nu|^2$
for both $\ell=1,2$. Using the pullback notation of
Subsection~\ref{subsec:neck-conservation}, define the vector fields on the corresponding rescaled disks by
\[
\begin{aligned}
 Y^1_{\alpha_j,\nu}(r_{\alpha_j}^1z)
 =\chi_R(z)P_{u_{\alpha_j}(r_{\alpha_j}^1z)}\nu,\qquad
 Y^2_{\alpha_j,\nu}((r_{\alpha_j}^2z)^{-1})
 =\chi_R(z)P_{u_{\alpha_j}((r_{\alpha_j}^2z)^{-1})}\nu,
 \qquad |z|<R^2,
\end{aligned}
\]
and extend them by 0 outside these disks. Recall that the second formula is interpreted at $z=0$ in the local coordinate
$1/\zeta$ at $\zeta=\infty$.
These extensions are global
$W^{1,2}$ sections of $u_{\alpha_j}^{-1}T\Sph^n$ and depend
linearly on $\nu$. The smooth rescaled convergence in
Subsection~\ref{subsec:neck-conservation} and
$P_{\omega_\ell}\nu=\nu$ imply that the rescaled vector fields converge strongly in $W^{1,2}(B_{R^2}(0))$ to $\chi_R\nu$,
uniformly for $|\nu|\leq1$. By conformal invariance of Dirichlet energy, we have
\[
 \|Y^\ell_{\alpha_j,\nu}\|_{L^\infty(\Sph^2)}\leq|\nu|,
 \qquad
 \int_{\Sph^2}|dY^\ell_{\alpha_j,\nu}|_{g_{\alpha_j}}^2
       \,\dd A_{g_{\alpha_j}}\leq C(R)|\nu|^2.
\]
By Cauchy-Schwarz inequality used as in \eqref{eq:alpha-correction-bridge}, we get
\[
\begin{aligned}
 &0\leq2\alpha_j(\alpha_j-1)\int_{\Sph^2}
 (\kappa+|du_{\alpha_j}|_{g_{\alpha_j}}^2)^{\alpha_j-2}
 \langle du_{\alpha_j},dY^\ell_{\alpha_j,\nu}\rangle_{g_{\alpha_j}}^2
 \,\dd A_{g_{\alpha_j}} \leq C(R)(\alpha_j-1)|\nu|^2\longrightarrow0\quad \text{as } j \to \infty.
\end{aligned}
\]
On the fixed rescaled disks $B_{R^2}(0)$, the rescaled maps converge smoothly and the vector
fields $Y_{\alpha_j,\nu}^\ell$ converge strongly in $W^{1,2}$, while the conductivity converges
uniformly to 1 by \eqref{eq:uniform-neck-conductivity}. Thus, for the fixed $R$ chosen above, we have 
\[
 D^2E_{\alpha_j,\kappa}(u_{\alpha_j})
       [Y^\ell_{\alpha_j,\nu},Y^\ell_{\alpha_j,\nu}]
 \longrightarrow
 D^2E(\omega_\ell)[\chi_R\nu,\chi_R\nu],
 \qquad \ell=1,2,
\]
uniformly for $|\nu|\leq1$.
For all sufficiently large $j$, each of these quadratic forms is
therefore at most $-2\pi|\nu|^2$.
This also proves that $\nu\mapsto Y^\ell_{\alpha_j,\nu}$ is
injective. Consequently, each space
$\mathcal E^\ell_{\alpha_j}:=\{Y^\ell_{\alpha_j,\nu}:\nu\in N\}$
has dimension $n-2$, and $D^2E_{\alpha_j,\kappa}(u_{\alpha_j})$
is negative definite on $\mathcal E^\ell_{\alpha_j}$ for $\ell=1,2$.

We next locate their supports relative to the fixed cylinder
$[-T,T]\times S^1$ used in Lemma~\ref{lem:negative-bridge}.
In the unshifted coordinate $s=t-t_{\alpha_j}^*$ from
\eqref{eq:neck-scales}, the first rescaled disk satisfies
$s\leq\log(\delta_{\alpha_j}R^2)$ away from its centre,
whereas the second satisfies
$s\geq-\log(\delta_{\alpha_j}R^2)$.
Since $\delta_{\alpha_j}\to0$, these bounds tend to
$-\infty$ and $+\infty$, respectively.
The fixed translation in Lemma~\ref{lem:odd-rotational-mode}
does not change this conclusion, so for large $j$ the supports of $Y_{\alpha_j,\nu}^1$ and $Y_{\alpha_j,\nu}^2$
lie in the interiors of the left and right caps, respectively.
They are disjoint, and by
\eqref{eq:hessian} we have 
\[
 D^2E_{\alpha_j,\kappa}(u_{\alpha_j})
       [Y^1_{\alpha_j,\nu},Y^2_{\alpha_j,\eta}]=0,
 \qquad \text{for any } \nu,\eta\in N.
\]
Consequently, the space
$\mathcal E_{\alpha_j}:=\mathcal E^1_{\alpha_j}\oplus\mathcal E^2_{\alpha_j}$
has dimension $2(n-2)$, and $D^2E_{\alpha_j,\kappa}(u_{\alpha_j})$
is negative definite on $\mathcal E_{\alpha_j}$.

Given $Y\in\mathcal E^\ell_{\alpha_j}$, by
\eqref{eq:signed-fields-caps} we see that  $V_{\alpha_j,\xi}=\pm K_\xi u_{\alpha_j}$ in a neighbourhood
of the support of $Y$, with the minus sign when $\ell=1$
and the plus sign when $\ell=2$.
Polarizing \eqref{eq:hessian} and applying
\eqref{eq:Killing-kernel-pairing}  yields
\begin{equation}\label{eq:exact-core-orthogonality}
\begin{aligned}
 D^2E_{\alpha_j,\kappa}(u_{\alpha_j})[Y,V_{\alpha_j,\xi}]
 &=\pm D^2E_{\alpha_j,\kappa}(u_{\alpha_j})
       [Y,K_\xi u_{\alpha_j}] = 0.
\end{aligned}
\end{equation}
Thus, $\mathcal E_{\alpha_j}$ is orthogonal to $\mathcal V_{\alpha_j}$ with respect to the Hessian $D^2E_{\alpha_j,\kappa}(u_{\alpha_j})$.
Lemma~\ref{lem:negative-bridge} shows that $\mathcal V_{\alpha_j}$
has dimension $n-2$ and that $D^2E_{\alpha_j,\kappa}(u_{\alpha_j})$
is negative definite on $\mathcal V_{\alpha_j}$.
The intersection of $\mathcal E_{\alpha_j}$ and $\mathcal V_{\alpha_j}$ is zero, as a vector field belonging to both
spaces would have zero Hessian value by orthogonality and
would therefore vanish by negative definiteness.
For $Y\in\mathcal E_{\alpha_j}$ and $\xi\in N$, orthogonality of $\mathcal E_{\alpha_j}$ and $\mathcal V_{\alpha_j}$
gives
\[
\begin{aligned}
 D^2E_{\alpha_j,\kappa}(u_{\alpha_j}) [Y+V_{\alpha_j,\xi},Y+V_{\alpha_j,\xi}] =D^2E_{\alpha_j,\kappa}(u_{\alpha_j})[Y,Y]
       +D^2E_{\alpha_j,\kappa}(u_{\alpha_j})
        [V_{\alpha_j,\xi},V_{\alpha_j,\xi}].
\end{aligned}
\]
This value is negative unless both fields vanish, so
$D^2E_{\alpha_j,\kappa}(u_{\alpha_j})$ is negative definite on
$\mathcal E_{\alpha_j}\oplus\mathcal V_{\alpha_j}$.
The dimension of this direct sum is
\[
 \dim(\mathcal E_{\alpha_j}\oplus\mathcal V_{\alpha_j})
 =2(n-2)+(n-2)=3(n-2).
\]
The $W^{1,2}$ definition of the Morse index in
Subsection~\ref{subsec:main-index} now gives
\eqref{eq:two-bubble-index-bound} on the subsequence of $u_{\alpha_j}$ used in
the preceding constructions. If \eqref{eq:two-bubble-index-bound} failed for infinitely many
$j$, we could restrict to that subsequence and then pass to a
further subsequence on which the constructions above apply.
The preceding construction would then give a subspace of dimension
$3(n-2)$ on which the Hessian is negative definite for all
sufficiently large indices in that further subsequence, contradicting
$\ind_{\Sph^n}E_{\alpha_j,\kappa}(u_{\alpha_j})<3(n-2)$ there.
Hence \eqref{eq:two-bubble-index-bound} holds for all sufficiently
large $j$. This completes the proof of Theorem \ref{thm:two-bubble-index}
\end{proof}

\subsection{Degree 2 minimizers and completion of the proof of Theorem \ref{thm:main-full-gap}}\label{subsec:degree-two-sharpness}
In this subsection, we construct minimizers of $E_{\alpha,\kappa}$
among maps $\Sph^2\to\Sph^2$ of degree 2 and prove in
Proposition~\ref{prop:degree-two-minimizers} that their equatorial
inclusions into $\Sph^n$ have Morse index exactly $3(n-2)$ for
$\alpha>1$ sufficiently close to 1. We then complete the proof of
Theorem~\ref{thm:main-full-gap}, by recalling that  Proposition~\ref{prop:two-bubble-reduction}
and Theorem~\ref{thm:two-bubble-index} give the rigidity assertion,
and the degree 2 minimizers obtained in Proposition \ref{prop:degree-two-minimizers} show that the strict Morse index bound is sharp.

\begin{proposition}\label{prop:degree-two-minimizers}
Fix $\kappa>0$. For every $\alpha>1$, the functional
$E_{\alpha,\kappa}$ attains its minimum among maps in
$W^{1,2\alpha}(\Sph^2,\Sph^2)$ of degree 2, and every minimizer
$v_\alpha$ is smooth and $\alpha$-harmonic. There exists
$\alpha_0=\alpha_0(\kappa)>1$ such that every such minimizer satisfies
\begin{equation}\label{eq:degree-two-spectral-data}
 \ind_{\Sph^2}E_{\alpha,\kappa}(v_\alpha)=0,
 \qquad I_\alpha(v_\alpha)=3,
 \qquad N_\alpha(v_\alpha)=3
\end{equation}
whenever $1<\alpha<\alpha_0$. Moreover, we have 
\begin{equation}\label{eq:degree-two-energy-limit}
 E(v_\alpha)>8\pi,
 \qquad
 E_{\alpha,\kappa}(v_\alpha)\longrightarrow2\pi\kappa+8\pi
 \quad\text{as }\alpha\searrow1,
\end{equation}
where the strict inequality holds for every $\alpha>1$ and the
limit holds for every choice of the minimizers.

For every equatorial inclusion $\iota:\Sph^2\hookrightarrow\Sph^n$,
with $n\geq3$, the maps $u_\alpha=\iota\circ v_\alpha$ are smooth
$\alpha$-harmonic maps for all $\alpha>1$ and satisfy
\begin{equation}\label{eq:degree-two-full-index}
 \ind_{\Sph^n}E_{\alpha,\kappa}(u_\alpha)=3(n-2)
 \qquad\text{for }1<\alpha<\alpha_0.
\end{equation}
In particular, none of these maps is an isometric equatorial parametrization.
\end{proposition}
\begin{proof}
For fixed $\alpha>1$, a minimizing sequence in the degree 2 class
of $W^{1,2\alpha}(\Sph^2,\Sph^2)$ is bounded in $W^{1,2\alpha}$
and therefore has a subsequence converging weakly in this space
and uniformly to a sphere-valued map $v_\alpha$. Uniform
convergence preserves degree, since sufficiently close maps into
$\Sph^2$ are homotopic through the shortest geodesics joining
their values. Weak lower semicontinuity of $E_{\alpha,\kappa}$
then shows that $v_\alpha$ is a minimizer, and hence weakly
$\alpha$-harmonic. By \cite[Proposition~2.3]{SacksUhlenbeck1981},
with the normalization \eqref{eq:kappa-normalization}, it is
smooth. For every smooth tangent field $Y$ along $v_\alpha$,
the variation $\exp_{v_\alpha}(tY)$ preserves degree, so
$D^2E_{\alpha,\kappa}(v_\alpha)[Y,Y]\geq0$. By density and
continuity of the Hessian, the same inequality holds for
$W^{1,2}$ tangent fields, giving
\begin{equation}\label{eq:degree-two-intrinsic-stability}
 \ind_{\Sph^2}E_{\alpha,\kappa}(v_\alpha)=0
 \qquad\text{for every }\alpha>1.
\end{equation}

To compare the energies as $\alpha\searrow1$, fix a holomorphic
map $w:\Sph^2\to\Sph^2$ of degree 2. Since $E(w)=8\pi$,
minimality of $v_\alpha$ and smoothness of $w$ give
\begin{equation}\label{eq:degree-two-comparison}
 E_{\alpha,\kappa}(v_\alpha)
 \leq E_{\alpha,\kappa}(w)
 =2\pi\kappa+8\pi+O(\alpha-1).
\end{equation}
The estimates below are uniform over all degree 2 minimizers
for $\alpha$ sufficiently close to 1. For every $t>0$, convexity
of the exponential function yields
\begin{equation}\label{eq:degree-two-convexity}
 t^\alpha
 =t\exp\bigl((\alpha-1)\log t\bigr)
 \geq t+(\alpha-1)t\log t
 \geq t-\frac{\alpha-1}{e}.
\end{equation}
Applying \eqref{eq:degree-two-convexity} with
$t=\kappa+|dv_\alpha|_g^2$, and using
\eqref{eq:degree-two-comparison}, $\operatorname{Area}_g(\Sph^2)=4\pi$,
and the degree-energy inequality $E(v_\alpha)\geq8\pi$, we obtain
\begin{equation}\label{eq:degree-two-energy-bounds}
\begin{aligned}
 2\pi\kappa+8\pi-\frac{2\pi}{e}(\alpha-1)
 &\leq E_{\alpha,\kappa}(v_\alpha)
 \leq2\pi\kappa+8\pi+C(\alpha-1),\\
 8\pi&\leq E(v_\alpha)\leq8\pi+C(\alpha-1).
\end{aligned}
\end{equation}
In particular, \eqref{eq:degree-two-energy-bounds} proves the
energy limit in \eqref{eq:degree-two-energy-limit} and gives
$E(v_\alpha)\to8\pi$.

The logarithmic term in \eqref{eq:degree-two-convexity}, together
with \eqref{eq:degree-two-comparison} and $E(v_\alpha)\geq8\pi$,
also gives
\begin{equation}\label{eq:degree-two-entropy-comparison}
\begin{aligned}
 \frac{\alpha-1}{2}\int_{\Sph^2}
 (\kappa+|dv_\alpha|_g^2)
 \log(\kappa+|dv_\alpha|_g^2)\,\dd A \leq E_{\alpha,\kappa}(v_\alpha)
       -2\pi\kappa-E(v_\alpha)
 \leq C(\alpha-1).
\end{aligned}
\end{equation}
Dividing \eqref{eq:degree-two-entropy-comparison} by $\alpha-1$
yields
\begin{equation}\label{eq:degree-two-entropy}
 \int_{\Sph^2}(\kappa+|dv_\alpha|_g^2)
 \log(\kappa+|dv_\alpha|_g^2)\,\dd A\leq C.
\end{equation}
Since $t\log t\geq-e^{-1}$, the negative part of the integrand
has integral at most $4\pi/e$, so \eqref{eq:degree-two-entropy}
also bounds the integral of its positive part. For $K>1$, the
integrand is at least $|dv_\alpha|_g^2\log K$ on
$\{|dv_\alpha|_g^2>K\}$, and hence
\begin{equation}\label{eq:degree-two-large-gradient-energy}
 \int_{\{|dv_\alpha|_g^2>K\}}|dv_\alpha|_g^2\,\dd A
 \leq\frac{C}{\log K}.
\end{equation}
Splitting any measurable set $A\subset\Sph^2$ into the regions
where $|dv_\alpha|_g^2\leq K$ and $|dv_\alpha|_g^2>K$, and
using \eqref{eq:degree-two-large-gradient-energy}, we obtain
\begin{equation}\label{eq:degree-two-set-energy}
 \int_A|dv_\alpha|_g^2\,\dd A
 \leq K|A|+\frac{C}{\log K},
\end{equation}
where $|A|$ denotes the area of $A$. Taking $K$ large and then
$|A|$ small in \eqref{eq:degree-two-set-energy} proves uniform
integrability of the Dirichlet energy densities.

Now let $\alpha_j\searrow1$ and choose an arbitrary degree 2
minimizer $v_{\alpha_j}$ for each $j$. The bound
\eqref{eq:degree-two-comparison} allows us to apply
Theorem~\ref{thm:compactness-input}. Every nonconstant harmonic
bubble has Dirichlet energy at least $4\pi$, whereas the uniform
integrability established by \eqref{eq:degree-two-set-energy}
makes the Dirichlet energy on sufficiently small domain balls
arbitrarily small, uniformly in $j$. Thus the concentration set $\mathcal{S}$
is empty, and there is a subsequence of $v_{\alpha_j}$ converges smoothly on 
$\Sph^2$ to a harmonic map $v:\Sph^2\to\Sph^2$ of $\deg v=2$, while
\eqref{eq:degree-two-energy-bounds} gives $E(v)=8\pi$.

Lemma~\ref{lem:harmonic-scalar-gap} therefore gives
$I_1(v)=N_1(v)=3$. By the smooth convergence of $v_{\alpha_j}$, we see that the potentials in \eqref{eq:conjugation} converge smoothly
to $-|dv|_g^2$, so the min-max principle on the fixed space
$L^2(\Sph^2,\dd A)$ gives convergence of the eigenvalues of the conjugated scalar Schr\"odinger operators, counted with multiplicity. Since the limiting operator of \eqref{eq:conjugation} has 3 negative eigenvalues and a
3-dimensional kernel, for all sufficiently large $j$ we have
\begin{equation}\label{eq:degree-two-spectral-bounds}
 I_{\alpha_j}(v_{\alpha_j})\geq3,
 \qquad
 I_{\alpha_j}(v_{\alpha_j})+
 N_{\alpha_j}(v_{\alpha_j})\leq6.
\end{equation}
By \eqref{eq:sphere}, the 3 coordinate functions of
$v_{\alpha_j}$ belong to $\ker J_{\alpha_j,v_{\alpha_j}}$. Surjectivity of a map of degree 2 makes them linearly independent, so $N_{\alpha_j}(v_{\alpha_j})\geq3$. Together with
\eqref{eq:degree-two-spectral-bounds}, this gives
$I_{\alpha_j}(v_{\alpha_j})=N_{\alpha_j}(v_{\alpha_j})=3$.
If no common $\alpha_0(\kappa)>1$ existed for these identities,
we could choose a sequence of minimizers violating them with
$\alpha_j\searrow1$, contradicting the conclusion just obtained.
Combining these scalar identities with
\eqref{eq:degree-two-intrinsic-stability} proves
\eqref{eq:degree-two-spectral-data}.

The coordinate equation \eqref{eq:sphere} is preserved under
equatorial inclusion, so $u_\alpha=\iota\circ v_\alpha$ is
$\alpha$-harmonic for every $\alpha>1$. The normal splitting
formula \eqref{eq:index-splitting} and
\eqref{eq:degree-two-spectral-data} now give
\[
\begin{aligned}
 \ind_{\Sph^n}E_{\alpha,\kappa}(\iota\circ v_\alpha)
 =\ind_{\Sph^2}E_{\alpha,\kappa}(v_\alpha)
      +(n-2)I_\alpha(v_\alpha) =3(n-2)
 \qquad\text{for }1<\alpha<\alpha_0,
\end{aligned}
\]
which proves \eqref{eq:degree-two-full-index}.

Finally, if $E(v_\alpha)=8\pi$ for some $\alpha>1$, equality
in the degree and energy inequality makes $v_\alpha$ holomorphic of degree 2. Its critical points are isolated branch points and
hence form a finite set. On the connected complement of this
set, $v_\alpha$ is harmonic and $dv_\alpha$ is invertible, so
\eqref{eq:drift} gives $\nabla|dv_\alpha|_g^2=0$.
The energy density is therefore a positive constant there and,
by continuity, on all of $\Sph^2$, which rules out every branch
point. Thus $v_\alpha$ is an unbranched covering of $\Sph^2$
by itself. Since the target is simply connected, this covering
is a diffeomorphism, contradicting $\deg v_\alpha=2$.
Consequently, $E(v_\alpha)>8\pi$ for every $\alpha>1$, while
degree 2 also excludes an isometric parametrization for
$u_\alpha=\iota\circ v_\alpha$. This completes the proof of Proposition \ref{prop:degree-two-minimizers}.
\end{proof}

\begin{remark}\label{rem:energy-degree-not-alpha-degree}
For the degree 2 minimizers $v_{\alpha}$ obtained in Proposition \ref{prop:degree-two-minimizers},
\eqref{eq:degree-two-spectral-data} gives
$I_\alpha(v_\alpha)+N_\alpha(v_\alpha)=6$, whereas
$E(v_\alpha)>8\pi$, which implies
\[
 I_\alpha(v_\alpha)+N_\alpha(v_\alpha)
 <\frac{E(v_\alpha)}{2\pi}+2.
\]
Thus the lower bound for the sum of scalar index and nullity in
\eqref{eq:harmonic-degree-estimates} does not extend to
$\alpha$-harmonic maps by substituting $E(u)/(4\pi)$ for the
harmonic degree. In the proofs of
Propositions~\ref{prop:two-bubble-reduction}
and~\ref{prop:degree-two-minimizers}, these harmonic estimates
are used through Lemma~\ref{lem:harmonic-scalar-gap} only for
the harmonic components obtained in the limit.
\end{remark}

\begin{proof}[Proof of Theorem~\ref{thm:main-full-gap}]
Fix $n\geq3$, $\kappa>0$, and $\Lambda<\infty$, and suppose
that the rigidity assertion in Theorem \ref{thm:main-full-gap} fails for exponents $\alpha$ arbitrarily close to 1. We can then choose a sequence $\alpha_j\searrow1$ and
nonconstant smooth $\alpha_j$-harmonic maps
$u_{\alpha_j}:\Sph^2\to\Sph^n$ satisfying
\eqref{eq:two-bubble-budget}, none of which is an isometric
equatorial parametrization. Proposition~\ref{prop:two-bubble-reduction}
therefore gives, after passing to a subsequence, a constant base
map and exactly 2 equatorial bubbles of harmonic degree 1 at
distinct antipodal points. Theorem~\ref{thm:two-bubble-index}
then gives the Morse index lower bound \eqref{eq:two-bubble-index-bound}
for all sufficiently large $j$, contradicting the strict Morse index upper
bound in \eqref{eq:two-bubble-budget}. This proves the rigidity
assertion part of Theorem \ref{thm:main-full-gap} for some $\alpha_0(n,\kappa,\Lambda)>1$.
Corollary~\ref{cor:equatorial-index} gives the Morse index $n-2$
and nullity $3n-3$ of the isometric equatorial parametrizations.

For sharpness part of Theorem \ref{thm:main-full-gap}, let $v_\alpha$ be degree 2 minimizers provided by
Proposition~\ref{prop:degree-two-minimizers}, and fix an equatorial
inclusion $\iota:\Sph^2\hookrightarrow\Sph^n$. The maps
$u_\alpha=\iota\circ v_\alpha$ are nonisometric
$\alpha$-harmonic maps and have Morse index exactly $3(n-2)$
for $\alpha>1$ sufficiently close to 1, by
\eqref{eq:degree-two-full-index}. Equatorial inclusion preserves
the energy, so \eqref{eq:degree-two-energy-limit} also gives
the energy limit in \eqref{eq:main-sharpness}. For every fixed
$\Lambda>2\pi\kappa+8\pi$, this convergence implies
$E_{\alpha,\kappa}(u_\alpha)<\Lambda$ for all $\alpha>1$
sufficiently close to 1. These nonisometric maps therefore
satisfy the energy upper bound and attain the Morse index
threshold, showing that the strict inequality in \eqref{eq:main-sharp-hyp} cannot be replaced by a nonstrict one. This completes the proof of Theorem \ref{thm:main-full-gap}.
\end{proof}

\appendix
\vskip1cm

\section{Pole expansions for the equivariant equation}\label{app:auxiliary}\label{app:poles}

In this section, we derive the pole expansions
\eqref{eq:north-pole} and~\eqref{eq:south-pole} that determine
the endpoint signs in the proof of Proposition~\ref{prop:equivariant}.
Throughout this appendix, $\alpha>1$ and $\kappa>0$, and all norms
and measures use the unit round domain metric $g$ unless otherwise
indicated.

As shown in the proof of Proposition~\ref{prop:equivariant}, smoothness at the poles and the second equation in \eqref{eq:twist-equations} imply $\gamma'=0$, so a target rotation removes the constant twist in \eqref{eq:twisted-map}. The first equation in \eqref{eq:twist-equations}, with $\gamma'=0$ and $|du|_g^2=(\partial_r f)^2+\sin^2f/\sin^2r$, therefore gives the profile equation
\begin{equation*}
\begin{aligned}
 &\partial_r\left[
 \left(\kappa+(\partial_r f)^2+\frac{\sin^2f}{\sin^2r}\right)^{\alpha-1}
 \partial_r f\sin r\right] = \left(\kappa+(\partial_r f)^2+\frac{\sin^2f}{\sin^2r}\right)^{\alpha-1}
 \frac{\sin f\cos f}{\sin r}.
\end{aligned}
\end{equation*}
Write $f(r)=cr+br^3+O(r^5)$ at the north pole. Smooth rotational equivariance gives an odd expansion, and diffeomorphicity gives $c>0$. We first expand the factors to get
\begin{align*}
 \sin f&=cr+\left(b-\frac{c^3}{6}\right)r^3+O(r^5),\\
 \cos f&=1-\frac{c^2}{2}r^2+O(r^4),\qquad
 \sin r=r-\frac{r^3}{6}+O(r^5),\\
 \frac{\sin^2f}{\sin^2r}
 &=\frac{c^2r^2+2c(b-c^3/6)r^4+O(r^6)}
 {r^2-r^4/3+O(r^6)}=c^2+\left(2bc+\frac{c^2-c^4}{3}\right)r^2+O(r^4).
\end{align*}
Since $(\partial_r f)^2=c^2+6bcr^2+O(r^4)$, we obtain
\begin{align*}
 (\partial_r f)^2+\frac{\sin^2f}{\sin^2r}
 &=2c^2+\left(8cb+\frac{c^2(1-c^2)}{3}\right)r^2+O(r^4),\\
 \left(\kappa+(\partial_r f)^2+\frac{\sin^2f}{\sin^2r}\right)^{\alpha-1}
 &=(\kappa+2c^2)^{\alpha-1}
 \left[1+\frac{(\alpha-1)\left(8cb+\frac{c^2(1-c^2)}{3}\right)}
 {\kappa+2c^2}r^2+O(r^4)\right].
\end{align*}
For the left-hand side of the profile equation, $\partial_r f\sin r=cr+(3b-c/6)r^3+O(r^5)$ gives
\begin{align*}
 &\partial_r\left[
 \left(\kappa+(\partial_r f)^2+\frac{\sin^2f}{\sin^2r}\right)^{\alpha-1}
 \partial_r f\sin r\right]\\*
 &\quad=(\kappa+2c^2)^{\alpha-1}
 \left[c+\left(9b-\frac c2+
 \frac{3(\alpha-1)c\left(8cb+\frac{c^2(1-c^2)}3\right)}
 {\kappa+2c^2}\right)r^2+O(r^4)\right].
\end{align*}
For the right-hand side of the profile equation, we have
\[
 \frac{\sin f\cos f}{\sin r}
 =c+\left(b-\frac{2c^3}{3}+\frac c6\right)r^2+O(r^4).
\]
Multiplication by the expanded weight and comparison of the coefficients of $r^2$ give
\begin{align*}
 9b-\frac c2+
 \frac{3(\alpha-1)c\left(8cb+\frac{c^2(1-c^2)}3\right)}{\kappa+2c^2}=b-\frac{2c^3}{3}+\frac c6
 +\frac{(\alpha-1)c\left(8cb+\frac{c^2(1-c^2)}3\right)}{\kappa+2c^2},
\end{align*}
which is equivalent to
\[
 8b+\frac23c(c^2-1)
 +\frac{2c(\alpha-1)\left(8cb+\frac{c^2(1-c^2)}3\right)}{\kappa+2c^2}=0.
\]
Collecting the terms containing $b$, we obtain
\begin{align*}
 b\left(1+\frac{2(\alpha-1)c^2}{\kappa+2c^2}\right)
 =\frac{c(1-c^2)}{12}
 -\frac{(\alpha-1)c^3(1-c^2)}{12(\kappa+2c^2)}=\frac{c(1-c^2)[\kappa+(3-\alpha)c^2]}{12(\kappa+2c^2)}.
\end{align*}
Using $1+2(\alpha-1)c^2/(\kappa+2c^2)=(\kappa+2\alpha c^2)/(\kappa+2c^2)$, we obtain
\[
 b=\frac{c(1-c^2)[\kappa+(3-\alpha)c^2]}{12(\kappa+2\alpha c^2)}.
\]
Finally, with $w=\log\bigl(\tan\frac f2/\tan\frac r2\bigr)$ we expand
\[
 \tan\frac r2=\frac r2+\frac{r^3}{24}+O(r^5),\qquad
 \tan\frac f2=\frac{cr}{2}+\Bigl(\frac b2+\frac{c^3}{24}\Bigr)r^3+O(r^5),
\]
which gives
\[
 \frac{\tan\frac f2}{\tan\frac r2}
 =c\Bigl[1+\Bigl(\frac bc+\frac{c^2-1}{12}\Bigr)r^2+O(r^4)\Bigr].
\]
Inserting the value of $b$,
\begin{align*}
 \frac bc+\frac{c^2-1}{12}
 =\frac{(1-c^2)[\kappa+(3-\alpha)c^2]}{12(\kappa+2\alpha c^2)}
 +\frac{c^2-1}{12} =\frac{(\alpha-1)c^2(c^2-1)}{4(\kappa+2\alpha c^2)},
\end{align*}
and therefore
\begin{equation}\label{eq:north-pole}
 w=\log c+\frac{(\alpha-1) c^2(c^2-1)}{4(\kappa+2\alpha c^2)}r^2+O(r^4).
\end{equation}
Applying the same calculation to $\widetilde f(s)=\pi-f(\pi-s)$, with initial slope $d>0$, yields
\begin{equation}\label{eq:south-pole}
 w=-\log d-\frac{(\alpha-1) d^2(d^2-1)}{4(\kappa+2\alpha d^2)}(\pi-r)^2
       +O((\pi-r)^4).
\end{equation}
Because $dr/dt=\sin r>0$, equations \eqref{eq:north-pole} and \eqref{eq:south-pole} give exactly the endpoint signs used in the proof of Proposition~\ref{prop:equivariant}.

\providecommand{\MR}{\relax\ifhmode\unskip\space\fi MR }
\providecommand{\MRhref}[2]{\href{http://www.ams.org/mathscinet-getitem?mr=#1}{#2}}
\bibliographystyle{amsalpha}
\bibliography{references}
\end{document}